%% file: phylo.tex
\documentclass{article}
\usepackage{a4wide}
\usepackage{graphicx}
\usepackage{amsmath}
\usepackage{amssymb}
\usepackage{amsthm,cite}
\usepackage{color}

\input local-defs.tex

\begin{document}

\title{The Reducts of the Homogeneous\\ Binary Branching C-relation\thanks{The first   and third author have received funding from the European Research Council under the European Community's Seventh Framework Programme (FP7/2007-2013 Grant Agreement no. 257039). The third author has received funding from the project P27600 of the Austrian Science Fund (FWF) and the project ``Models on graphs: enumerative combinatorics and algebraic structures'' of the Vietnam National Foundation for Science and Technology Development (NAFOSTED). The second author is partially supported by the {\em Swedish Research Council} (VR) under Grant 621-2012-3239.}
}
\author{Manuel Bodirsky \\
Institut f\"ur Algebra, TU Dresden, Germany \\
{\tt manuel.bodirsky@tu-dresden.de}\\
\and Peter Jonsson\\
Department of Computer and System Science\\ Link\"{o}pings Universitet, Link\"{o}ping, Sweden\\
{\tt peter.jonsson@liu.se}\\
 \and Trung Van Pham\\
Institut f\"ur Computersprachen, Theory and Logic Group\\ TU Wien, Austria\\
{\tt pvtrung@logic.at}}
\maketitle

\begin{abstract}
Let $(\mL;C)$ be the (up to isomorphism unique) countable homogeneous 
structure carrying a binary branching C-relation.
We study the reducts of $(\mL;C)$, i.e., 
the structures with domain $\mL$
that are first-order definable in $(\mL;C)$. We show that up to existential interdefinability, there are finitely many
such reducts.
This implies that there are finitely many reducts up to 
first-order interdefinability, thus confirming a conjecture of Simon Thomas 
for the special case of $(\mL;C)$. 
We also study the endomorphism monoids of such reducts
and show that they fall into four categories.

\end{abstract}
\tableofcontents
\input intro.tex
\input phylogeny.tex
\input autos.tex

\input ramsey.tex
\input canonical.tex
\input two-constants.tex
\input classification.tex

\input ending.tex


\bibliographystyle{abbrv} 
\bibliography{local}


\end{document}

%% file: local-defs.tex
 \usepackage{epsf}
\usepackage{epsfig}
\usepackage{graphicx}
\usepackage{url}
\usepackage{cite}
\usepackage{amsxtra, amssymb}
\usepackage{eucal, url}
\usepackage{xspace, varioref}
\usepackage[matrix, curve, arrow, tips, frame]{xy}
\usepackage{mdwtab}
\usepackage{dcolumn}

\DeclareMathOperator{\Aut}{Aut}

\DeclareMathOperator{\End}{End}

\DeclareMathOperator{\Csp}{CSP}

\DeclareMathOperator{\id}{id}

\DeclareMathOperator{\rer}{rer}

\DeclareMathOperator{\cut}{cut}

\DeclareMathOperator{\lin}{lin}

\DeclareMathOperator{\lrl}{nil}
\DeclareMathOperator{\nil}{nil}
\DeclareMathOperator{\Nil}{Nil}
\DeclareMathOperator{\yca}{yca}

\newcommand\set[1]{\{#1\}}

\newtheorem{theorem}{Theorem}
\newtheorem{proposition}[theorem]{Proposition}
\newtheorem{definition}[theorem]{Definition}

\newtheorem{corollary}[theorem]{Corollary}
\newtheorem{lemma}[theorem]{Lemma}

\theoremstyle{remark}
\newtheorem{example}{Example}

\newcommand{\mL}{{\mathbb L}}
\newcommand{\mQ}{{\mathbb Q}}

\renewcommand{\L}{\mathbb{L}}

\newcommand{\Fresse}{Fra\"{i}ss\'{e}}
\newcommand{\Nesetril}{Ne\v{s}et\v{r}il}

\newcommand{\ignore}[1]{}

\def\Com{{\rm pref}}

%% file: intro.tex
\section{Introduction}
A structure $\Gamma$ is called \emph{homogeneous} (or sometimes \emph{ultra-homogeneous} in order to 
distinguish it from other notions of homogeneity that are used in adjacent areas of mathematics)
if every isomorphism between finite substructures of $\Gamma$ can be extended to an automorphism of $\Gamma$. Many classical structures in mathematics are homogeneous 
such as
$({\mathbb Q};<)$, the random graph, and  
the homogeneous universal poset.

\emph{C-relations} are central for the structure
 theory of Jordan permutation groups~\cite{AdelekeNeumann,JordanSurvey,AdelekeMacPherson,AdelekeNeumannPrimitive}. They also appear frequently
 in model theory. For instance, there is a substantial literature
 on \emph{C-minimal structures} which are analogous
 to $o$-minimal structures but where a C-relation plays the role of the order in an $o$-minimal structure~\cite{C-minimal,MacphersonSteinhorn}.
In this article we study the \emph{universal homogeneous binary branching} $C$-relation $(\mL;C)$. 
This structure is one of the  
 fundamental homogeneous structures~\cite{Oligo,AdelekeNeumann,MacphersonSurvey} and 
can be defined in several different 
ways---we present two distinct definitions in Section~\ref{sect:c-rel}. 
We mention that $(\mL;C)$ is the up to
isomorphism unique countable $C$-relation
which is \emph{existential positive complete} in the class of all $C$-relations -- see~\cite{BodHilsMartin} for the notion of existential positive completeness.


If $\Gamma$ has a finite relational signature (as in the examples mentioned above), then homogeneity implies 
that $\Gamma$ is \emph{$\omega$-categorical}, that is, every countable model of the first-order theory of $\Gamma$ 
is isomorphic to $\Gamma$. 
A relational structure $\Delta$ is called a {\em reduct} of $\Gamma$ if
$\Delta$ and $\Gamma$ have the same domain and
every relation in $\Delta$ has a first-order definition (without parameters) in $\Gamma$. 
It is well known that reducts of $\omega$-categorical structures 
are again $\omega$-categorical~\cite{Hodges}. 
Two reducts $\Delta_1$ and $\Delta_2$ are said to be \emph{first-order interdefinable} if $\Delta_1$ is first-order definable in $\Delta_2$, and vice versa.
\emph{Existential} and \emph{existential positive}\footnote{A first-order formula is \emph{existential} if it is of the form $\exists x_1,\ldots,x_m \: . \: \psi$ where $\psi$
is quantifier-free, and \emph{existential-positive} if it is existential and does not contain the negation symbol $\neg$.}
interdefinability are defined analogously. 

It turns out that several fundamental homogeneous
structures with finite relational signatures
have only \emph{finitely many reducts up to first-order interdefinability}.
This was shown for $(\mathbb Q;<)$ by Cameron~\cite{Cameron5} (and, independently and in somewhat different language, by Frasnay~\cite{Frasnay}), by Thomas for
the the random graph~\cite{RandomReducts},
 by Junker and Ziegler for the expansion of $(\mathbb Q;<)$ by a constant~\cite{JunkerZiegler}, by Pach, Pinsker, Pluh\'ar, Pongr\'acz, and Szab\'o for the homogeneous universal poset~\cite{Poset-Reducts}, and by Bodirsky, Pinsker and Pongr\'acz for the random ordered graph \cite{42}. Thomas has conjectured that \emph{all}
homogeneous structures with a finite relational
signature have finitely many reducts~\cite{RandomReducts}. In this paper, we study the reducts of $(\mL;C)$ up to first-order, and even up to existential and existential positive interdefinability. Our results for reducts up to first-order interdefinability confirm Thomas' conjecture for the case of $(\mL;C)$.    

Studying reducts of $\omega$-categorical structures 
has an additional motivation coming from permutation group 
theory. We write $S_\omega$ for the group of all permutations on a countably infinite set. 
The group $S_\omega$ is naturally equipped with the topology of pointwise convergence.
By the fundamental theorem of Engeler, Ryll-Nardzewski, and Svenonius, the reducts of an $\omega$-categorical
 structure $\Gamma$ are one-to-one correspondence 
 with the \emph{closed} subgroups of $S_\omega$ that contain the automorphism group of
 $\Gamma$.
 The automorphism groups of
 $\omega$-categorical structures are important and
 well-studied groups in permutation group theory,
 and classifications of reducts up to first-order interdefinability shed light on their nature.
 Indeed, all the classification results mentioned
 above make extensive use of the 
 group-theoretic perspective on reducts.


Let us also mention that
reducts of $(\mL;C)$ are
used for modeling various computational 
problems studied in phylogenetic reconstruction~\cite{Steel,BryantSteel,Bryant,NgSteelWormald,HenzingerKingWarnow,phylo-long}. When $\Gamma$ is such a structure
with a finite relational signature, then the \emph{constraint satisfaction problem (CSP)} 
for the \emph{template} $\Gamma$ is the problem to
decide for a finite structure $\Delta$ with the same signature as $\Gamma$ whether there exists a homomorphism from 
$\Delta$ to $\Gamma$ or not. For example, the CSP for $(\mL;C)$ itself has been called the \emph{rooted triple consistency problem} and it is known to be solvable in polynomial time by a non-trivial 
algorithm~\cite{ASSU,HenzingerKingWarnow,phylo-long}. 
Other phylogeny problems that can be
modeled as CSPs for reducts of $(\mL;C)$ are
the NP-complete \emph{quartet consistency problem}~\cite{Steel} and the NP-complete \emph{forbidden triples problem}~\cite{Bryant}. To classify the complexity of CSPs of reducts of an $\omega$-categorical structure,
a good understanding of the endomorphism monoids of these reducts is important;
for example, such a strategy has been used successfully in \cite{tcsps,BodPin-Schaefer-both,equiv-csps}. In this paper, we show that the endomorphism monoids of (L;C) fall into four categories.
In \cite{Phylo-Complexity} the authors give a full complexity classification for CSPs for reducts of (L;C)
and make essential use of this result.

\section{Results}\label{sect:results}
We show that there are only three
reducts of $(\mL;C)$ up to existential interdefinability (Corollary~\ref{cor:ex-reducts}).
In particular, there are only
three reducts of $(\mL;C)$ up to first-order interdefinability. 
The result concerning reducts up to first-order interdefinability can also be shown 
with a proof 
based on known results on Jordan permutation groups (Section~\ref{sect:autos}). 
However, we do not know how to obtain our stronger statement concerning reducts up to existential interdefinability
using Jordan group techniques. 

Our proof of Corollary~\ref{cor:ex-reducts} uses Ramsey theory for studying endomorphism
monoids of reducts of $(\mL;C)$.
More specifically, we use a Ramsey-type result for C-relations which is a special case of Miliken's theorem~\cite{Mil79} (see~\cite{BodirskyPiguet} for a short proof). 
We use it to show
that endomorphisms of reducts of $(\mL;C)$ must behave
\emph{canonically} (in the sense of Bodirsky \& Pinsker~\cite{BP-reductsRamsey}) on large parts of the domain and this
enables us to perform a combinatorial analysis of the endomorphism monoids. 
This approach provides additional insights which we describe next.

Assume that $\Gamma$ is a homogeneous structure
with a finite relational signature
whose age\footnote{The {\em age} of a relational structure $\Gamma$ is the set of finite structures that
are isomorphic to some substructure of $\Gamma$.}
has the Ramsey property (all examples mentioned above are reducts of such a structure).
Then, there is a 
general approach to analyzing reducts up to first-order interdefinability via the transformation monoids 
that contain $\Aut(\Gamma)$
instead of the closed permutation groups 
that contain $\Aut(\Gamma)$. 
This Ramsey-theoretic approach has been described in~\cite{BP-reductsRamsey}. 
We write $\omega^\omega$ for the transformation monoid of all unary functions 
on a countably infinite set. The monoid $\omega^\omega$ is naturally equipped with the topology of pointwise convergence
and the closed submonoids of $\omega^\omega$ that contain $\Aut(\Gamma)$ are in one-to-one correspondence with the
reducts of $\Gamma$ considered up to existential positive interdefinability.
We note that giving a complete description
of the reducts up to existential
positive interdefinability is usually difficult. For instance,
already the structure $(\mathbb N;=)$ admits infinitely many such reducts~\cite{BodChenPinsker}. However, it is often feasible to describe all reducts up to existential interdefinability; here, the Random Graph provides a good illustration~\cite{RandomMinOps}. In this paper, we show that it is feasible to describe all reducts of $(\mL;C)$ up to existential positive interdefinability. In particular, we show that the reducts of $(\mL;C)$ fall into four categories. An important category is when a reduct $\Gamma$ of $(\mL;C)$ has the same endomorphisms as the reduct $(\mL;Q)$. This reduct is a natural $D$-relation which is associated to $(\mL;C)$ (see Subsection \ref{sect:d-rel}), and its known complexity allows us to derive the complexity of the CSP for a large class of the reducts of  $(\mL;C)$. Those four categories are stated in the following main result of our paper.  
  
\begin{theorem}\label{thm:endos}
Let $\Gamma$ be a reduct of $(\mL;C)$. 
Then one of the following holds.
\begin{enumerate}
\item $\Gamma$ has the same endomorphisms as
$(\mL;C)$,
\item $\Gamma$ has a constant endomorphism,
\item $\Gamma$ is homomorphically equivalent to a reduct of $(\mL;=)$, or
\item $\Gamma$ has the same endomorphisms
as $(\mL;Q)$. 
\end{enumerate}
\end{theorem}

We use this result to identify in Corollary~\ref{cor:ex-reducts} below the reducts of $(\mL;C)$ up to existential interdefinability. 
The proof of Corollary~\ref{cor:ex-reducts} is based on a connection 
between existential and existential positive definability on the one hand, and the endomorphisms of $\Delta$ on the other hand.

\begin{proposition}[Proposition~3.4.7 
in~\cite{Bodirsky-HDR}]\label{prop:galois}
For every $\omega$-categorical 
structure $\Gamma$, it holds that
\begin{itemize}
\item a relation $R$ has an existential positive definition in $\Gamma$ if and only if $R$ is preserved by the endomorphisms of $\Gamma$ and
\item a relation $R$ has an existential definition in $\Gamma$ if and only if $R$ is preserved by the embeddings of $\Gamma$ into $\Gamma$. 
\end{itemize}
\end{proposition}

 \begin{corollary}\label{cor:ex-reducts}
 Let $\Gamma$ be a reduct of $(\mL;C)$. Then $\Gamma$ is existentially interdefinable
 with $(\mL;C)$, with $(\mL;Q)$, or with $(\mL;=)$. 
 \end{corollary}


Our result has important consequences for the study
of CSPs for reducts of $(\mL;C)$. 
To see this, note that when two structures $\Gamma$ and $\Delta$ are homomorphically equivalent, then they have the same CSP. 
Since the complexity of 
$\Csp(\Gamma)$ has been classified for all reducts
$\Gamma$ of $(\mL;=)$ (see Bodirsky and K\'{a}ra~\cite{ecsps}) and
since $\Csp(\Gamma)$ is trivial if $\Gamma$ has
a constant endomorphism, 
 our result
shows that we can focus on the case when
$\Gamma$ has the same endomorphisms
as $(\mL;C)$ or $(\mL;Q)$. This kind of simplifying assumptions have proven to be extremely important
in complexity classications of CSPs: examples include
Bodirsky \& K\'{a}ra~\cite{tcsps-journal} and Bodirsky \& Pinsker~\cite{BodPin-Schaefer-both}.




This article is organized as follows.
The structure $(\mL;C)$ is formally defined in Section~\ref{sect:c-rel}. 
We then show (in Section~\ref{sect:autos})  how to classify the reducts
of $(\mL;C)$ up to first-order interdefinability by using
known results about Jordan permutation groups. 
For the stronger classification up to existential definability, 
we investigate transformation monoids. 
The Ramsey-theoretic approach works well for studying transformation monoids and will be described in Section ~\ref{sect:tramsey}.
The main result is proved in 
Section~\ref{sect:canonical}.

%% file: phylogeny.tex
\section{Preliminaries}
\label{sect:c-rel}
We will now present some important definitions and results. We begin, in Section~\ref{sect:modtheory},
by providing a few preliminaries from model theory. Next, we define the
universal homogeneous binary branching C-relation $(\mL;C)$. 
There are several equivalent ways to do this and we consider two of them in Sections~\ref{sect:lc-amalgamation} and
\ref{sect:lc-axiom}, respectively.
The first approach is via
\Fresse-amalgamation and the second approach is an axiomatic approach based on
Adeleke and Neumann~\cite{AdelekeNeumann}. 
In Section~\ref{sect:d-rel}, we also give an axiomatic treatment of
an interesting reduct of $(\mL;C)$. 
In Section~\ref{sect:ordering}, we continue by introducing an 
ordered variant of the binary branching C-relation~\cite{PJC87} which
will be important in the later sections.

\subsection{Model theory}
\label{sect:modtheory}
We follow standard terminology as, for instance, used by Hodges~\cite{Hodges}. 
Let $\tau$ be a relational signature (all signatures in this paper will be relational)
and $\Gamma$ a $\tau$-structure.
When $R \in \tau$, we write $R^{\Gamma}$ 
for the relation denoted by $R$ in $\Gamma$;
we simply write $R$ instead of $R^{\Gamma}$ when
the reference to $\Gamma$ is clear. 
Let $\Gamma_1$ and $\Gamma_2$ be two $\tau$-structures
with domains $D_1$ and $D_2$, respectively, and
let $f \colon D_1 \to D_2$ be a function.
If $t = (t_1,\dots,t_k) \in (D_1)^k$,
then we write $f(t)$ for $(f(t_1),\dots,f(t_k))$, i.e. we
extend single-argument functions pointwise to sequences of arguments.
We say that $f$ \emph{preserves} $R$ iff $f(t) \in R^{\Gamma_2}$ whenever $t \in R^{\Gamma_1}$. 
If $X \subseteq D_1$ and $R \in \tau$ is a $k$-ary relation, then 
we say that $f$ \emph{preserves $R$ on $X$} if $f(t) \in R^{\Gamma_2}$ whenever $t \in R^{\Gamma_1} \cap X^k$. 
If $f$ does not preserves $R$ (on $X$), then we say that
$f$ \emph{violates} $R$ (on $X$). 

A function $f \colon D_1 \rightarrow D_2$ is an \emph{embedding} of $\Gamma_1$ into $\Gamma_2$ if 
$f$ is injective and has the property that for all
$R \in \tau$ (where $R$ has arity $k$) and all $t \in (D_1)^k$, we have $f(t) \in R^{\Gamma_2}$ if and only if $t \in R^{\Gamma_1}$. 

A \emph{substructure} of a structure
$\Gamma$ is a structure $\Delta$ with domain $S = D_{\Delta}\subseteq D_\Gamma$
and $R^{\Delta}=R^\Gamma \cap S^n$ for each $n$-ary $R \in \tau$; we also write $\Gamma[S]$ for $\Delta$. 
The \emph{intersection} $\Delta$ of two $\tau$-structures $\Gamma,\Gamma'$ is the structure with domain $D_{\Gamma}\cap D_{\Gamma'}$ and relations
$R^{\Delta}=R^{\Gamma}\cap R^{\Gamma'}$ for all $R\in\tau$; we also write $\Gamma \cap \Gamma'$ for $\Delta$. 


Let $\Gamma_1,\Gamma_2$ be $\tau$-structures such
that $\Delta=\Gamma_1\cap \Gamma_2$ is a substructure of both $\Gamma_1$ and $\Gamma_2$.  
A $\tau$-structure $\Delta'$ is an \emph{amalgam of $\Gamma_1$ and $\Gamma_2$ over
  $\Delta$} if for $i \in \{1,2\}$ there are embeddings $f_i$ of $\Gamma_i$ to $\Delta'$ such that $f_1(a)=f_2(a)$ for all $a\in D_\Delta$.  We assume that classes of structures are closed under isomorphism. 
  A class $\mathcal A$ of
$\tau$-structures has the \emph{amalgamation property} if for all
$\Delta,\Gamma_1,\Gamma_2\in\mathcal A$ with $\Delta=\Gamma_1\cap \Gamma_2$, 
there is a $\Delta' \in\mathcal A$ that is an
amalgam of $\Gamma_1$ and $\Gamma_2$ over $\Delta$. 
A class of finite
$\tau$-structures that  has the
  amalgamation property, is closed under isomorphism and
closed under taking substructures is called an
\emph{amalgamation class}. 

A relational structure $\Gamma$ is called \emph{homogeneous}
if all isomorphisms between finite  substructures can be
extended to automorphisms of $\Gamma$. A class $\mathcal{K}$ of $\tau$-structures has the joint embedding property if for any $\Gamma,\Gamma'\in \mathcal{K}$, there is $\Delta\in \mathcal{K}$ such that $\Gamma$ and $\Gamma'$ embed into $\Delta$. An amalgamation class has the joint embedding property since it always contains an empty structure. The following basic result is known as \Fresse's theorem.

\begin{theorem}[see
Theorem~6.1.2 in Hodges~\cite{Hodges}]\label{thm:fraisse}
Let $\mathcal A$ be an amalgamation class with countably many non-isomorphic members. Then there is a countable homogeneous $\tau$-structure $\Gamma$ such that $\mathcal A$ is the class of structures that embeds into $\Gamma$.
The structure $\Gamma$, which is unique up to isomorphism, is called the \emph{\Fresse-limit} of $\mathcal A$. 
\end{theorem}

\subsection{The structure $(\L;C)$: \Fresse-amalgamation}
\label{sect:lc-amalgamation}

We will now define the structure $(\L;C)$ as the \Fresse-limit of an appropriate amalgamation class. 
We begin by giving some standard terminology concerning rooted trees.
Throughout this article, a \emph{tree} is a simple, undirected, acyclic, and connected graph. 
A \emph{rooted tree} is a tree $T$ together with a distinguished vertex $r$ which is called 
the \emph{root} of $T$. 
The vertices of $T$ are denoted by $V(T)$.
The \emph{leaves} $L(T)$ of a rooted tree $T$ are the vertices of degree one that
are distinct from the root $r$. In this paper, a rooted tree is often drawn downward from the root.

For $u,v \in V(T)$, we say that $u$ \emph{lies below} $v$
if the path from $u$ to $r$ passes through $v$. 
We say that $u$ \emph{lies strictly below} $v$ if
$u$ lies below $v$ and $u \neq v$.
All trees in this article will be rooted and \emph{binary},
i.e., all vertices except for the root have either degree $3$ or $1$,
and the root has either degree $2$ or $0$.
A \emph{subtree} of $T$ is a tree $T'$ with $V(T') \subseteq V(T)$ and $L(T') \subseteq L(T)$. If the root
of $T'$ is different from the root of $T$, the 
subtree is called \emph{proper subtree}.
The \emph{youngest common ancestor (yca)} of a non-empty
finite set of vertices $S \subseteq V(T)$ is the (unique) node $w$
that lies above all vertices in $S$ and has maximal distance from $r$.

\begin{definition}\label{def:leaf-struct}
The \emph{leaf structure} of a binary rooted tree $T$ is the relational structure $(L(T);C)$ 
where $C(a,bc)$ holds in $C$ if and only if $\yca(\{b,c\})$ lies strictly below $\yca(\{a,b,c\})$ in $T$.
We call $T$ the \emph{underlying tree} of the leaf structure. 
\end{definition}

We mention that the definition of C-relation on binary rooted trees can also be obtained from the relation $\mid$ on trees with a distinguished leaf~\cite{PJC87}. The slightly non-standard way of writing the arguments of the relation $C$
has certain advantages that will be apparent in forthcoming sections.

\begin{definition}
For finite non-empty $S_1,S_2 \subseteq L(T)$, we write $S_1|S_2$ if 
neither of $\yca(S_1)$ and $\yca(S_2)$ lies below
the other. For sequences of (not necessarily distinct) vertices $x_1,\dots,x_n$ and $y_1,\dots,y_m$ we write $x_1,\dots,x_n|y_1,\dots,y_m$ if $\{x_1,\dots,x_n\}\ | \{y_1,\dots,y_m\}$.
\end{definition}

In particular, $xy|z$ (which is the notation that is typically used in the literature on phylogeny problems) is equivalent to $C(z,xy)$;
it will be very convenient to have both notations available. Note that 
if $xy|z$ then this includes the possibility that $x=y$; however, $xy|z$ implies that
$x \neq z$ and $y \neq z$. 
Hence, for every triple $x,y,z$ of leaves in a rooted binary tree, we either have $xy|z$, $yz|x$, $xz|y$, or $x=y=z$. 
Also note that $x_1,\dots,x_n|y_1,\dots,y_m$ if and only if $x_ix_j|y_k$ and $x_i|y_ky_l$ for
all $i,j \leq n$ and $k,l \leq m$. 
The following result is known but we have been unable to find an explicit proof in the literature. 
Hence, we give a proof for the convenience of the reader.

\begin{proposition}\label{prop:amalgam}
The class $\cal C$ of all finite leaf structures is an amalgamation class.
\end{proposition}

\begin{proof}
Arbitrarily choose 
$B_1,B_2 \in \cal C$ such that $A=B_1 \cap B_2$ is a substructure of both $B_1$ and $B_2$. 
We inductively assume that the statement has been shown for all triples
$(A,B_1',B_2')$ where $D(B_1') \cup D(B_2')$ is a proper subset
of $D(B_1) \cup D(B_2)$.

Let $T_1$ be the rooted binary tree underlying $B_1$ and $T_2$ the rooted binary tree underlying $B_2$.
Let $B_1^1 \in \cal C$ be the substructure of $B_1$ induced by the vertices below the left child of $T_1$
and $B_1^2 \in \cal C$ be the substructure of $B_1$ induced by the vertices below the right child of $T_1$. 
The structures $B_2^1$ and $B_2^2$ are defined
analogously for $B_2$.

First consider the case when there is a vertex $u$ that lies in both $B_1^1$ and $B_2^1$ and a vertex $v$ that lies in both $B_2^1$ and $B_1^2$. We claim that in this case no vertex $w$ from $B_2^2$ can lie inside $B_1$.
Assume the contrary and note
that $w$ is either in $B_1^1$, in which case we have $uw|v$ in $B_1$,
or in $B_1^2$, in which case we have $u|vw$ in $B_1$.
But since $u,v,w$ are in $A$, this contradicts the fact that $uv|w$ holds
in $B_2$. Let $C' \in \cal C$ be the amalgam of $B_1$ and $B_2^1$ over $A$ (which exists by the inductive 
assumption) and let $T'$ be its underlying tree.
Consider a tree $T$ with root $r$, $T'$ as its left subtree, and the
underlying tree of $B_2^2$ as its right subtree.
It is straightforward to verify that the leaf structure of $T$ is in $\cal C$
and that it is an amalgam of $B_1$ and $B_2$ over $A$.

The above argument can also be applied to the cases where the role of $B_1$ and $B_2$, or the role of $B_1^1$ with $B_1^2$,
or the role of $B_2^1$ with $B_2^2$ are exchanged.  
Hence, the only remaining essentially different 
case we have to consider is
when $D(B_1^1) \cup D(B_2^1)$ and $D(B_1^2) \cup D(B_2^2)$ are disjoint.
In this case, it is straightforward to first amalgamate $B_1^1$ with $B_2^1$
and $B_1^2$ with $B_2^2$ to obtain the amalgam of $B_1$ and $B_2$;
the details are left to the reader.
\end{proof}


We write $(\mL;C)$ for the \Fresse-limit of $\cal C$. 
Obvious reducts of $(\mL;C)$ are $(\mL;C)$ itself and $(\mL;=)$. 
To define a third reduct, consider
the 4-ary relation $Q(xy,uv)$ with the following first-order definition over $(\L;C)$: 
\begin{align*}
(xy|u \wedge xy|v) \vee (x|uv \wedge y|uv)
\end{align*}
This relation is often referred to as the {\em quartet} relation~\cite{Steel}.


\subsection{The structure $(\L;C)$: an axiomatic approach}
\label{sect:lc-axiom}

The structure $(\mL;C)$ that we defined in Section~\ref{sect:lc-amalgamation} 
is an important example of a so-called
\emph{$C$-relation}. This concept was introduced by Adeleke \& Neumann~\cite{AdelekeNeumann}
and we closely follow their definitions in the following.
A ternary relation $C \subseteq X^3$ is said to be a C-relation on $X$ if 
the following conditions hold:
\begin{enumerate}
\item [C1.] $\forall a,b,c \: \big(C(a,bc) \Rightarrow C(a,cb)\big)$
\item [C2.] $\forall a,b,c \: \big( C(a,bc) \Rightarrow \neg C(b,ac) \big)$ 
\item [C3.] $\forall a,b,c,d \: \big( C(a,bc) \Rightarrow C(a,dc) \vee C(d,bc) \big)$ 
\item [C4.] $\forall a,b \: \big(  a \neq b \Rightarrow C(a,b,b) \big)$
\end{enumerate}
A C-relation is called \emph{proper}
if it satisfies two further properties:
\begin{enumerate}
\item [C5.] $\forall a,b \, \exists c \: \big( C(c,ab)\big)$
\item [C6.] $\forall a,b \: \big( a \neq b \Rightarrow \exists c (c \neq b \wedge C(a,bc))\big)$
\end{enumerate}

These six axioms do \emph{not} describe
the \Fresse-limit $(\mL;C)$ up to isomorphism. 
To completely axiomatize the theory of $(\mL;C)$,
we need two more axioms.

\begin{enumerate}
\item [C7.] $\forall a,b,c \: \big( C(c,ab) \Rightarrow \exists e \: (C(c,eb) \wedge C(e,ab)) \big)$
\item [C8.] $\forall a,b,c \: \big(  (a \neq b \vee a \neq c \vee b \neq c) \Rightarrow (C(a,bc) \vee C(b,ac) \vee C(c,ab)) \big)$
\end{enumerate}

\noindent C-relations that satisfy C7 are called \emph{dense} and
C-relations that satisfy C8 are called \emph{binary branching}. Note that C1-C8 are satisfiable since $(\mL;C)$ is a countable model of C1-C8.

We mention that the structure $(\mL;C)$ 
is \emph{existential positive complete} within the class of all $C$-relations, as defined in~\cite{BodHilsMartin}: for every homomorphism $h$ of $(\mL;C)$ into another C-relation and every existential positive formula $\phi(x_1,\dots,x_n)$ and all $p_1,\dots,p_n \in \mL$ such that 
$\phi(h(p_1),\dots,h(p_n))$ holds we have that 
$\phi(p_1,\dots,p_n)$ holds in $(\mL;C)$, too. 
It is also easy to see that 
every existential positive complete C-relation 
must satisfy C7 and C8. These facts about
existential positive completeness of $(\mL;C)$ are not needed in the remainder of the article, but 
together with Lemma~\ref{lem:c-homo} below they demonstrate that the structure $(\mL;C)$ can
be seen as the (up to isomorphism unique) \emph{generic} countable C-relation.

The satisfiability of C1-C8 can also be shown using the idea of constructing C-relations in \cite[page 123]{Infinite-perm-groups}. Let $\mathcal{F}$ be the set of functions $f\colon(0,\infty)\to\{0,1\}$, where $(0,\infty)$ denotes the set of positive rational numbers with the standard topology, such that the following conditions hold. 
\begin{itemize}
\item there exists $a\in (0,\infty)$ such that $f(x)=0$ for every $x\in (0,a)$.
\item $f$ has finitely many points of discontinuity and for each point $b$ of discontinuity, there exists $\epsilon\in (0,b)$ such that $f(x)\neq f(b)$ for every $x\in (b-\epsilon,b)$.
\end{itemize}
 For every $f,g\in \mathcal{F}$ such that $f\neq g$, let $\Com(f,g)$ denote the interval $(0,c)$ such that $f(x)=g(x)$ for every $x\in (0,c)$, and $f(c)\neq g(c)$. If $f=g$, let $\Com(f,g):=(0,\infty)$. Note that $c$ is a point of discontinuity of either $f$ or $g$. We define a relation $C$ on $\mathcal{F}$ by $C(f,gh)$ if $\Com(f,h)\subsetneq \Com(g,h)$. We can easily to verify that $(\mathcal{F};C)$ is a countable model of C1-C8.
 
We will now prove (in Lemma~\ref{lem:c-homo}) that there is a unique countable model of C1-C8 up to isomorphism. It suffices to show that if $\Gamma$ is a countable structure with signature $\{C\}$ 
satisfying C1-C8, then $\Gamma$ is isomorphic to $(\mL;C)$.  To do so, we need a number of observations (Lemmas~\ref{consequences}, \ref{lem:split}, and \ref{lem:pres-c}).

The following consequences of C1-C8 are used in the proofs without further notice. 
\begin{lemma}[C-consequences]
\label{consequences}
Let $C$ denote a C-relation. Then

\begin{itemize}
  \item[1. ] $\forall x,y,z,t \: \big((C(x,yz)\land C(x,yt))\Rightarrow C(x,zt)\big)$,
  \item[2. ] $\forall x,y,z,t \: \big((C(x,zt)\land C(z,xy))\Rightarrow (C(t,xy)\land C(y,zt))\big)$, and 
  \item[3. ] $\forall x,y,z,t \: \big((C(z,xy)\land C(y,xt))\Rightarrow (C(z,yt)\land C(z,xt))\big)$.
\end{itemize}  
\end{lemma}
\begin{proof}
We prove the first consequence. The others can be shown analogously. 
Assume to the contrary that $C(x,zt)$ does not hold. 
By applying C3 to $x,y,z,t$, we get that $C(t,yz)$
and C2 implies that $C(z,yt)$ does not hold. 
By applying C3 to $x,y,t,z$, it follows that $C(x,zt)$ holds
and we have a contradiction.
\end{proof}
For  two subsets $Y,Z$ of $X$, we write $C(Y,Z)$ if

\begin{enumerate}
\item
 $C(y,z_1 z_2)$ for 
arbitrary $y \in Y$ and $z_1,z_2\in Z$, and

\item
$C(z,y_1y_2)$ for arbitrary $y_1,y_2 \in Y$ and $z\in Z$.
\end{enumerate}

\begin{lemma}
\label{lem:split}
Let $C$ be a ternary relation on a countably 
infinite set $X$ that satisfies \emph{C1-C8}. 
Then for every finite subset $Y$ of $X$ of size at least $2$ there are two non-empty subsets $A,B$ of $Y$ such that $A\cup B=Y$ and $C(A,B)$.
\end{lemma}
\begin{proof}
We prove the lemma by induction on $|Y|$. Clearly, the claim holds if $|Y|=2$ 
so we assume that the lemma holds when $|Y| = k-1$ for some $k > 2$.
Henceforth, assume $|Y|=k$.
Arbitrarily choose $y \in Y$ and let $Y'=Y\backslash \set{y}$. By the induction hypothesis, 
there are two non-empty subsets $A',B'$ of $Y'$ such that $A'\cup B'=Y'$ and 
$C(A',B')$. Pick $a'\in A'$ and $b'\in B'$. One of the following holds. 
\begin{itemize}
  \item $C(y,a'b')$. Arbitrarily choose $c,d \in Y'$. We show that $C(y,cd)$ holds. If $c,d \in A'$, then we have 
$C(y,a'b')$, $C(b',a'c)$, and $C(b',a'd)$. It follows immediately 
from Lemma~\ref{consequences} that $C(y,cd)$. 
Analogously, $C(y,cd)$ holds if $c,d\in B'$. It remains to consider the case $c\in A',d\in B'$. 
Here, we have $C(y,a'b'),C(b',a'c)$ and $C(a',b'd)$. Once again, it 
follows from Lemma~\ref{consequences} 
that $C(y,cd)$ holds. By setting $A=\{y\}$ and $B=A'\cup B'$, the lemma of the lemma follows.
 \item $C(b',ya')$. We first show that for arbitrary $a''\in A'$ and $b''\in B'$, we have $C(b'',a''y)$. This follows from 
Lemma \ref{consequences} and the fact that $C(b',ya'),C(b',a'a'')$, and $C(a',b'b'')$ hold. 
We can now show that for arbitrary 
$b'',b'''\in B'$, we have that $C(y,b''b''')$ holds. 
This follows from Lemma \ref{consequences} and the fact that
$C(b',a'y)$, $C(a',b'b'')$, and $C(a',b'b''')$ hold. 
This implies that $C(A'\cup \set{y},B')$ and
we have proved
the induction step by setting $A=A' \cup \set{y}$ and $B=B'$. 
 \item $C(a',yb')$. This case can be proved analogously to the previous case: we get that
$A=A'$ and $B=B' \cup \set{y}$.
\end{itemize} 
The case distinction is exhaustive because of C8. 
\end{proof}

We would like to point out an important property of maps that preserve $C$.

\begin{lemma}\label{lem:pres-c}
Let $e \colon X \to \mL$ for $X \subseteq \mL$
be a function that preserves $C$. Then $e$ is injective and preserves the relation 
$\{t \in \mL^3 : t \notin C\}$.
\end{lemma}
\begin{proof}
Clearly, $e$ preserves the binary relation $\{(x,y) \in \mL^2 : x \neq y\} = \{(x,y) : \exists z. C(x,y,z)\}$, and so
$e$ is injective. 
Arbitrarily choose $u_1,u_2,u_3 \in \mL$ such that 
$u_1|u_2u_3$ does \emph{not} hold. 
If $|\{u_1,u_2,u_3\}|=1$ then $e(u_1)|e(u_2)e(u_3)$
does not hold and there is nothing to show.
If $|\{u_1,u_2,u_3\}|=2$ then by C4 either $u_1=u_2$ or $u_1 = u_3$, and $e(u_1)|e(u_2)e(u_3)$
does not hold. 
If $|\{u_1,u_2,u_3\}|=3$ then
by C6 we have either
$u_2|u_1u_3$,
or $u_3|u_1u_2$. It follows that 
$e(u_2)|e(u_1)e(u_3)$
or $e(u_3)|e(u_1)e(u_2)$. 
In both cases, 
$e(u_1)|e(u_2)e(u_3)$ does not hold by C2.
\end{proof}

We will typically use the contrapositive version of Lemma~\ref{lem:pres-c} in the sequel. This allows to draw the conclusion $a|bc$ under the assumption $e(a)|e(b)e(c)$.

\begin{lemma}\label{lem:c-homo}
Let $\Gamma$ be a countable structure with signature $\{C\}$ that satisfies \emph{C1-C8}. Then $\Gamma$ is isomorphic to $(\mL;C)$. 
\end{lemma}
\begin{proof}
It is straightforward (albeit a bit tedious) to show that $(\mL;C)$ satisfies C1-C8. It then remains to show that if 
$\Gamma_1$ and $\Gamma_2$ are two countably infinite $\{C\}$-structures that satisfy C1-C8, then the two structures are isomorphic. Let $X_1,X_2$ denote the domains of $\Gamma_1$ and $\Gamma_2$, respectively. 
This can be shown by a back-and-forth argument based on the
following claim. 

\noindent
\textbf{Claim:} Let $A$ be a non-empty finite subset of $X_1$ and let $f$ denote a 
map from $A$ to $X_2$ 
that preserves $C$. Then for every $a \in X_1$, the map $f$ can be extended to a map $g$ from 
$A\cup \set{a}$ to $X_2$ that preserves $C$. 

It follows from Lemma~\ref{lem:pres-c} that $f$
also preserves $\set{(x,y,z):\neg C(x,yz)}$. 
We prove the claim by induction 
on $|A|$. Clearly, we are done if $a \in A$ or $|A|=1$. Hence, assume that $a \not \in A$ 
and $|A|\geq 2$. Let $A_1,A_2$ be subsets of $A$ such that $A_1\cup A_2=A$ and 
$C(A_1,A_2)$, which exist due to Lemma~\ref{lem:split}. 
Note that $C(f(A_1),f(A_2))$ holds 
in $(X_2;C)$. Pick $a_1\in A_1$ and $a_2\in A_2$. We  construct the map $g$ in each of the following three cases.
\begin{itemize}
  \item $C(a,a_1a_2)$. 
 We claim that $C(\set{a},A)$ holds.
Arbitrarily choose $u,v \in A$. Then
either $C(ua_1,a_2)$ or $C(a_1,a_2u)$ by the choice of $a_1$ and $a_2$. Similarly, either $C(va_1,a_2)$ or $C(a_1,a_2v)$. 
So there are four cases to consider; we only treat the case $C(ua_1,a_2)$ and $C(a_1,va_2)$ since the other cases are similar or easier. Now, $C(ua_1,a_2)$ and $C(a_1a_2,a)$ 
imply that $C(a,u a_1)$ by item 3 of Lemma~\ref{consequences}.
Similarly, we have $C(a, v a_2)$.
Now $C(a_1a_2,a)$ 
and two applications of item 1 of Lemma~\ref{consequences} give 
$C(u v,a)$, which proves the subclaim. 

It follows from C5 that there exists an $a'\in X_2$ such that $C(a',f(a_1) f(a_2))$.
Once again, we obtain
$C(\set{a'},f(A))$ as a consequence of 
Lemma~\ref{consequences}.
This implies that 
the map $g \colon A \cup \{a\} \to X_2$, defined by $g|_A=f$ and $g(a)=a'$, preserves $C$.
  \item $C(a_2,aa_1)$. It follows from Lemma~\ref{consequences} that $C(\set{a}\cup A_1,A_2)$ holds. 
We consider the following cases. 

\smallskip

\underline{$|A_1|=1$.} 
There exists $a' \in X_2$ such that $C(f(a_2),a' f(a_1))$ by C6, and 
Lemma~\ref{consequences} implies that $C(\set{f(a_1),a'},f(A_2))$ holds. Since we also have $C(\set{a,a_1},A_2)$, the map $g \colon A \cup \set{a} \to X_2$ defined by 
$g|_{A}=f$ and $g(a)=a'$, preserves $C$. 

\underline{$|A_2|=1$.} This case can be treated analogously to the previous case. 

\underline{$|A_1|\geq 2$ and $|A_2| \geq 2$.} 
Let $B_1,B_2$ be non-empty such that $C(B_1,B_2)$ and $A_1=B_1\cup B_2$. 
Arbitrarily choose $b_1\in B_1,b_2\in B_2$. The following cases are exhaustive by C8. 
\begin{itemize}
  \item $C(a,b_1b_2)$. It is a direct consequence of C7 that there exists an $a'\in X_2$ such that both 
$C(\set{f(a_2)},\set{f(b_1),f(b_2),a'})$ 
and $C(a',f(b_1) f(b_2))$ hold. Furthermore, Lemma \ref{consequences} implies 
that $C(\set{a},A_1)$,$C(\set{a}\cup A_1,A_2)$, $C(\set{a'},f(A_1))$, and 
$C(\set{a'}\cup f(A_1),f(A_2))$. Hence, the map $g \colon A\cup \set{a}\to X_2$, 
defined by $g|_{A}=f$ and $g(a)=a'$, preserves $C$.
  \item $C(b_2,a b_1)$. By assumption we know that $|A_2| \geq 2$, and since $A_1 \cup A_2 = A$ it follows that $|A_1 \cup \set{a}| < |A|$. Hence, by the induction hypothesis there exists a map $h \colon A_1\cup \set{a}\to X_2$ such that 
$h|_{A_1}=f|_{A_1}$ and $h$ preserves $C$ on $A_1\cup \set{a}$. Since $h$ preserves $C$ on 
$A_1\cup \set{a}$, we see that $C(h(b_2),h(a) h(b_1))$, and consequently that 
$C(f(b_2),h(a) f(b_1))$ holds. Since 
both
$C(A_1,A_2)$ and $C(f(A_1),f(A_2))$ hold, it follows from Lemma \ref{consequences} that 
$C(\set{a}\cup A_1,A_2)$ and $C(\set{h(a)}\cup f(A_1),f(A_2))$ hold. 
This implies that the map 
$g \colon A\cup \set{a}\to X_2$, defined by $g|_{A_1\cup \set{a}}=h,g|_{A_2}=f|_{A_2}$, preserves $C$.
  \item $C(b_1,a b_2)$.
  The proof is analogous to the case above.
\end{itemize} 
\item $C(a_1,a a_2)$. The proof is analogous to the case when $C(a_2,aa_1)$.
\end{itemize}  
The case distinction is exhaustive because of C8.
\end{proof}


\subsection{The reduct $(\mL;Q)$}
\label{sect:d-rel}
The reduct $(\mL;Q)$ of $(\mL;C)$ can
be treated axiomatically, too. 
A 4-ary relation $D$ is said to be a \emph{D-relation} on $X$
if 
the following conditions hold:
\begin{enumerate}
\item [D1.] $\forall a,b,c,d \: \big( D(ab,cd) \Rightarrow D(ba,cd) \wedge D(ab,dc) \wedge D(cd,ab) \big)$
\item [D2.] $\forall a,b,c,d \: \big(  D(ab,cd) \Rightarrow \neg D(ac,bd) \big)$
\item [D3.] $\forall a,b,c,d,e \: \big( D(ab,cd) \Rightarrow D(eb,cd) \vee D(ab,ce) \big)$ 
\item [D4.] $\forall a,b,c \: \big( (a \neq c \wedge b \neq c) \Rightarrow D(ab,cc) \big)$
\end{enumerate}
A D-relation is called \emph{proper}
if it additionally satisfies the following condition:
\begin{enumerate}
\item [D5.] For pairwise distinct $a,b,c$ there is $d \in X \setminus \{a,b,c\}$ with $D(ab,cd)$.
\end{enumerate}


As with $(\mL;C)$, it is possible
to axiomatize the theory of $(\mL;Q)$
by adding finitely many axioms. 

\begin{itemize}
\item [D6.] $\forall a,b,c,d \: \big( D(ab,cd) \Rightarrow \exists e \: (D(eb,cd) \wedge D(ae,cd) \wedge D(ab,ed) \wedge D(ab,ce)) \big )$ 
\item [D7.] $\forall a,b,c,d \: \big(  |\{a,b,c,d\}|\geq 3 \Rightarrow (D(ab,cd) \vee D(ac,bd) \vee D(ad,bc)) \big)$
\end{itemize}

D-relations satisfying D6 are called \emph{dense},
and D-relations satisfying D7 are called \emph{binary branching}. 

We will continue by proving that if two countable
structures with signature $\{D\}$ satisfy D1-D7,
then they are isomorphic.
For increased readability, we write
$D(xyz,uv)$ when $D(xy,uv)\land D(xz,uv)\land D(yz,uv)$, and we write $D(xy,zuv)$ when 
$D(xy,zu)\land D(xy,zv)\land D(xy,uv)$. One may note, for instance, that
$D(xy,zuv)$ is equivalent to $D(xy,uvz)$.

\begin{lemma}[D-consequences]
\label{D-consequences}
If $D$ is a D-relation, then
\begin{itemize}
\item $\forall x,y,z,u,v \big((D(xy,zu)\land D(xy,zv))\Rightarrow D(xy,uv)\big)$, and
\item $\forall x,y,z,u,v \big(D(xy,zu)\Rightarrow (D(xyv,zu)\lor D(xy,zuv))\big)$.
\end{itemize}
\end{lemma}
\begin{proof}
We prove the first item. By applying D1 and D3 to $D(xy,zu)\land D(xy,zv)$, we get that
\[(D(yv,uz) \vee D(xy,uv)) \wedge (D(yu,vz) \vee D(xy,uv)).\]
If $D(xy,uv)$ does not hold, then $D(yv,uz) \wedge D(yu,vz)$ must
hold. However, this immediately leads to a contradiction via
D2: $D(yu,vz) \Rightarrow \neg D(yv,uz)$.

To prove the second item, assume that $D(xy,zu)$ holds and arbitrarily choose $v$. By D3,
we have $D(vy,zu) \vee D(xy,zv)$.
Assume that $D(vy,zu)$ holds; the other case can be proved in a similar
way. By definition $D(xyv,zu)$ if and only if $D(xy,zu) \wedge D(xv,zu) \wedge D(yv,zu)$. 
We know that $D(xy,zu)$ holds and that $D(vy,zu)$ implies
$D(yv,zu)$ via D1. It remains to show that $D(xv,zu)$ holds, too. By once again applying D1, we
see that
$D(zu,yx) \wedge D(zu,yv)$. It follows that $D(zu,xv)$ holds by the claim above
and we conclude that $D(xv,zu)$ holds by D1.
\end{proof}

\begin{lemma}\label{lem:pres-q}
Let $e \colon X \to \mL$ for $X \subseteq \mL$
be a function that preserves $Q$. Then $e$ is injective
and preserves the relation $\{q \in \mL^4 : q \notin Q\}$.
\end{lemma}
\begin{proof}
The proof is very similar to the proof of Lemma~\ref{lem:pres-c} and is left to the reader.
\end{proof}
We will typically use the contrapositive version of Lemma~\ref{lem:pres-q} in the sequel. This allows to draw the conclusion $Q(ab,cd)$ under the assumption $Q(e(a)e(b),e(c)e(d))$.

\begin{lemma}\label{lem:q-homo}
Let $D$ be a 4-ary relation on a countably infinite set $X$ that satisfies \emph{D1-D7}. Then $(X;D)$
is isomorphic to $(\mL;Q)$, and homogeneous. 
\end{lemma}
The proof of Lemma~\ref{lem:q-homo} is based on Lemma~\ref{lem:c-homo} and the idea of rerooting at a fixed leaf to create a C-relation from a D-relation. The idea of rerooting was already discussed in \cite{PJC87}.  

\begin{proof}
It is straightforward to verify that $(\mL;Q)$ satisfies D1-D7.
Let $(X_1;D)$ and $(X_2;D)$ be two countably infinite 
sets that satisfy D1-D7, let $Y_1$ be a finite subset of $X_1$, and $\alpha$ an embedding of the structure induced by $Y_1$ in $(X;D)$ into $(X_2,D)$. We will show that $\alpha$
can be extended to an isomorphism between $(X_1;D)$
and $(X_2,D)$. This can be applied
to $(X_1;D) = (X_2;D) = (\mL;Q)$ and hence also shows homogeneity of $(\mL;Q)$. 
  
Arbitrarily choose $c \in Y_1$.
We define a relation $C$ on $X_1' := X_1 \setminus \{c\}$ as follows:
for every $(x,y,z) \in (X_1')^3$, let $(x,y,z)\in C$ if and only if $D(cx,yz)$ holds. Similarly, we define a relation $C$ 
on $X_2' := X_2 \setminus \{\alpha(c)\}$ as follows: for every $(x,y,z) \in (X_1')^3$, let $(x,y,z)\in C$ if and only if $D(\alpha(c)x,yz)$ holds. 

One can verify that both $(X_1',C)$ and $(X_2',C)$ 
satisfies C1-C8. It follows from Lemma~\ref{lem:c-homo} that $(X_1';C)$ 
and $(X_2';C)$ are isomorphic to $(\mL;C)$, 
and it follows from homogeneity of $(\mL;C)$ 
that the restriction of $\alpha$ to 
$Y_1 \setminus \{c\}$ can be extended
to an isomorphism $\alpha'$ between $(X_1';C)$ and $(X_2';C)$. 

We conclude the proof by showing that the map $\beta \colon X_1\to X_2$, defined by 
$\beta(c) := \alpha(c)$ and $\beta|_{X_1'}=\alpha'$, is an isomorphism between $(X_1;D)$ and $(X_2;D)$.
Arbitrarily choose $x,y,u,v \in X_1$ satisfying $D(xy,uv)$. By Lemma~\ref{lem:pres-q} it is sufficient to show that 
$D(\beta(x)\beta(y),\beta(u)\beta(v))$. Clearly, we are done if $x,y,u,v$ are not 
pairwise distinct, or if $c \in \set{x,y,u,v}$, so assume otherwise. 
By Lemma~\ref{D-consequences} we have $D(xy c,uv)$
or $D(xy,c uv)$. In the former case, 
it follows from the definition of $C$ 
  on $X_1'$ and $X_2'$ that $D(\alpha(x) \alpha(c),\alpha(u) \alpha(v))$ and 
$D(\alpha(y) \alpha(c),\alpha(u)\alpha(v))$. Lemma~\ref{D-consequences} 
implies that $D(\alpha(x) \alpha(y),\alpha(u)\alpha(v))$, which is equivalent to 
$D(\beta(x)\beta(y),\beta(u)\beta(v))$.
The case that  $D(xy,c uv)$ can be shown analogously to the previous case.
\end{proof}

\begin{corollary}\label{prop:rer-exists}
There exists an operation $\rer \in \Aut(\mL;Q)$ that violates $C$. 
\end{corollary}
\begin{proof}
It follows from Lemma \ref{lem:q-homo} that $\Aut(\mL;Q)$ is 3-transtive. Since  $\Aut(\mL;C)$ is $2$-transitive, but not $3$-transitive, it follows that $\Aut(\mL;C)\neq \Aut(\mL;Q)$. Since $\Aut(\mL;C)\subseteq \Aut(\mL;Q)$, there is $\rer\in \Aut(\mL;Q)$ which violates $C$.
\end{proof}
The name $\rer$ may seem puzzling at first sight: it is short-hand for {\em rerooting}. The choice of terminology will be clarified in the next section. 


\subsection{Convex orderings of C-relations}
\label{sect:ordering}
In the proof of our main result, it will be useful to work with an expansion 
$(\L;C,\prec)$ of $(\L;C)$ by a certain linear order 
$\prec$ on $\mathbb L$. We will next describe how this linear order
is defined as a \Fresse-limit.  
A linear order $\prec$ on the elements of a leaf structure
$(L;C)$ is called \emph{convex}
if for all $x,y,z \in L$ with $x \prec y \prec z$ we have that either $xy|z$ or that $x|yz$ (but not $xz|y$). The concept of convex linear order was already discussed in \cite{PJC87} and in \cite[page 162]{HMMNT00}.

\begin{proposition}\label{prop:convex-a}
Let $(L(T);C)$ be the leaf structure
of a finite binary rooted tree $T$ and arbitrarily choose
$a \in L(T)$. Then there exists 
a convex linear order $\prec$ of $L(T)$ 
whose maximal element is $a$. In particular,
every leaf structure can be expanded
to a convexly ordered leaf structure. 
\end{proposition}
\begin{proof}
Perform a depth-first search of $T$, starting at the root, such that vertices that lie above $a$ in $T$
are explored latest possible during the search.
Let $\prec$ be the order on $L(T)$ in which
the vertices have been visited during the search. Clearly, $\prec$
is convex and $a$ is its largest element.
\end{proof}


\begin{proposition}\label{prop:convex-order}
The class $\mathcal C'$ of all finite convexly ordered leaf structures is an amalgamation class and
its \Fresse-limit is isomorphic to an 
expansion $(\L;C,\prec)$ of $(\L;C)$ by a convex linear 
ordering $\prec$. The structure $(\mL;C,\prec)$ is described uniquely up to isomorphism by
the axioms C1-C8 and by the fact that $\prec$ is a dense and unbounded linear order which is
convex with respect to $(\mL;C)$. 
\end{proposition}
\begin{proof}
The proof that $\cal C'$ is an amalgamation class
is similar to the proof of Proposition~\ref{prop:amalgam}.
The \Fresse-limit of $\mathcal C'$ clearly 
satisfies C1-C8, it is equipped with a
convex linear order, and all countable structures with these properties are 
 in fact isomorphic; this can be shown by a back-and-forth argument.
By Lemma~\ref{lem:c-homo}, the structure
obtained by forgetting the order 
is isomorphic to $(\mL;C)$ and the statement follows. 
\end{proof}

By the classical result of Cantor~\cite{CantorDicht},
all countable dense unbounded linear orders are isomorphic to 
$(\mQ;<)$, and hence Proposition~\ref{prop:convex-order} implies that $(\mL;\prec)$ is isomorphic to $(\mQ;<)$.

%% file: autos.tex
\section{Automorphism groups of reducts}
\label{sect:autos}
We will now 
show that the structure $(\mL;C)$ has precisely three reducts up to first-order interdefinability. 
Our proof
uses a result by Adeleke and Neumann~\cite{AdelekeNeumannPrimitive} about primitive permutation groups
with \emph{primitive Jordan sets}.
The link between reducts of $(\mL;C)$
and permutation groups is given by the theorem of Engeler, Ryll-Nardzewski, and Svenonius, which we briefly recall
in Section~\ref{sect:perm-groups}. 
We continue in Section~\ref{sect:rerooting} by presenting some important
lemmata about functions that preserve $Q$ but
violate $C$. 
With these results in place,
we finally prove the main result of this section
in Section~\ref{sect:group-class}.

\subsection{Permutation group preliminaries}
\label{sect:perm-groups}
Our proof will utilize links between homogeneity, $\omega$-categoricity, and permutation groups
so we begin by discussing these central concepts.
A structure $\Gamma$ is {\em $\omega$-categorical} if
all countable structures that satisfy the same first-order sentences as 
$\Gamma$ are isomorphic (see e.g.\,Cameron~\cite{Oligo} or Hodges~\cite{Hodges}). 
Homogeneous structures with finite relational signatures are $\omega$-categorical, so the structure
$(\mL;C)$ is $\omega$-categorical.
Moreover, all structures with a first-order definition
in an $\omega$-categorical structure are $\omega$-categorical (see again Hodges~\cite{Hodges}).
This implies, for instance, that $(\mL;Q)$ is $\omega$-categorical.

\begin{figure}
\begin{center}
\includegraphics[scale=.7]{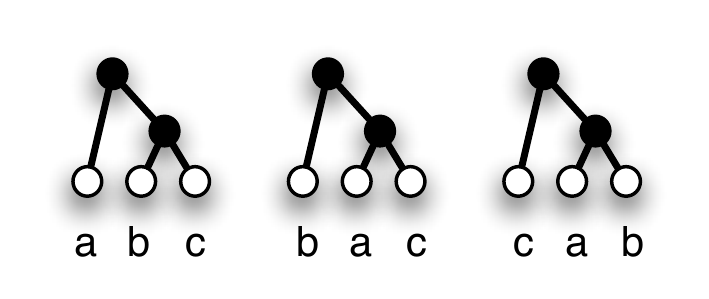} 
\end{center}
\caption{Illustration of the 3 orbits of triples $(a,b,c)$
with pairwise distinct entries of $\Aut(\L;C)$.}
\label{fig:3-orbits}
\end{figure}

The fundamental theorem by Engeler, Ryll-Nardzewski, and Svenonius 
is a 
characterization of $\omega$-categoricity
in terms of permutation groups.
When $G$ is a permutation group on a set $X$,
then the \emph{orbit} of a $k$-tuple $t$
is the set $\{\alpha(t) \; | \; \alpha \in G\}$. 
We see that homogeneity of $(\mL;C)$ implies that $\Aut(\mL;C)$ has precisely three 
orbits of triples with pairwise distinct entries; an illustration of these orbits can be found
in Figure~\ref{fig:3-orbits}. 
We now state the Engeler-Ryll-Nardzewski-Svenonius theorem and its proof
can be found in, for instance, Hodges~\cite{Hodges}.

\begin{theorem}
\label{thm:ryll}
A countable relational structure $\Gamma$ is $\omega$-categorical if and only if 
the automorphism group of $\Gamma$ is \emph{oligomorphic}, that is,
 if for each $k \geq 1$ there are finitely many orbits of $k$-tuples 
under $\Aut(\Gamma)$. A relation 
$R$ has a first-order definition in an $\omega$-categorical
structure $\Gamma$ if and only if $R$ is preserved by all 
automorphisms of $\Gamma$.
\end{theorem}

This theorem implies that a structure $({\mathbb L}; R_1,R_2,\dots)$ is 
first-order definable in $(\mL;C)$
if and only if its automorphism group 
contains the automorphisms of $(\L;C)$. 

Automorphism groups $G$ of relational structures
carry a natural topology, namely the topology of \emph{pointwise convergence}.  
Whenever we refer to topological properties of groups it will be with respect to this topology.
To define this topology, we begin by giving 
the domain $X$ of the relational structure the discrete topology. We then
view $G$ as a subspace of the Baire space $X^X$ which carries the 
product topology; see e.g.\ Cameron~\cite{Oligo}. A set of permutations is called \emph{closed} if it is closed in the subspace ${\rm Sym}(X)$ of $X^X$, where ${\rm Sym}(X)$ is the set of all bijections from $X$ to $X$. 
The \emph{closure} of
a set of permutations $P$ is the 
smallest closed set of permutations that contains $P$
and it will be denoted by $\bar P$. 
Note that $\bar P$ equals the set of all permutations $f$ 
 such that for every finite subset $A$ of the domain there is a $g \in P$ such that $f(a)=g(a)$ for all $a \in A$.

We write $\left< P \right>$ for the smallest permutation group that contains a given set of permutations $P$.
Note that the smallest closed permutation group that 
contains a set of permutations $P$ equals $\overline{\left< P \right>}$. 
It is easy to see that a set of permutations $G$ on a set $X$ is a closed subgroup of the group of all permutations of $X$
if and only if $G$ is the automorphism group of a relational structure~\cite{Oligo}. 

We need some terminology from permutation group theory and we mostly follow
 Bhattacharjee, Macpherson, M\"oller and Neumann~\cite{Infinite-perm-groups}.
A permutation group $G$ on a set $X$ is called 
\begin{itemize}
\item \emph{$k$-transitive}
if for any two sequences $a_1,\dots,a_k$ and $b_1,\dots,b_k$ of $k$ distinct points of $X$ there exists $g$ in $G$ such that $g(a_i)=b_i$ for all $1 \leq i \leq k$,
\item \emph{transitive}
if $G$ is $1$-transitive,
\item \emph{highly transitive} if it is $k$-transitive for all natural numbers $k$,
\item \emph{primitive} if it is transitive and all equivalence relations that are preserved
by all operations in $G$ are either the equivalence relation with one equivalence class
or the equivalence relation with equivalence classes of size one. 
\end{itemize}

The following simple fact illustrates the link between model theoretic and permutation group theoretic concepts.

\begin{proposition}\label{prop:equality}
For an automorphism group $G$ of a relational structure $\Gamma$ with domain $D$, the following are equivalent. 
\begin{itemize}
\item $G$ is highly transitive.
\item $G$ equals the set of all permutations of $D$.
\item $\Gamma$ is a reduct of $(D;=)$. 
\end{itemize}
\end{proposition}

The \emph{pointwise stabilizer} at $Y \subset X$ of a permutation group $G$ on $X$ is the permutation group on $X$ consisting of all permutations $\alpha \in G$ such that $\alpha(y)=y$ for all $y \in Y$.
A subset $X'$ of $X$ is said to be a \emph{Jordan set} (for $G$ in $X$) if $|X'| > 1$ and the pointwise stabilizer
$H$ of $G$ at $X \setminus X'$ is transitive on $X'$.

If the group $G$ is $(k+1)$-transitive and $X'$ 
is any co-finite subset with $|X \setminus X'| = k$, then $X'$ is automatically a Jordan set.
Such Jordan sets will be said to be \emph{improper} while all other will be called \emph{proper}. 
We say that the Jordan set $X'$ is \emph{$k$-transitive} 
if the pointwise stabilizer $H$ is $k$-transitive on $X'$. 
The permutation group $G$ on the set $X$ is said to be a \emph{Jordan group} if $G$ is transitive on $X$ and there exists a proper Jordan set for $G$ in $X$. 
The main result that we will use in Section~\ref{sect:group-class}
is the following.

\begin{theorem}[Note 7.1 in Adeleke and Neumann~\cite{AdelekeNeumannPrimitive}]
\label{thm:jg1}
If $G$ is primitive and has 2-transitive proper Jordan sets, then
$G$ is either highly transitive or it preserves a C- or D-relation on $X$.
\end{theorem}

Note that $\Aut(\mathbb L;C)$ is 2-transitive by homogeneity and
that 2-transitivity implies primitivity. The following proposition
shows that Theorem~\ref{thm:jg1} applies in our setting. 

%

\begin{proposition}\label{prop:Jordan-sets}
For two arbitrary distinct elements $a,b \in \mL$,
the set $S := \{x \in \mL : ax|b \}$ is a 2-transitive proper primitive
Jordan set of $\Aut(\L;C)$.
\end{proposition}
\begin{proof}
The pointwise stabilizer of $\Aut(\L;C)$ at $\L \setminus S$ acts 2-transitively on $S$; this can be shown via a simple back-and-forth argument.  
\end{proof}

\subsection{The rerooting lemma}
\label{sect:rerooting}
We will now prove some fundamental
lemmata concerning
functions that preserve $Q$. 
They will be needed to prove Theorem~\ref{thm:reducts} which is the main result of Section~\ref{sect:autos}.
They will also be used in subsequent sections: we emphasize that these results are not restricted
to permutations.
The most important lemma is the \emph{rerooting lemma} (Lemma~\ref{lem:rer}) 
about functions that preserve $Q$ and violate $C$. The following notation 
will be convenient in the following.
\begin{definition}
We write $x_1\dots x_n:y_1 \dots y_m$
if $Q(x_ix_j,y_ky_l)$ for all $i,j \leq n$ and $k,l \leq m$. 
\end{definition}

\begin{lemma}\label{lem:q-to-c}
Let $A_1,A_2 \subseteq \mL$ 
be such that $A_1 | A_2$ and let $f \colon A_1 \cup A_2 \to \mL$ be
a function that preserves $Q$ and satisfies $f(A_1)|f(A_2)$.
Then $f$ also preserves $C$.
\end{lemma}
\begin{proof}
Since $A_1|A_2$, we have $A_1\cup A_2\geq 2$. Clearly, the claim of lemma holds if $|A_1\cup A_2|=2$. It remains to consider the case $|A_1\cup A_2|\geq 3$. Let $a_1,a_2,a_3 \in A_1 \cup A_2$ be three distinct elements such that $a_1a_2|a_3$. We have to verify that
$f(a_1)f(a_2)|f(a_3)$ and we do this by considering four different cases.
 \begin{itemize}
\item $a_1,a_2 \in A_1$ and $a_3 \in A_2$.
In this case, since $f(A_1)|f(A_2)$, we have in particular that
$f(a_1)f(a_2)|f(a_3)$.
 \item $a_1,a_2 \in A_2$ and $a_3 \in A_1$. Analogous to the previous case. 
 \item $a_1,a_2,a_3 \in A_1$. Let $b \in A_2$. Clearly $a_1a_2:a_3b$, and $f(a_1)f(a_2):f(a_3)f(b)$ since $f$ preserves $Q$. 
Moreover, we have $f(a_1)f(a_2)f(a_3)|f(b)$, and thus $f(a_1)f(a_2)|f(a_3)$. 
 \item $a_1,a_2,a_3 \in A_2$. Analogous to the previous case. 
 \end{itemize}
Since we have assumed that $A_1|A_2$, these cases are in fact exhaustive.
One may, for instance, note that if $a_1,a_3 \in A_1$ and $a_2 \in A_2$,
then $a_1a_3|a_2$ which immediately contradicts that
$a_1a_2|a_3$. 
 \end{proof}

\begin{lemma}\label{lem:invq}
Let $A \subset \mL$ be finite of size at least two and let $f \colon A \to \mL$
be a function which preserves $Q$. 
Then there exists a non-empty $B \subsetneq A$ such that the following conditions hold:
\begin{itemize}
\item $f(B)|f(A \setminus B)$
\item $B|x$ for all $x \in A \setminus B$. 
\end{itemize}
\end{lemma}
\begin{proof}
Let $B_1,B_2$ be non-empty such that $B_1 \cup B_2 = A$ and 
$f(B_1)|f(B_2)$. We see that $B_1,B_2$ is a partitioning of $A$
since $f(B_1)|f(B_2)$ implies $B_1 \cap B_2 = \emptyset$.
If $B_1|x$ for all $x \in B_2$, then we can choose $B = B_1$ and we are done. Otherwise there are $u,v \in B_1$
and $w \in B_2$ such that $u|vw$. We claim that in this case $x|B_2$ for all $x \in B_1$. Since $f$ preserves $Q$ on 
$A$ and $f(u) f(v):f(w) f(x)$ holds for every $x \in B_2$, we have $uv:wx$ by Lemma~\ref{lem:pres-q}. Therefore $u|wx$ and $v|wx$ hold. 
This implies that $u|B_2$ holds. Let $w',w''$ be two
arbitrary elements from $B_2$ and $u'$ an arbitrary element from $B_1$. We thus have
$f(w')f(w''):f(u')f(u)$ and, once again by Lemma~\ref{lem:pres-q}, we have $uu':w'w''$. This implies
$u|w'w''$ and consequently $u'|w'w''$. Hence,
$u'|B_2$ for arbitrary $u' \in B_2$. 
\end{proof}

We will now introduce the idea of $c$-{\em universality}. 
This seemingly simple concept is highly important throughout the article
and it will be encountered in several different contexts.

\begin{definition}
Arbitrarily choose $c \in \mL$.
A set $A \subseteq \mL \setminus \{c\}$ is called \emph{$c$-universal}
if for every finite $U \subset \mL$ and for every $u \in U$,
there exists an $\alpha \in \Aut(\mL;C)$ such that 
$\alpha(u)=c$ and $\alpha(U) \subseteq A \cup \{c\}$. 
\end{definition}

We continue by presenting the rerooting lemma which identifies permutations $g$ of $\mL$ 
that preserve $Q$ and 
can be used for generating all automorphisms of $(\mL;Q)$ when combined with $\Aut(\mL;C)$. The idea is based 
on the
following observation: the finite substructures of $(\mL;C)$ provide information about the root of the 
underlying tree whereas the
finite substructures of $(\mL;Q)$ only provide information about the underlying unrooted trees. 
Intuitively, we use the function $g$ to change the position of the root
in order to generate all automorphisms of $(\mL;Q)$. 

\begin{lemma}[Rerooting Lemma]\label{lem:rer}
Arbitrarily choose $c \in \mL$ and assume that $A \subseteq \mL \setminus \{c\}$
is $c$-universal. 
If $g$ is a permutation of $\mL$ that preserves $Q$ on $A \cup \{c\}$
and satisfies $g(A)|g(c)$, then
$$\Aut(\mL;Q) \subseteq \overline{ \left< \Aut(\mL;C) \cup \{g\} \right>}.$$ 
\end{lemma}

\begin{proof}
Arbitrarily choose $f \in \Aut(\mL;Q)$
and let $X$ be an arbitrary finite subset of $\mL$. 
We have to show that $\left < \Aut(\mL;C) \cup \{g\} \right >$ contains
an operation $e$ such that $e(x) = f(x)$ for all $x \in X$. 
This is trivial when $|X|=1$ so we assume that $|X| \geq 2$. 
By Lemma~\ref{lem:invq}, there exists a non-empty proper subset $Y$ of $X$ such that 
$f(Y) | f(X \setminus Y)$ and $Y|x$ for all $x \in X \setminus Y$. 
By the homogeneity of $(\mL;C)$, we can choose an element $c' \in \mL \setminus X$ such that $c'|Y$ and
$(Y \cup \{c'\})|x$ for all $x \in X \setminus Y$. 
By $c$-universality, there exists an 
$\alpha \in \Aut(\mL;C)$ such that $\alpha(X \cup \{c'\}) \subseteq A \cup \{c\}$ and $\alpha(c') = c$. Let $h := g \circ \alpha$. Note that $h$ preserves $Q$ on $X$ and that $h$ is a permutation. 
We continue by proving a particular property of $h$.

\medskip

{\bf Claim.} $h(Y)|h(X \setminus Y)$. 

To prove this, we first show that $h(y_1)h(y_2)|h(y_3)$ 
for every $y_1, y_2 \in Y$ and $y_3 \in X \setminus Y$. 
We have $y_1y_2:y_3c'$ by the choice of $c'$ and this implies that 
$h(y_1)h(y_2):h(y_3)h(c')$. Since $\alpha(X)\subseteq A$, it follows from the definition of $h$ that $h(y_1),h(y_2),h(y_3)\in g(A)$. Since $g(c)|g(A)$ and $\alpha(y_i)\in A$ for every $i\in \{1,2,3\}$, we have $g(c)|h(y_1)h(y_2)h(y_3)$. Since $h(c')=g(c)$ and $h(y_1)h(y_2):h(y_3)h(c')$, it follows that $h(y_1)h(y_2)|h(y_3)$. 
In the same vein, we show that $h(y_1)|h(y_2)h(y_3)$ for every $y_1 \in Y$ and $y_2, y_3 \in X \setminus Y$. 
In this case, we have $y_1c' : y_2y_3$ by the choice of $c'$ and this implies $h(y_1)h(c'):h(y_2)h(y_3)$. 
Since $h(c') = g(c)$ and $g(c)|h(y_1)h(y_2)h(y_3)$, we see that $h(y_1)|h(y_2)h(y_3)$. 

\medskip

Let $\beta \colon h(X) \to f(X)$ be defined by
$\beta(x)=f(h^{-1}(x))$. Note that $h^{-1}$ is well-defined
since $h$ is an injective function.
Since both $h$ and $f$ preserve
$Q$, we have that $\beta$
preserves $Q$ by Lemma~\ref{lem:pres-q}. 

Note that $\beta(h(Y))|\beta(h(X \setminus Y))$ since $\beta(h(x))=f(x)$
and we have assumed that $f(Y)|f(X \setminus Y)$.
Hence, the conditions of Lemma~\ref{lem:q-to-c}
apply to $\beta$ for $A_1 := h(Y)$ and $A_2 := h(X \setminus Y)$ if we use the claim above. 
It follows that $\beta$ preserves $C$. 
By the homogeneity of $(\mL;C)$, there exists an $\gamma \in \Aut(\mL;C)$ that extends $\beta$.
Then $e := \gamma \circ h$ has the desired property.
\end{proof}
Observe the following important consequence of
Lemma~\ref{lem:rer}. 

\begin{corollary}\label{cor:q}
Assume $f \in \Aut(\mL;Q)$ violates $C$. 
Then $$\overline{\left< \Aut(\mL;C) \cup \{f\} \right>} = 
\Aut(\mL;Q) \; .$$
\end{corollary}
\begin{proof}
The relation $Q$ is first-order definable over $(\mL;C)$ so $\Aut(\mL;C) \subseteq \Aut(\mL;Q)$.
Furthermore, $f$ preserves $Q$ and it follows that $$\overline{\left< \Aut(\mL;C) \cup \{f\} \right>} \subseteq \Aut(\mL;Q) \; .$$

For the converse, choose $f \in \Aut(\mL;Q)$ 
such that there are $a_1,a_2,a_3 \in \mL$
with $a_1|a_2a_3$ and $f(a_1)f(a_2)|f(a_3)$. 
Let $A = \{ x \; | \; xa_1:a_2a_3\}$. We will show that $f(A)|f(a_3)$. Let $x,y\in A$ be arbitrary. Since $f$ preserves $Q$, we have $f(x)f(a_1):f(a_2)f(a_3)$ and $f(y)f(a_1):f(a_2)f(a_3)$. It follows from the condition $f(a_1)f(a_2)|f(a_3)$ that $$f(x)f(a_1)|f(a_2)\wedge f(x)f(a_1)|f(a_3)\wedge f(y)f(a_1)|f(a_2)\wedge f(y)f(a_1)|f(a_3).$$ Since $f(x)f(a_1)|f(a_3)\wedge f(y)f(a_1)|f(a_3)$, we have $f(x)f(y)|f(a_3)$. Thus $f(A)|f(a_3)$.

Clearly, $A$ is $a_3$-universal. Applying Lemma \ref{lem:rer} to $c=a_3$ we have $$\overline{\Aut(\mL;Q)\subseteq \langle\Aut(\mL;C)\cup \{f\}\rangle}.$$ 
\end{proof}

\subsection{Automorphism group classification}
\label{sect:group-class}
We are now ready to prove the main result concerning automorphism groups of the reducts of $(\mL;C)$.

\begin{theorem}\label{thm:reducts}
Let $G$ be a closed permutation group on the set $\mL$ that contains 
$\Aut(\mL;C)$. Then $G$ is
either $\Aut(\mL;C)$, $\Aut(\mL;Q)$, or $\Aut(\mL;=)$.
\end{theorem}
\begin{proof}
Because $G$ satisfies the conditions of Theorem~\ref{thm:jg1}, it is either highly transitive
or it preserves a C- or D-relation. 
If $G$ is highly transitive, then $G$ equals $\Aut(\mL;=)$
by Proposition~\ref{prop:equality}. 
Assume instead that $G$ preserves a C-relation $C'$. 
We begin by making an observation.

\medskip

{\bf Claim 0.} All tuples $(o,p,q) \in C'$ with pairwise distinct entries satisfy $o|pq$. 

Suppose for contradiction that $p|oq$. 
Then, $(o,p,q)$ is in the same orbit as
$(q,p,o)$ in $\Aut(\mL;C)$ and therefore also in $G$. Since $C'$ is preserved by $G$, 
we have $C'(q,p,o)$ which contradicts C2. 
Similarly, it is impossible that $q|op$. 
Thus, the only remaining possibility is $o|pq$
since $C$ satisfies C8.

\medskip

Arbitrarily choose $a,b,c \in \mL$ such that $a|bc$ 
and some $\alpha \in G$.
If $|\{a,b,c\}| = 2$, then (by 2-transitivity of $\Aut(\mL;C)$)
we have that $(\alpha(a),\alpha(b),\alpha(c))$
is in the same orbit as $(a,b,c)$ of $\Aut(\mL;C)$. Consequently,
$(\alpha(a),\alpha(b),\alpha(c)) \in C$. 
Suppose instead that $|\{a,b,c\}| = 3$. 
Observe that $C'$ contains a triple with pairwise
distinct entries. Arbitrarily choose two distinct elements $u,v \in \mL$. 
Axiom C6 implies that there exists a $w \in \L$ such that $C'(u,vw)$ and $w \neq v$. 
In fact, we also have $w \neq u$ since otherwise $C'(u,vu)$ which is impossible
due to C2 and C4. 
In particular, it follows that $u|vw$ and therefore
$(u,v,w)$ is in the same orbit
as $(a,b,c)$ in $\Aut(\mL;C)$. 
It follows that $(a,b,c) \in C'$. 
Since $G$ preserves $C'$ we have
$C'(\alpha(a),\alpha(b)\alpha(c))$. 
By Claim~0, 
$\alpha(a)|\alpha(b)\alpha(c)$. 
We conclude that $\alpha$ preserves $C$ and that $G = \Aut(\mL;C)$.


\medskip

Finally, we consider the case when $G$ preserves a D-relation $D$. 
We begin by making three intermediate observations.


\medskip

{\bf Claim 1.} 
Every tuple $(a,b,c,d) \in D$ with pairwise distinct entries satisfies $ab:cd$. 


Suppose for contradiction that $ac:bd$. 
Then either $ac|b \wedge ac|d$ or 
$a|bd \wedge c|bd$ by the definition of relation $Q$. 
In the first case, $(a,b,c,d)$ is in the same orbit
as $(c,b,a,d)$ in $\Aut(\mL;C)$ so $(c,b,a,d) \in D$. Axiom D1 implies that $(a,d,b,c) \in D$
and this contradicts D2. If $a|bd \wedge c|bd$, then we can obtain a contradiction in
a similar way.
Finally, the case when $ad:bc$ can be treated analogously.
It follows that $ab:cd$ since $Q$ satisfies D7.

\medskip

{\bf Claim 2.} 
$D$ contains a tuple $(o,p,q,r)$
with pairwise distinct entries 
such that $op|qr$ holds. 

Let $u,v,w \in \mL$ be three distinct elements such that $uv|w$. There is an 
$x \in \mL \setminus \{u,v,w\}$ such that $D(uv,wx)$ by D5.
Claim~1 immediately implies that $uv:wx$. 
We consider the following cases.

\begin{itemize}
\item $uv|wx$. There is nothing to prove in this case.
\item $uvw|x$. Choose $y \in \mL$ be such that $y \neq w$ and $uv|yw$. It follows from D3 that $D(uv,yw)$
or $D(yv,wx)$. The second case is impossible since
$yv:wx$ does not hold. We see that $(u,v,y,w) \in D$
and we are done.
\item $uv|x$ and $uvx|w$. One may argue similarly as in the previous case by
choosing $y \in \mL \setminus \{u,v,w,x\}$
such that $uv|yx$ and observe that $(u,v,x,w) \in D$ by D1.
\end{itemize}

\medskip

{\bf Claim 3. } $D$ contains a tuple $(a,b,c,d)$
with pairwise distinct entries such that
$ab|c \wedge abc|d$.  

It follows from Claim~2
that there exists a tuple $(o,p,q,r)$ with pairwise distinct
entries such that $op|qr$ holds.
Choose $s \in \mL$ such that $opqr|s$ holds.
 Axiom D3 implies that
$D(sp,qr)$ or $D(op,qs)$. We are done if the second case holds. If the first case holds, 
then we have $D(qr,ps)$ by D1 and we are once again done. 

\medskip

Now, we show that every $f \in G$ preserves $Q$.
Arbitrarily choose $a_1,a_2,a_3,a_4 \in \mL$ such that $a_1a_2:a_3a_4$.
We show that $(a_1,a_2,a_3,a_4) \in D$ (and, consequently, that
$(f(a_1),f(a_2),f(a_3),f(a_4)) \in D$) by an exhaustive case analysis.
Claim~2 implies that $D$ contains
a tuple $(o,p,q,r)$ with pairwise distinct entries 
and $op|qr$. Consequently, $D$ contains all tuples
in the same orbit as $(o,p,q,r)$ in $\Aut(\mL;C)$. 

If $a_1,a_2,a_3,a_4$ are pairwise distinct
and satisfy $a_1a_2|a_3a_4$, then $(a_1,a_2,a_3,a_4) \in D$ by Claim~1.
Similarly, if $a_1,a_2,a_3,a_4$ are pairwise distinct
and satisfy $a_1a_2|a_3$ and $a_1a_2a_3|a_4$, then 
Claim~3 implies that $(a_1,a_2,a_3,a_4) \in D$.
If
$a_2|a_3a_4$ and $a_1|a_2a_3a_4$, then 
$(a_1,a_2,a_3,a_4) \in D$ by D1.
If 
$a_1a_2a_4|a_3$ and $a_1a_2|a_4$, 
or if $a_1|a_3a_4$ and $a_2|a_1a_3a_4$, then
$(a_1,a_2,a_3,a_4) \in D$ by D1.
If $a_3=a_4$, $a_1 \neq a_3$, $a_2 \neq a_3$, 
then $(a_1,a_2,a_3,a_4) \in D$ by D4. 
The only remaining possibility to satisfy $a_1a_2:a_3a_4$ is that $a_1=a_2$, $a_3 \neq a_1$, $a_4 \neq a_1$. In this case, $(a_1,a_2,a_3,a_4) \in D$ by D4 and D1. 
Hence, in all cases we have $(a_1,a_2,a_3,a_4) \in D$ and, consequently,
$(f(a_1)f(a_2),f(a_3)f(a_4)) \in D$.

We can now conclude this part of the proof.
If $f(a_1),f(a_2),f(a_3),f(a_4)$
are pairwise distinct, then $f(a_1)f(a_2):f(a_3)f(a_4)$ by Claim~1. Otherwise, 
one of the following cases hold:
\begin{itemize}
\item $f(a_1)=f(a_2)$, $f(a_3) \neq f(a_1)$, and 
$f(a_4) \neq f(a_1)$, 
\item $f(a_3)=f(a_4)$, $f(a_1) \neq f(a_3)$, and $f(a_1) \neq f(a_4)$, or
\item $f(a_1)=f(a_2)$ and $f(a_3)=f(a_4)$.
\end{itemize}
In all three cases, we have that $f(a_1)f(a_2):f(a_3)f(a_4)$ and 
$G \subseteq \Aut(\mL;Q)$. If $G = \Aut(\mL;C)$, then
we are done. Otherwise, pick one $f \in G \setminus \Aut(\mL;C)$. 
Corollary~\ref{cor:q} asserts that 
$$\overline{\left < \Aut(\mL;C) \cup \{f\} \right >} = \Aut(\mL;Q) \subseteq G,$$ and it
follows that $G = \Aut(\mL;Q)$. 
\end{proof}

The following is an immediate consequence of Theorem~\ref{thm:reducts} in combination with 
Theorem~\ref{thm:ryll}. 
 
 \begin{corollary}\label{cor:fo-reducts} 
 Let $\Gamma$ be a reduct of $(\mL;C)$. Then
 $\Gamma$ is first-order interdefinable
 with $(\mL;C)$, $(\mL;Q)$, or $(\mL;=)$. 
 \end{corollary}
 \begin{proof}
Since $\Gamma$ is a reduct of $(\mL;C)$, $\Aut(\Gamma)$ is a closed group that contains $\Aut(\mL;C)$ and therefore equals $\Aut(\mL;C)$, $\Aut(\mL;Q)$, or $\Aut(\mL;=)$ by Theorem~\ref{thm:reducts}. Theorem~\ref{thm:ryll} implies that
 $\Gamma$ is first-order interdefinable with 
 $(\mL;C)$, with $(\mL;Q)$, or with $(\mL;=)$. 
 \end{proof}
  
Corollary~\ref{cor:fo-reducts} will be refined to a classification
up to existential interdefinability in the forthcoming sections.

%% file: ramsey.tex
\section{Ramsey theory for the C-relation}\label{sect:tramsey}
To analyze endomorphism monoids of reducts of $(\mL;C)$, 
we apply Ramsey theory; 
a survey on this technique can be found in Bodirsky \& Pinsker~\cite{BP-reductsRamsey}. 
The basics of the Ramsey approach are presented in
Section~\ref{sect:ramseyclass} 
and we introduce the important concepts of {\em canonicity}
and the {\em ordering property} in Sections~\ref{sect:behavior}
and~\ref{sect:orderingproperty}, respectively. We would like to mention that none of the results from the previous sections that use the theory of Jordan permutation groups is needed in the subsequent parts.

We will frequently use topological methods 
when studying transformation monoids.
The definition of the topology of pointwise convergence
 for transformations monoids is analogous
to the definition for groups:
the \emph{closure $\overline F$} of $F \subseteq \mL^{\mL}$ 
is the set of all
functions $f \in \mL^{\mL}$ with the property that for
every finite subset $A$ of $\mL$, there is a $g \in F$
such that $f(a)=g(a)$ for all $a \in A$.
A set of functions is \emph{closed} if $F = \overline F$. We write $\left < F \right >$ for 
the smallest transformation monoid that contains $F$.
The smallest closed transformation monoid
that contains a set of functions $F$ equals $\overline{\left< F \right >}$. The closed transformation monoids are precisely those
 that are endomorphism monoids of relational structures.
We say that a function $f$ is \emph{generated} by a set of operations $F$ is 
$f$ is in the smallest closed monoid that contains $F$. 
A more detailed introduction to these concepts can be found in Bodirsky~\cite{Bodirsky-HDR}.

\subsection{Ramsey classes}
\label{sect:ramseyclass}
 Let $\Gamma, \Delta$ be finite $\tau$-structures. We write ${\Delta \choose \Gamma}$ for the set of all substructures of $\Delta$ that are isomorphic to $\Gamma$. When $\Gamma, \Delta, \Theta$ are $\tau$-structures, then we write
$\Theta \to (\Delta)^{\Gamma}_r$
if for all functions $\chi \colon {\Theta \choose \Gamma} \to \{1,\dots,r\}$ there exists 
$\Delta' \in {\Theta \choose \Delta}$  such that $\chi$ is constant on ${\Delta' \choose \Gamma}$.

\begin{definition}
\label{def:Ramseyclass}
A class of finite relational structures $\cal C$ that is closed under isomorphisms and substructures is called
\emph{Ramsey} if for all $\Gamma, \Delta \in \cal C$ and arbitrary $k \geq 1$, 
there exists a $\Theta \in \cal C$
such that $\Delta$ embeds into $\Theta$ and 
$\Theta \to (\Delta)^{\Gamma}_k$.
\end{definition}
 
A homogeneous structure $\Gamma$ is called \emph{Ramsey} if the class
 of all finite 
structures that embed into $\Gamma$ is Ramsey. We refer the reader to
Kechris, Pestov and Todorcevic~\cite{Topo-Dynamics} or \Nesetril~\cite{RamseyClasses} for
more information about 
the links between Ramsey theory and homogeneous structures. 
An example of a Ramsey structure is $(D;=)$---the fact that the 
class of all finite structures that embed
into $(D;=)$ is Ramsey can be seen as a reformulation of Ramsey's 
classical result~\cite{Ram:On-a-problem}.

The Ramsey result that is relevant 
in our context (Theorem~\ref{thm:tramsey}) is a consequence of a 
more powerful theorem due to Miliken~\cite{Mil79}. 
The theorem in the form presented below and a direct proof
of it can be found in Bodirsky \& Piguet~\cite{BodirskyPiguet}. We mention that a weaker version of this theorem (which was shown by the academic grand-father of the first author of this article~\cite{DeuberTreeRamsey}) has been known for a long time.

\begin{theorem}[see Bodirsky~\cite{BodirskyRamsey} or Miliken~\cite{Mil79}]
\label{thm:tramsey}
The structure $(\mL;C,\prec)$ is Ramsey.
\end{theorem}

We also need the following result. 

\begin{theorem}[see Bodirsky, Pinsker \& Tsankov~\cite{BPT-decidability-of-definability}]
\label{thm:cramsey}
If $\Gamma$ is homogeneous and Ramsey, then every expansion
of $\Gamma$ by finitely many constants is Ramsey, too.
\end{theorem}

\subsection{Canonical functions}
\label{sect:behavior}
The typical usage of Ramsey theory in this article is for showing
that the endomorphisms
of $\Gamma$ behave {\em canonically} on large parts of the 
domain; this will be formalized below.  
A wider introduction to canonical operations can be found 
in Bodirsky~\cite{Bodirsky-HDR} and
Bodirsky \& Pinsker
~\cite{BP-reductsRamsey}.
The definition of canonical functions given below
is slightly different from the one given in~\cite{Bodirsky-HDR} and~\cite{BP-reductsRamsey}. It is easy to see that they are equivalent, though.

\begin{definition}
\label{def:canonical}
Let $\Gamma,\Delta$ be structures
and let $S$ be a subset of the domain $D$ of $\Gamma$. A function 
$f \colon \Gamma \to \Delta$ 
is \emph{canonical on $S$} as a function from $\Gamma$ to
$\Delta$ if for all $s_1,\dots,s_n \in S$ and all $\alpha \in \Aut(\Gamma)$,
there exists a $\beta \in \Aut(\Delta)$
such that $f(\alpha(s_i)) = \beta(f(s_i))$ for all $i \in \{1,\dots,n\}$. 
\end{definition}

In Definition~\ref{def:canonical}, we might omit the set $S$
if $S=D$ is clear from the context. 
Note that a function $f$ from $\Gamma$ to $\Delta$ is
canonical if and only if 
for every $k \geq 1$ and every $t \in D^k$,
the orbit of $f(t)$ in $\Aut(\Delta)$ only depends on the orbit of $t$ in $\Aut(\Gamma)$. 

\begin{example}\label{ex:minus}
Write $x \succ y$ if $y \prec x$.
The structure $(\mL;C,\succ)$ is isomorphic to 
$(\mL;C,\prec)$; 
let $-$ be such an isomorphism. 
Note that $-$ is canonical as a function from
$(\mL;C,\prec)$ to $(\mL;C,\prec)$.
\end{example}

When $\Gamma$ is Ramsey, then the following theorem allows us to work with canonical endomorphisms of $\Gamma$. 
It can be shown with the same proof as presented in Bodirsky, Pinsker, \& Tsankov~\cite{BPT-decidability-of-definability}.
\begin{theorem}
\label{thm:can}
Let $\Gamma,\Delta$ denote finite relational structures
such that $\Gamma$ is homogeneous and Ramsey while
$\Delta$ is $\omega$-categorical. 
Arbitrarily choose a function $f \colon \Gamma \to \Delta$. 
Then, 
there exists a function
$$g \in \overline{ \{ \alpha_1 f \alpha_2 :  \alpha_1 \in \Aut(\Delta), \alpha_2 \in \Aut(\Gamma) \}}$$ that is canonical as a function from 
$\Gamma$ to $\Delta$. 
\end{theorem}

Note that expansions of homogeneous structures with constant symbols are again homogeneous.
We obtain the following by combining
the previous theorem and Theorem~\ref{thm:cramsey}.

\begin{corollary}
\label{cor:const-can}
Let $\Gamma,\Delta$ denote finite relational structures
such that $\Gamma$ is homogeneous and Ramsey while
$\Delta$ is $\omega$-categorical. 
Arbitrarily choose a function $f \colon \Gamma \to \Delta$
and elements $c_1,\dots,c_n$ of $\Gamma$. Then, there exists a function
$$g \in \overline{ \{ \alpha_1 f \alpha_2 :  \alpha_1 \in \Aut(\Delta), \alpha_2 \in \Aut(\Gamma,c_1,\dots,c_n) \}}$$ 
that is canonical as a function from 
$(\Gamma,c_1,\dots,c_n)$ to $\Delta$. 
\end{corollary}



\subsection{The ordering property}
\label{sect:orderingproperty}
Another important concept from Ramsey theory
that we will exploit in the forthcoming proofs is the
\emph{ordering property}. We will next prove that the class of ordered leaf structure has this property. 

\begin{definition}[See Kechris, Pestov \& Todorcevic~\cite{Topo-Dynamics} or \Nesetril~\cite{RamseyClasses}]
Let $\mathcal C'$ be a class of finite structures over the signature $\tau \cup \{\prec\}$, where $\prec$ denotes a linear order, and let $\mathcal C$ be the class
of all $\tau$-reducts of structures from $\mathcal C'$. Then $\mathcal C'$ has the \emph{ordering property} if for every $\Delta_1 \in \mathcal C$
there exists a $\Delta_2 \in \mathcal C$ such that
for all expansions $\Delta_1' \in \mathcal C'$ of $\Delta_1$
and $\Delta_2'\in \mathcal C'$ of $\Delta_2$ 
there exists an embedding of $\Delta_1'$ into $\Delta_2'$.
\end{definition}

\begin{proposition}\label{prop:ordering}
Let $\Gamma$ be a
homogeneous relational
$\tau$-structure with domain $D$ and suppose that $\Gamma$ has an $\omega$-categorical 
homogeneous expansion $\Gamma'$ with
signature $\tau \cup \{\prec\}$ where $\prec$
denotes a linear order. 
Then, the following are equivalent.
\begin{itemize}
\item the class $\mathcal C'$ of finite structures that embed into $\Gamma'$ has the
ordering property and
\item for every finite
$X \subseteq D$ there exists a finite $Y \subseteq D$ such that for every $\beta \in \Aut(\Gamma)$ there exists an $\alpha \in \Aut(\Gamma')$ such that $\alpha(X) \subseteq \beta(Y)$. 
\end{itemize}
\end{proposition}
\begin{proof}
First suppose that $\mathcal C'$ has the ordering property and let $X \subseteq D$ be finite.
Let $\Delta_1$ be the structure induced by $X$
in $\Gamma$. 
Then, there exists $\Delta_2 \in \mathcal C$ such that for
for all expansions $\Delta_1' \in \mathcal C'$ of $\Delta_1$
and for all expansions $\Delta_2' \in \mathcal C'$ of $\Delta_2$,
there exists an embedding of $\Delta_1'$ into $\Delta_2'$. Since
every structure in $\mathcal C'$ can be embedded into
$\Gamma'$, we may
assume that $\Delta_2'$ is a substructure of $\Gamma'$ with domain $Y$.
Arbitrarily choose $\beta \in \Aut(\Gamma)$. Then, there exists an embedding from the structure
induced by $X$ in $\Gamma'$
to the structure induced by $\beta(Y)$ 
in $\Gamma'$. By homogeneity of $\Gamma'$,
this embedding can be extended to an automorphism $\alpha$ of $\Gamma'$ which has the desired property. 

\medskip

\noindent
For the converse direction, 
let $\Delta_1$ be the $\tau$-reduct of an arbitrary structure from $\mathcal C'$
and let $n$ denote the cardinality of $\Delta_1$.
Since $\Gamma'$ is $\omega$-categorical, 
there is a finite number $m$ of orbits of $n$-tuples.
Hence, there exists a set $Z$ of cardinality $n \cdot m$
such that for every 
embedding $e$ of $\Delta_1$ into $\Gamma$,
 there exists an automorphism $\alpha$ of $\Gamma'$
such that the image of $\alpha \circ e$ is a subset
of $Z$. 
By assumption, there exists a set $Y \subseteq D$ such that for every 
$\beta \in \Aut(\Gamma)$, 
there exists an $\alpha \in \Aut(\Gamma')$ 
and $\alpha(Z) \subseteq \beta(Y)$. 
Let $\Delta_2$ be the structure induced by $Y$ in
$\Gamma$. 
Now, let $\Delta_1' = (\Delta_1,\prec)$ and arbitrarily choose 
$\Delta_2' = (\Delta_2,\prec) \in \mathcal C'$.
By the choice of $Z$, there is an
embedding $f$ of $\Delta_1'$ into the substructure induced by $Z$ in $\Gamma'$. 
Since $\Gamma'$ embeds all structures from $\mathcal C'$, we can assume
that $\Delta_2'$ is a substructure of $\Gamma'$. 
By homogeneity of $\Gamma$, there is
a $\beta \in \Aut(\Gamma)$ that maps $\Delta_2$
to $\Delta_2'$. By the choice of $Y$, 
there exists an $\alpha \in \Aut(\Gamma')$ 
such that $\alpha(Z) \subseteq \beta(Y)$.
Now, $\alpha \circ f$
is an embedding of $\Delta_1'$ into $\Delta_2'$ which concludes the proof.
\end{proof}

\begin{theorem}\label{thm:ordering}
The class of all ordered leaf structures has the 
ordering property. 
\end{theorem}
\begin{proof}
By Proposition~\ref{prop:ordering}, 
it is sufficient to show that for every finite
$X \subseteq \mL$, there exists a finite $Y \subseteq \mL$ such that for 
every $\beta \in \Aut(\mL;C)$ there exists an $\alpha \in \Aut(\mL;C,\prec)$ satisfying
$\alpha(X) \subseteq \beta(Y)$. 
Let $X$ be an arbitrary finite subset of $\mL$,
let $Z=X \cup -X$ (where $-$ is defined as in Example~\ref{ex:minus}),
and let $\Delta$ be the structure induced by $Z$
in $(\mL;C,\prec)$. 
Let $\Gamma$ be the structure induced by
a two-element subset of $\mL$ in $(\mL;C,\prec)$.
The exact choice is not important since
all such structures are isomorphic. 
Since $(\mL;C,\prec)$ is Ramsey by Theorem~\ref{thm:tramsey}, 
there exists a leaf structure $\Theta$ such that
$\Theta \rightarrow (\Delta)^\Gamma_2$. 
Let $Y$ be the domain of $\Theta$.

Now, choose some $\beta \in \Aut(\mL;C)$ arbitrarily.
Define the following 2-coloring of 
${\Theta \choose \Gamma}$:
suppose that 
$x,y \in Y$ satisfy $x \prec y$.
Color the copy of $\Gamma$ induced by $\{x,y\}$ 
red iff 
$\beta(x) \prec \beta(y)$ and blue otherwise. 
Then, there exists a copy $\Delta'$ of $\Delta$ in $\Theta$ such that all 
copies of $\Gamma$ in $\Delta'$ have the same color.
If the color is red, clearly there is an automorphism 
$\alpha \in \Aut(\mL;C,\prec)$ such that $\alpha(X) \subseteq \beta(\Delta') \subseteq Y$. If 
the color is blue, then there is also an automorphism 
$\alpha \in \Aut(\mL;C,\prec)$ such that $\alpha(X) \subseteq \beta(\Delta') \subseteq Y$ since
$Z$ also contains $-X$. 
\end{proof}

\ignore{
Let $n:=|X|$. 
Since $(\mL;C,\prec)$ is $\omega$-categorical, $\Aut(\mL;C,\prec)$ has a finite number $m$ of
orbits of $n$-tuples $\bar x = (x_1,\dots,x_n)$
with the property that $x_1 \prec \cdots \prec x_n$ and $\gamma(\{x_1,\dots,x_n\})=X$ 
for some $\gamma \in \Aut(\mL;C)$. 
Let $\bar x^1,\dots,\bar x^m$ be orbit representatives for those orbits.
Let $Y$ be the set of all entries from all those
representatives.
To see that $Y$ has the required property, let $\beta \in \Aut(\mL;C)$ be arbitrary. Let $\bar y = (y_1,\dots,y_n)$ be the elements of $\beta(X)$ such that $y_1 \prec \cdots \prec y_n$; then $\bar y$ lies in the same orbit as 
$\bar x^i$ in $\Aut(\mL;C,\prec)$, for some $i \leq m$. Let $\alpha \in \Aut(\mL;C,\prec)$
be such that $\alpha(\bar x^i) = \bar y$. 
Then
$\beta(X) = \{y_1,\dots,y_n\} \subseteq \alpha(Y)$. 
and thus we obtain an automorphism $\beta$ of 
$(\mL;C,\prec)$ as required. 
}

%% file: canonical.tex
\section{Endomorphism monoids of reducts}
\label{sect:canonical}
In this section we prove the remaining results that were stated in Section~\ref{sect:results}.
We start with a description of the basic idea how to use 
the Ramsey theoretic tools introduced in the previous section.
In our proof, we can exclusively focus on analyzing
\emph{injective} endomorphisms, because of a fundamental
lemma which we describe next. 
Since $(\mL;C)$ has a 2-transitive automorphism group,
all reducts $\Gamma$ of $(\mL;C)$ 
also have a 2-transitive automorphism group. 
We can thus apply the following result.

\begin{lemma}[see, e.g., Bodirsky\cite{Bodirsky-HDR}]\label{lem:const}
Let $\Gamma$ be a relational structure with a 2-transitive
automorphism group. If $\Gamma$ has a non-injective endomorphism, then it also has a constant endomorphism. 
\end{lemma}
 
Let $\Gamma$ be a reduct of $(\mL;C)$.
Suppose that $\Gamma$ has an endomorphism $e$
that does not preserve $C$, i.e., there is $(o,p,q) \in C$
such that $(e(o),e(p),e(q)) \notin C$. 
 If $e$ is not injective, 
then $\Gamma$ also has a constant endomorphism by Lemma~\ref{lem:const}.
 In this case, the third item in Theorem~\ref{thm:endos} applies and we are done. So suppose in the following
 that $e$ is injective. 
By Theorem~\ref{thm:tramsey}, the structure 
$(\mL;C,\prec)$ is Ramsey. Hence, 
Corollary~\ref{cor:const-can} implies that $\{e\} \cup \Aut(\mL;C)$ generates 
an injective function $f$ that equals $e$ on $o,p,q$ and therefore still violates $C$, 
but is canonical as a function from
$(\mL;C,\prec,o,p,q)$ to $(\mL;C,\prec)$.

As we have noted above, 
a canonical function $f$ from $\Gamma$
to $\Delta$ induces
a function from the orbits of $k$-tuples in $\Aut(\Gamma)$ to the orbits of $k$-tuples in $\Aut(\Delta)$;
we will refer to those functions as the \emph{behavior} of $f$. 
There are finitely many behaviors of canonical injections from
$(\mL;C,\prec,o,p,q)$ 
to $(\mL;C,\prec)$: since the pre-image is homogeneous in 
a ternary language with three constants, 
their number is bounded by the number of functions from ${\mathcal O}_6 \rightarrow {\mathcal O}_3$,
where ${\mathcal O}_k$ denotes the set of 
orbits of $k$-tuples of distinct elements
in $(\mL;C,\prec)$. The function $s \colon k \mapsto |{\mathcal O}_k|$ is well-known in combinatorics (see Sloane's Integer Sequence A001813 and see \cite{PJC87} for various enumerative results for leaf structures on trees), and we have $s(n) = (2n)!/n!$. In particular, $s(3) = 12$ and $s(6) = 30240$. So the number of
canonical behaviors of functions from $(\mL;C,\prec,o,p,q)$ to $(\mL;C,\prec)$ is bounded by
$12^{30240}$.  
For every function with one of those behaviors, we prove that $\Gamma$ must be as described in item 3 and 4 of Theorem~\ref{thm:endos}.
Since $12^{30240}$ is a somewhat large number of cases,
the way we treat these cases in the following is important. 
We then repeat the same strategy 
for the structure $(\mL;Q)$ but here we have
to expand with four constants, that is,
we analyze canonical functions from 
$(\mL;C,\prec,c_1,\dots,c_4) \to (\mL;C,\prec)$. 

In the following, several arguments hold for the expansion 
of $(\mL;C,\prec)$ by any finite number of constants $\bar c = (c_1,\dots,c_n)$. 
 The following equivalence relation plays an important role. 

\begin{definition}\label{def:Ec}
Let $\bar c = (c_1,\dots,c_n) \in \mL^n$. Then $E_{\bar c}$ denotes the equivalence relation defined on $\mL \setminus \{c_1,\dots,c_n\}$ by $$E_{\bar c}(x,y) \Leftrightarrow \bigwedge_{i=1}^n xy|c_i$$ 
The equivalence classes of $E_{\bar c}$ are called
\emph{cones} (of $(\mL;C,\bar c)$). 
We write $S^{\bar c}_a$ for the cone that contains $a \in \mL \setminus \{c_1,\dots,c_n\}$.
\end{definition}
Note that each cone induces in $(\mL;C,\prec)$ a structure that
is isomorphic to $(\mL;C,\prec)$. 
\vspace{.5cm}

In Section~\ref{sect:canonical-without-constants}, Section~\ref{sect:one-constant}, 
and Section~\ref{sect:two-constants}
we study the behavior of canonical functions with zero, one, and two constants, respectively. 
Finally, in Section~\ref{sect:classification},
we put the pieces together and prove Theorem~\ref{thm:endos} and Corollary~\ref{cor:ex-reducts}.

\subsection{Canonical behavior without constants}
\label{sect:canonical-without-constants}
In this section we analyze the behavior
of canonical functions from $(\mL;C,\prec)$ to $(\mL;C,\prec)$.
In particular, we discuss
possible behaviors on cones (Corollary~\ref{cor:exhaust}) 
and close with a useful lemma 
(Lemma~\ref{lem:lin})
that shows that when a reduct $\Gamma$ of $(\mL;C)$ is preserved by functions with certain behaviors, then $\Gamma$ is homomorphically equivalent to a reduct of $(\mL;=)$. 

\begin{definition}
Let $A \subseteq \mL$ and $e \colon \mL \to \mL$
a function. Then we say that $e$ \emph{has on $A$ the behavior}
\begin{itemize}
\item[$\id$] if for all $x,y,z \in A$ with $xy|z$ we have that $e(x)e(y)|e(z)$. 
\item[$\lin$] if for all $x,y,z \in A$ with $x \prec y \prec z$ we have that $e(x)|e(y)e(z)$.
\item[$\nil$] if for all $x,y,z \in A$ with $x \prec y \prec z$ we have that $e(x)e(y)|e(z)$.
\end{itemize}
\end{definition}
In this case, 
we will also say that $e$ \emph{behaves as}
$\id$, $\lin$, or $\nil$ on $A$, respectively. 
When $f$ behaves as $\lin$ on $A = \mL$, then we do not mention $A$
and simply say that $f$ behaves as $\lin$;
we make the analogous convention for all other behaviors that we define. 
We first prove that functions with behavior $\lin$ and $\lrl$ really exist. 

\begin{lemma}\label{lem:realize}
There are
functions from $\mL \rightarrow \mL$ which preserve $\prec$ and have the 
behavior $\lin$ and $\lrl$. 
\end{lemma}
\begin{proof}
A function $f$ with behavior $\lin$ can be constructed as follows. Let $v_1,v_2,\dots$ be an enumeration of $\mL$. 
Inductively suppose that there exists a function $f \colon \{v_1,\dots,v_n\} \to \mL$ such that for all $x,y,z \in \{v_1,\dots,v_n\}$ with $x \prec y \prec z$ it holds that
$f(x)|f(y)f(z)$ and $f(x) \prec f(y) \prec f(z)$. 
This is clearly true for $n = 1$. 
We prove that $f$ has an extension $f'$
to $v_{n+1}$ with the same property. 
Let $w_1,\dots,w_{n+1}$ be such that
$\{w_1,\dots,w_{n+1}\} = \{v_1,\dots,v_{n+1}\}$
and $w_1 \prec \cdots \prec w_{n+1}$. 
We consider the following cases. 
\begin{itemize}
\item $v_{n+1} = w_1$. There exists a $c \in \mL$ such that $c|f(w_2)f(w_3)$ (see axiom C4). Note that if $n=1$, then let $c$ be such that $c\neq f(w_2)$. Pick $c$ such that $c \prec f(w_2)$, and define $f'(v_{n+1})=c$. 
\item $v_{n+1} = w_i$ for $i \in \{2,\dots,n-1\}$. 
There exists a $c \in \mL$ such that $f(w_{i-1})|cf(w_{i+1})$ and $c|f(w_{i+1})f(w_{i+2})$ (see Axiom C7). Pick $c$ such that $f(w_{i-1}) \prec c \prec f(w_{i+1})$, and define $f'(v_{n+1})=c$.
\item $v_{n+1} = w_{n}$. There exists a $c \in \mL$ such that $c \neq f(w_{n+1})$ and $f(w_{n-1})|f(w_{n+1})c$ (see Axiom C6). Pick $c$ such that $f(w_{n-1}) \prec c \prec f(w_{n+1})$, and define $f'(v_{n+1})=c$. 
\item $v_{n+1} = w_{n+1}$. There exists a $c \in \mL$ such that $c \neq f(w_{n})$ and $f(w_{n-1})|f(w_{n})c$ (see Axiom C6). Pick $c$ such that $f(w_n) \prec c$, and define $f'(v_{n+1})=c$. 
\end{itemize}
The function defined on all of $\mL$ in this way has the behavior $\lin$.
The existence of a function with behavior $\lrl$ can be shown analogously. 
\end{proof}

The functions $\lin$ and $\nil$ constructed in Lemma \ref{lem:realize} preserve the linear order $\prec$. In general, a function $f:\mL\to \mL$ with behavior $\lin$ or $\nil$ may not preserve $\prec$, however together with $\Aut(\mL;C)$ generates $\lin$ or $\nil$, respective. In the following, we will use $\lin$ and $\nil$
also to denote the functions with behavior $\lin$ and $\nil$ that have been constructed in Lemma~\ref{lem:realize}; whether we mean the behavior or the function $\lin$ and $\nil$ will always be clear from the context. As we see in the following proposition, the two functions are closely related. 

\begin{proposition}\label{prop:nil-lin}
$\Aut(\mL;C) \cup \{\nil\}$ generates $\lin$,
and $\Aut(\mL;C) \cup \{\lin\}$ generates $\nil$. 
\end{proposition}
\begin{proof}
Let $n \geq 1$ and arbitrarily choose $t \in \mL^n$.
Then $-\lrl(-t)$ and $\lin (t)$ induce isomorphic substructures
in $(\mL;C,\prec)$, and by the homogeneity of $(\mL;C,\prec)$ there is an $\alpha \in \Aut(\mL;C,\prec)$ such that $\alpha(-\lrl(-t)) = \lin (t)$.
It follows that $\lin \in \overline{\left< \Aut(\mL;C) \cup \{\lrl\} \right>}$. The fact 
that $\Aut(\mL;C) \cup \{\lin\}$ generates $\lrl$ can be shown in the same way.
\end{proof}

The following lemma classifies the behavior of canonical
injective functions from $(\mL;C,\prec)$ to $(\mL;C,\prec)$ on sufficiently 
large subsets of $\mL$. 
\begin{lemma}\label{lem:exhaust}
Let $S \subseteq \mL$ be a set that contains four elements $x, y, u, v$ such that $xy|uv$, and let $f \colon D \rightarrow D$ be injective and canonical
on $S$ as a function from $(\mL;C,\prec)$ to $(\mL;C,\prec)$.  
Then $f$ behaves as $\id$, $\lin$, or $\lrl$
on $S$. 
\end{lemma}
\begin{proof}
Since $f$ is canonical on $S$ as a function from 
$(\mL;C,\prec)$ to $(\mL;C,\prec)$,
it either preserves or reverses the order $\prec$ on $S$. We focus on the case that $f$ preserves $\prec$ on $S$, since the order-reversing case is analogous. Without loss of generality, we assume that $x\prec y\prec u\prec v$. Since $f$ preserves $\prec$, we have $f(x)\prec f(y)\prec f(u)\prec f(v)$. The following cases are exhaustive. 

\begin{itemize}
\item $f(x)f(y)|f(u)$ and $f(y)f(u)|f(v)$. By canonicity, $f$ behaves as $\nil$ on $S$.  
\item $f(x)f(y)|f(u)$ and $f(y)|f(u)f(v)$. 
By canonicity, $f$ behaves as $\id$ on $S$. 
\item $f(x)|f(y)f(u)$ and $f(y)|f(u)f(v)$. 
By canonicity, $f$ behaves as $\lin$ on $S$. 
\item $f(x)|f(y)f(u)$ and $f(y)f(u)|f(v)$. 
By canonicity, $f(x)|f(y)f(v)$ and $f(x)f(u)|f(v)$. 
It is easy to see that those conditions
are impossible to satisfy over $(\mL;C)$. 
\end{itemize}
\end{proof}

\begin{corollary}\label{cor:exhaust}
Let $\bar c \in \mL^n$ for $n \geq 0$, let $S$ be a cone of $(\mL;C,\bar c)$,
and let $f \colon \mL \to \mL$ be an injection 
that is canonical on $S$ as a function from 
$(\mL;C,\prec,\bar c)$ to $(\mL;C,\prec)$. 
Then $f$ behaves as $\id$, $\lin$, or $\lrl$
on $S$. 
\end{corollary}
\begin{proof}
Note that $f$ is on $S$ canonical as a function from $(\mL;C,\prec)$ to $(\mL;C,\prec)$; also note that every cone contains 
elements $x, y, u, v$ such that $xy|uv$.
Hence, the statement follows from
Lemma~\ref{lem:exhaust}. 
\end{proof}

We finally show 
that if $\Gamma$ is preserved by 
$\lin$, then $\Gamma$ is homomorphically equivalent to a reduct of
$(\mL;=)$; this will be a consequence of
the stronger Lemma~\ref{lem:lin} below. 

\begin{definition}
For $a_1,a_2,\dots,a_k\in \mL$, $k \geq 2$, we write $\Nil(a_1,a_2,\dots,a_k)$ if $a_1\prec a_2\prec\dots\prec a_k$ and $a_1 a_2 \dots a_{i-1}|a_i$ for all
$i\in\{2,\dots,k\}$. 
\end{definition}

Observe that for all $a_1,\dots,a_k \in \mL$ such that $a_1\prec \dots\prec a_k$ we have $\Nil(\nil(a_1,\dots,a_k))$ (recall from Section~\ref{sect:modtheory} that we apply functions to tuples componentwise). Also observe that all $k$-tuples in 
$\Nil$ lie in the same orbit of $k$-tuples. 

\begin{figure}
\begin{center}
\includegraphics[scale=.7]{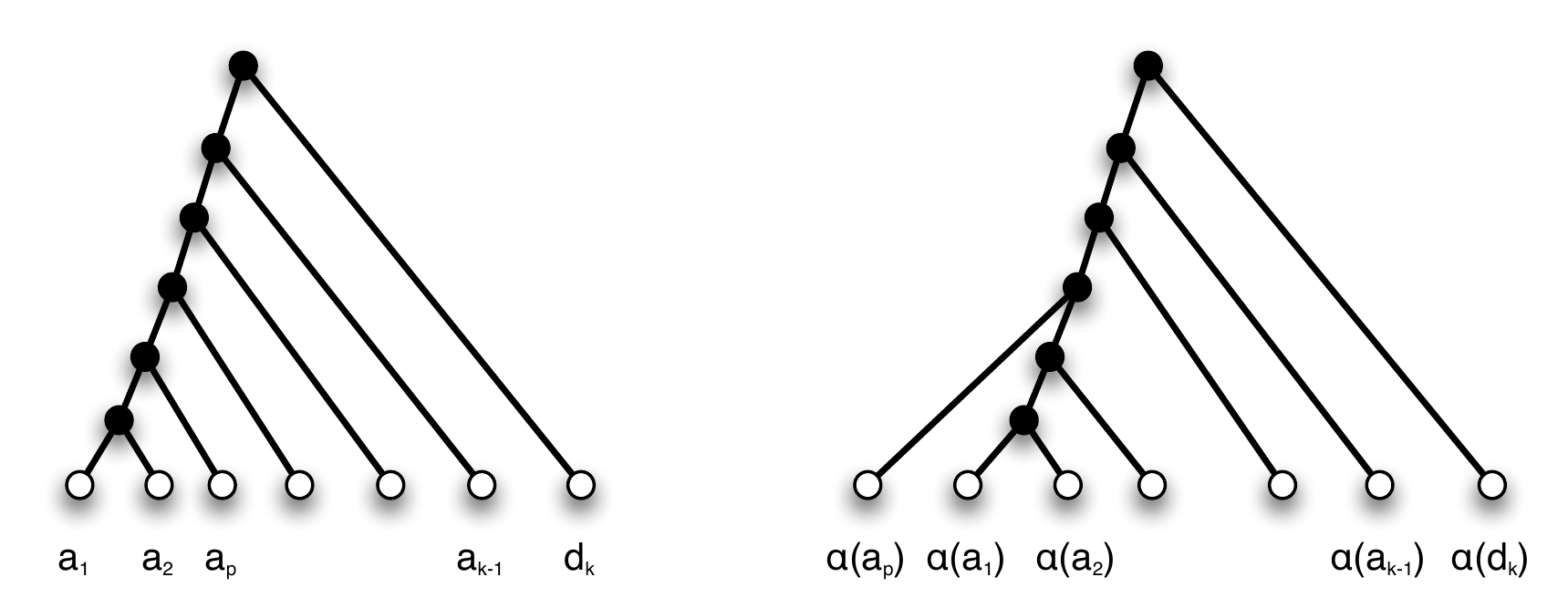} 
\end{center}
\caption{Illustration of the re-ordering of 
Lemma~\ref{lem:gotohead}.}
\label{fig:re-order}
\end{figure}

\begin{lemma}\label{lem:gotohead}
Let $a_1,\dots,a_k\in \mL$ be such that $\Nil(a_1,\dots,a_k)$. Then for every $p\in \{1,\dots,k\}$ there is an $e \in \left<\Aut(\mL;C)\cup \{\nil\}\right>$ such that $\Nil(e(a_p,a_1,\dots,a_{p-1},a_{p+1},a_{p+2},\dots,a_k))$.
\end{lemma}
\begin{proof}
By Proposition~\ref{prop:convex-a} and the homogeneity of $(\mL;C)$ there exists an $\alpha \in \Aut(\mL;C)$ such that $\alpha(a_p)\prec \alpha(a_1)\prec \alpha(a_2)\prec\dots\prec \alpha(a_{p-1})\prec \alpha(a_{p+1})\prec \dots\prec \alpha(a_k)$; see Figure~\ref{fig:re-order}. 
Define $e:=\nil\circ \alpha$. 
By the observation above we have $\Nil(e(a_p,a_1,a_2,\dots,a_{p-1},a_{p+1},a_{p+2},\dots,a_k))$, as desired.
\end{proof}

\begin{figure}
\begin{center}
\includegraphics[scale=.7]{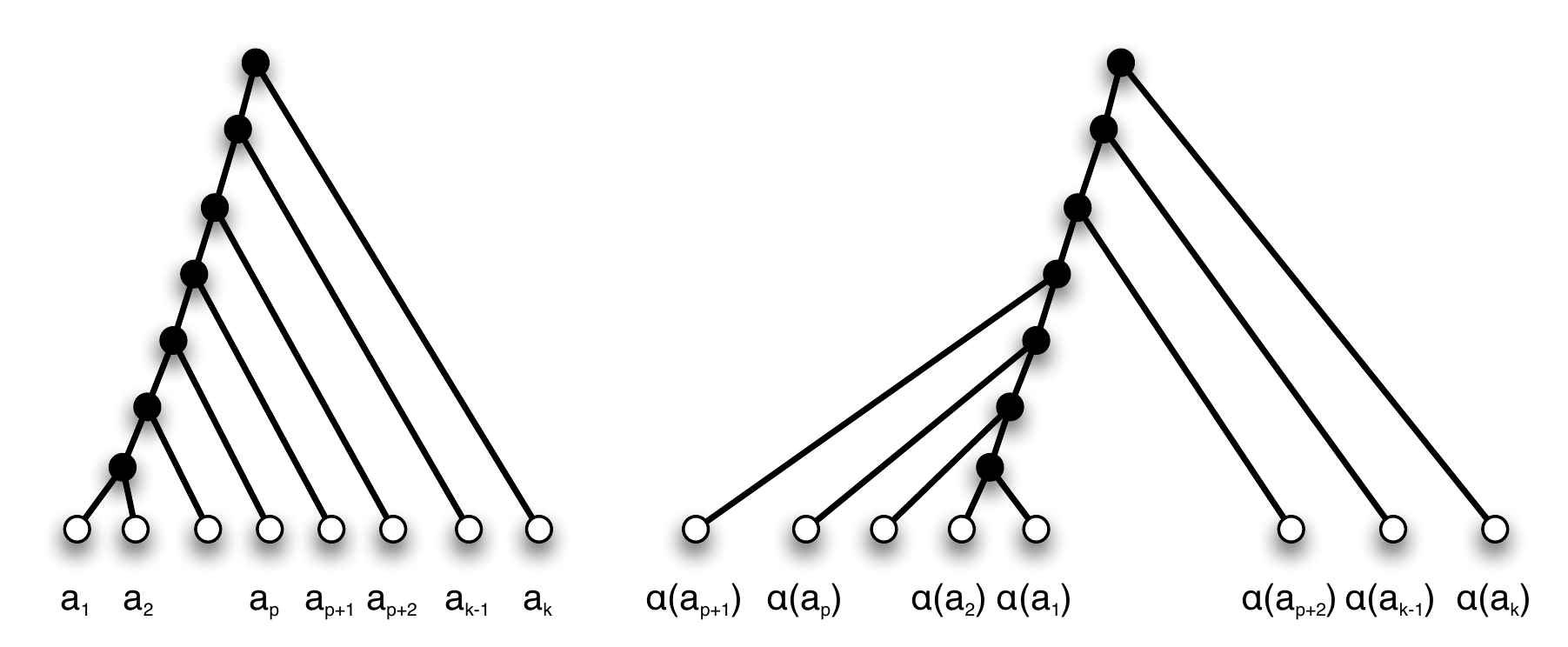} 
\end{center}
\caption{Illustration of the re-ordering of 
Lemma~\ref{lem:localexchange}.}
\label{fig:shuffle}
\end{figure}

\begin{lemma}\label{lem:localexchange}
Let $a_1,\dots,a_k \in \mL$ be such that 
$\Nil(a_1,\dots,a_k)$. Then for any $p\in \{1,2,\dots,k-1\}$, there is an $e \in \left<\Aut(\mL;C)\cup \{\nil\}\right>$ such that 
$e(a_p)=a_{p+1}$, $e(a_{p+1})=a_p$, 
and $e(a_i)=a_i$ for every $i \in \{1,\dots,k\} \setminus \{p,p+1\}$.
\end{lemma}
\begin{proof}
By the homogeneity of $(\mL;C)$ there is an $\alpha \in \Aut(\mL;C)$ such that $\alpha(a_{p+1})\prec \alpha(a_p) \prec \alpha(a_{p-1}) \prec \dots \prec \alpha(a_2) \prec \alpha(a_1)\prec \alpha(a_{p+2})\prec \alpha(a_{p+3})\prec \dots\prec \alpha(a_k)$; see Figure~\ref{fig:shuffle}.  Let $z_i := \nil(\alpha(a_i))$ for $i \in \{1,2,\dots,k\}$. Clearly, we have $\Nil(z_{p+1},z_{p},z_{p-1},\dots,z_2,z_1,z_{p+2},z_{p+3},\dots,z_k)$. 
Starting with the tuple $(z_{p+1},z_{p},z_{p-1},\dots,z_2,z_1,z_{p+2},z_{p+3},\dots,z_k)$ we repeatedly apply Lemma \ref{lem:gotohead} to the resulting tuple at the positions $p := 3,4,\dots,p-1$ in this order. In this way, we obtain in the first step
an $e_1 \in M := \left<\Aut(\mL;C)\cup \{\nil\}\right>$ such that 
$$\Nil(e_1(z_{p-1},z_{p+1},z_{p},z_{p-2},\dots,z_2,z_1,z_{p+2},z_{p+3},\dots,z_k)) \, .$$ 
In the second step, we obtain an $e_2 \in M$ 
such that $$\Nil(e_2(z_{p-2},z_{p-1},z_{p+1},z_{p},z_{p-3},\dots,z_2,z_1,z_{p+2},z_{p+3},\dots,z_k))\, .$$ 
In the $i$-th step, we obtain an $e_i \in M$ such that
$$\Nil(e_i(z_{p-i},z_{p-i+1},...,z_{p-2},z_{p-1},z_{p+1},z_p,
z_{p-i-1},\dots,z_2,z_1,
z_{p+2},z_{p+3},...,z_k)) \, .$$ For $i = p-1$, we therefore 
obtain an $e' \in M$ 
such that $$\Nil(e'(z_1,z_2,\dots,z_{p-2},z_{p-1},z_{p+1},z_p,z_{p+2},\dots,z_k)) \, .$$ 

Define $f:= e' \! \circ \nil \circ \, \alpha$ and observe that 
$\Nil(f(a_1,a_2,\dots,a_{p-1},a_{p+1},a_p,a_{p+2},\dots,a_k))$. 
Therefore, 
$f(a_1,\dots,a_{p-1},a_{p+1},a_p,a_{p+2},\dots,a_k)$
 and $(a_1,a_2,\dots,a_k)$ are in the same orbit in $(\mL;C)$, and there is $\gamma \in \Aut(\mL;C)$ such that 
 $\gamma(f(a_{p},a_{p+1}))=(a_{p+1},a_p)$ and $\gamma(f(a_i))=a_i$ for all $i
 \in \{1,\dots,k\} \setminus \{p,p+1\}$. Then $e:=\gamma\circ f \in M$ has the desired property.
\end{proof}

We write $S_k$ for the symmetric group on $\{1,\dots,k\}$. 

\begin{lemma}\label{lem:permutation}
Let $(x_1,x_2,\dots,x_k),(y_1,\dots,y_k) \in \Nil$, and arbitrarily choose $\delta \in S_k$. 
Then there exists an
$e \in M := \left <\Aut(\mL;C)\cup \{\nil\} \right>$ such that 
$e(x_i) = y_{\delta(i)}$. 
\end{lemma}
\begin{proof}
By Lemma~\ref{lem:localexchange}, for each $p \in \{1,2,\dots,k-1\}$ there is an $e_p \in M$ such that $e_p(x_p)=x_{p+1}$, $e_p(x_{p+1}) = x_p$,
and $e_p(x_i)=x_i$ for all $i \in \{1,\dots,k\} \setminus \{p,p+1\}$. 
Since $S_k$ is generated by the transpositions $(1,2),(2,3),\dots,(k-1,k)$, it follows that there exists $e' \in \left < \{e_p : 1\leq p \leq k\} \right > \subseteq M$ 
such that $e'(x_i)=x_{\delta(i)}$ for all $i \in \{1,\dots,k\}$. 
By the homogeneity of $(\mL;C)$, there
exists an $\alpha \in \Aut(\mL;C)$ such
that $\alpha(x_i) = y_i$ for all $i \in \{1,\dots,k\}$.  Then $e:=\alpha \circ e'$ satisfies $e(x_i) = y_{\delta(i)}$.  
\end{proof}

We can finally prove the announced result. 

\begin{lemma}\label{lem:lin}
Let $\Gamma$ be a reduct of $(\mL;C)$.
Let $\bar c \in \mL^n$ and suppose that
$\Gamma$ has an endomorphism behaving as $\lin$
on a cone of $(\mL;C,\bar c)$. 
Then $\Gamma$ is homomorphically equivalent to
a reduct of $(\mL;=)$.
\end{lemma}
\begin{proof}
Recall that each cone $S$ of $(\mL;C,\bar c)$
induces in $(\mL;C,\prec)$ a structure
that is isomorphic to $(\mL;C,\prec)$.
Let $e$ be an endomorphism of $\Gamma$ that 
behaves as $\lin$ on $S$ and arbitrarily choose a finite
set $A \subseteq \mL$.  By the homogeneity of $(\mL;C,\prec)$ 
there are automorphisms
$\alpha,\beta$ of $(\mL;C,\prec)$ such that $\lin(x) = \alpha(e(\beta(x)))$ for all $x \in A$. 
Since $\End(\Gamma)$ is closed
we have that $\lin \in \End(\Gamma)$ and $\nil \in \End(\Gamma)$ by Proposition~\ref{prop:nil-lin}. 

Our proof has two steps: we first prove
that the structure $\Delta$ 
induced by $D := \nil(\mL)$ is isomorphic
to a reduct of $({\mathbb Q};<)$ and
then we prove in the next step that
$\Delta$ is in fact isomorphic to a reduct
of $(\mL;=)$. 
Clearly, this implies the statement since
$\Gamma$ and $\Delta$ are homomorphically
equivalent. 
Let $\tau$ be the signature of $\Gamma$.
We first show that for every $R \in \tau$,
the relation $R^\Delta$ 
has a first-order definition in $(D;\prec)$. 
The relation $R^\Gamma$ has a first-order definition $\phi$ in $(\mL;C)$. 
An atomic
sub-formula $C(x;yz)$ of $\phi$ holds in $\Delta$ if and only if
$y\prec z\prec x$, $z\prec y\prec x$ or $y=z\wedge x\neq y$. Hence, if we replace in $\phi$ all occurrences of $C(x;yz)$
by $y\prec z\prec x\vee z\prec y\prec x\vee (y=z\wedge x\neq y)$,
we obtain a formula that defines $R^\Delta$ over $(D;\prec)$.
Since $(\mL;\prec)$ is isomorphic to $({\mathbb Q};<)$
and $\nil$ preserves $\prec$, it follows that $(D;\prec)$ is isomorphic to $({\mathbb Q};<)$, too.
Hence, $\Delta$ is isomorphic to a reduct of $({\mathbb Q};<)$.

To show that $\Delta$ is isomorphic to a reduct of $(\mL;=)$, let $X$ be a
finite subset of $D$ 
and $\alpha$ be a permutation of $D$. 
By Proposition~\ref{prop:equality} and the fact that $\End(\Delta)$ is a closed subset of $D^D$, it suffices to find an $e \in \End(\Delta)$ 
such that $e(x) = \alpha(x)$ for all $x \in X$.
Since $X \subseteq \nil(\mL)$, the elements of $X$
can be enumerated by $x_1,\dots,x_n$ such that $\Nil(x_1,\dots,x_n)$.  
Since $\alpha(X) \subseteq D = \nil(\mL)$,
there is a $\gamma \in S_n$ such that 
$\Nil(\alpha(x_{\gamma(1)}),\dots,\alpha(x_{\gamma(n)}))$. 
We apply Lemma~\ref{lem:permutation}
to $(x_1,\dots,x_n)$, $(y_1,\dots,y_n):=(\alpha(x_{\gamma(1)}),\dots,\alpha(x_{\gamma(n)}))$, and $\delta = \gamma^{-1}$, 
and obtain an $f \in \End(\Gamma)$ such 
that $f(x_i)=y_{\delta(i)} = \alpha(x_{\gamma \gamma^{-1} (i)}) = \alpha(x_i)$ for all $i \in \{1,\dots,k\}$. 
The restriction of $\nil \circ f$ to $D$
is an endomorphism $f'$ of $\Delta$.
Since $\nil$ preserves $\prec$ and $\Delta$ is a reduct of $(D;\prec)$ which is isomorphic to $({\mathbb Q};<)$, we have that
$(\alpha(x_1),\dots,\alpha(x_n)) = (f(x_1),\dots,f(x_n))$
and $(f'(x_1),\dots,f'(x_n))$ lie in the same orbit
in $\Delta$. Hence, there exists an 
$\beta \in \Aut(\Delta)$ such that $\beta(f' (x_i)) = \alpha(x_i)$ for all $i \in \{1,\dots,n\}$, and $\beta \circ f'$ is an endomorphism of
$\Delta$ as required. 
\end{proof}

\subsection{Canonical behavior with one constant}
\label{sect:one-constant}
In this section 
we study the behavior of canonical
functions from $(\mL;C,\prec,c_1)$
to $(\mL;C,\prec)$. Some important behaviors are introduced
in Definition~\ref{def:cut-rer}. 
We then show in Sections~\ref{sect:cut}, \ref{sect:rer}, and \ref{sect:tilde-rer} 
that when a function $f$ has some of those behaviors on a $c$-universal set, then $\{f \} \cup \Aut(\mL;C)$ generates
$\lin$ 
or $\End(\mL;Q)$. 
Finally, Section~\ref{sect:one-const-class} classifies behaviors 
of canonical functions from $(\mL;C,\prec,c_1)$
to $(\mL;C,\prec)$. 

\begin{definition}\label{def:cut-rer}
Let $c \in \mL$ and  $A \subseteq \mL \setminus \{c\}$. Let $e \colon \mL \to \mL$ be a function such that 
\begin{itemize}
\item[1.] for any $a\in A$ we have that $e(c)|e(A\cap S^c_a)$, 
\item[2.] for any $a\in A$, $e$ preserves $C$ on $A \cap S^c_a$, and 
\item[3.] for any $a,b\in A$ we have either $S^c_a=S^c_b$ or $e(A\cap S^c_a)|e(A\cap S^c_b)$. 
\end{itemize}
Then we say that $e$ has on $A$ the behavior 
\begin{itemize}
\item[$\id_c$] iff for all $x,y,z \in A$ 
with $x|yzc$ and $y|zc$ we have that 
$e(x)|e(y)e(z)e(c)$ and $e(y)|e(z)e(c)$.
\item[$\cut_c$] iff for all $x,y,z \in A$ 
with $x|yzc$ and $y|zc$ 
we have that 
$e(x)e(y)e(z)|e(c)$ and $e(x)|e(y)e(z)$.
\item[$\rer_c$]
iff for all $x,y,z \in A$ 
with $x|yzc$ and $y|zc$ we have that 
$e(x)e(y)e(z)|e(c)$ and $e(x)e(y)|e(z)$.
\item[$\tilde\rer_c$]
iff for all $x,y \in A$ 
with $x|yc$  we have that $ e(y)|e(x)e(c)$.
\end{itemize}
\end{definition}


\subsubsection{The behavior $\cut_c$}
\label{sect:cut}
Recall that a set $A\subseteq \mL\backslash \{c\}$ is called $c$-universal if for any finite $U\subset \mL$ and $u\in U$, there is $\alpha\in \Aut(\mL;C)$ such that $\alpha(u)=c$ and $\alpha(U)\subseteq A\cup \{c\}$. In this section we prove that for all $c \in \mL$, functions with behavior $\cut_c$ on a $c$-universal set together with $\Aut(\mL;C)$ generate $\lin$. This follows from the following more general fact. 

\begin{lemma}[Cut Lemma]\label{lem:cut}
Let $\Aut(\mL;C) \subseteq M \subseteq \mL^\mL$ be such that
for any finite $U \subset \mL$ and $u \in U$, there exists $g \in M$ that behaves as $\cut_u$ on $U \setminus \{u\}$. 
Then $M$ generates $\nil$ and $\lin$. 
\end{lemma}

\begin{proof}
By Proposition \ref{prop:nil-lin}, it suffices to show that $M$ generates $\nil$. We show that for all $k$ 
and all $x_1,\dots,x_k \in \mL$
 there is an $f \in \left<M\right>$
 such that $f(x_i)=\nil(x_i)$ for all $i \leq k$. 
 We prove this by induction on $k$. 
 Clearly, the claim holds for $k=1$ so assume that $k>1$. 
Suppose without loss of generality that $x_1\prec x_2\prec\dots\prec x_k$. We inductively
assume that there is an $f' \in \left<M\right>$ 
such that $f'(x_i) = \nil(x_i)$
for all $i < k$. 
Let $U := f'(\{x_1,\dots,x_{k}\})$ and
$u := f'(x_k)$. 
By assumption, there exists a $g\in M$ that behaves as $\cut_{u}$ on $U \setminus \{u\}$.

Set $f'' = g \circ f'$. 
We claim that $f''(x_1) \cdots f''(x_i)|f''(x_{i+1})$ for all
$1 \leq i < k$. For $i<k-1$, this follows
from the inductive assumption that 
$f'(x_1) \cdots f'(x_i)|f'(x_{i+1})$, and the assumption
that $g$ behaves as $\cut_u$
on $U \setminus \{u\}$ and in particular preserves $C$ on $U \setminus \{u\}$. 
For $i=k-1$, note that 
that $g(u)|g(U \setminus \{u\})$ 
since $g$ behaves as $\cut_{u}$ on $U \setminus \{u\}$, 
and $f''(x_1),\dots,f''(x_i) \subseteq g(U \setminus \{u\})$. 
Therefore, $f''(x_1) \cdots f''(x_{k-1})|f''(x_k)$ which concludes
the proof of the claim. 

Since $\nil(x_1) \cdots \nil(x_i)|\nil(x_{i+1})$ for all $1 \leq i < k$,
the homogeneity of $(\mL;C)$ implies that there exists an $\alpha \in \Aut(\mL;C) \subseteq M$ such that 
$\alpha(f''(x_i)) = \nil(x_i)$ for all $i \leq k$.
Then $f := \alpha \circ f'' \in \left<M\right>$ has the desired property which concludes the proof.
\end{proof}

\begin{corollary}\label{cor:cut}
Let $c \in \mL$, let $A \subseteq \mL \setminus \{c\}$ be $c$-universal, and let
$g$ be a function that behaves as $\cut_c$ on $A$.
Then $\{g\} \cup \Aut(\mL;C)$ generates $\lin$. 
\end{corollary}
\begin{proof}
The $c$-universality of $A$ implies that 
for every finite $U \subset \mL$ and $u \in U$ there exists an $\alpha \in \Aut(\mL;C)$ such that 
$\alpha(U \setminus \{u\}) \subseteq A$ and $\alpha(u) = c$.
Then $g \circ \alpha$ behaves as $\cut_u$ on $U$
and the statement follows from Lemma~\ref{lem:cut}.
\end{proof}

\subsubsection{The behavior $\rer_c$}
\label{sect:rer}
We will next prove that for all $c \in \mL$,
functions with behavior
$\rer_{c}$ on a $c$-universal set together with $\Aut(\mL;C)$
generate $\End(\mL;Q)$. We need the following lemmas in some later proofs. 

\begin{lemma}\label{lem:q-preservation}
Arbitrarily choose $X \subseteq \mL$ and $c \in X$.
If $f \colon \mL \to \mL$
preserves $Q$ on every 4-element subset of $X$ that contains $c$, then $f$ preserves $Q$ on all of $X$.
\end{lemma} 
\begin{proof}
Let $a_1,a_2,a_3,a_4 \in X$ be such that
$a_1a_2:a_3a_4$. We
show $f(a_1)f(a_2):f(a_3)f(a_4)$. 
It is easy to see that this holds if $a_1,a_2,a_3,a_4$ are not pairwise distinct and it holds by assumption when
$c \in \{a_1,a_2,a_3,a_4\}$. Suppose that this is not the case. Then $a_1a_2c:a_3a_4$ or
$a_1a_2:ca_3a_4$. The latter case can be treated analogously to the former so we assume 
that $a_1a_2c:a_3a_4$. In particular, $a_1c:a_3a_4$
and $a_2c:a_3a_4$ so $f(a_1)f(c):f(a_3)f(a_4)$
and $f(a_2)f(c):f(a_3)f(a_4)$.
It follows from Lemma \ref{D-consequences} that $f(a_1)f(a_2):f(a_3)f(a_4)$ as desired. 
\end{proof}
\begin{lemma}\label{lem:rercqpreservation}
Let $A\subseteq \mL\backslash \{c\}$ and $f:\mL\to \mL$ be such that $f$ behaves as $\rer_c$ on $A$. Then $f$ preserves $Q$ on $A\cup \{c\}$.
\end{lemma}
\begin{proof}
By Lemma \ref{lem:q-preservation} it suffices to show that $f$ preserves $Q$ on $\{x,y,z,c\}$, where $x,y,z$ are pairwise distinct elements in $A$. The following cases are essential.
\begin{itemize}
\item $xyz|c$. \\
\hspace*{0.5cm} It follows from Item 1 and Item 2 of Definition \ref{def:cut-rer} that $f(x)f(y)f(z)|f(c)$ and $f$ preserves $C$ on $\{x,y,z\}$. This implies that $(x,y,z,c)$ and $(f(x),f(y),f(z),f(c))$ are in the same orbit of $\Aut(\mL;C)$. Thus $f$ preserves $Q$ on $\{x,y,z,c\}$.
\item $xy|zc$\\
\hspace*{0.5cm} Clearly, we have $xy:zc$. It follows from Item 1 of Definition of \ref{def:cut-rer} that $f(x)f(y)|f(c)$. It follows from Item 3 of  Definition \ref{def:cut-rer} that $f(x)f(y)|f(z)$. It implies that $f(x)f(y):f(z)f(c)$. Thus $f$ preseves $Q$ on $\{x,y,z,c\}$. 
\item $x|yzc\wedge yz|c$\\
\hspace*{0.5cm} Clearly, we have $xc:yz$. It follows from Item 1 of Definition \ref{def:cut-rer} that $f(y)f(z)|f(c)$. It follows from Item 3 of Definition \ref{def:cut-rer} that $f(x)|f(y)f(z)$. It implies that $f(y)f(z):f(x)f(c)$. Thus $f$ preserves $Q$ on $\{x,y,z,c\}$.
\item $x|yzc\wedge y|zc$\\
\hspace*{0.5cm} Clearly, we have $xy:zc$. By the definition of $\rer_c$ that $f(x)f(y)f(z)|f(c)\wedge f(x)f(y)|f(z)$. It implies that $f(x)f(y):f(z)f(c)$. Thus $f$ preserves $Q$ on $\{x,y,z,c\}$.
\end{itemize}
The other cases can be obtained from the above cases by exchanging the roles of $x,y,z$.
\end{proof}
The following generation lemma is flexible and will be useful later.
\begin{lemma}[Rerooting Lemma, general form]\label{lem:rer2}
Let $\Aut(\mL;C) \subseteq M \subseteq \mL^\mL$ be such that
for any finite $U \subset \mL$ and $u \in U$
there exists $g \in M$ that behaves as 
$\rer_u$ on $U \setminus \{u\}$.
Then $M$ generates $\End(\mL;Q)$. 
\end{lemma}

\begin{proof}
We follow almost literally the proof of Lemma~\ref{lem:rer} (the rerooting lemma).
Arbitrarily choose $f \in \End(\mL;Q)$
and let $A$ be an arbitrary finite subset of $\mL$. 
We have to show that $\left < M \right >$ contains
an operation $e$ such that $e(x) = f(x)$ for all $x \in A$. This is trivial when $|A|=1$ so we assume that $|A| \geq 2$. By Lemma~\ref{lem:invq}, there exists a non-empty proper subset $B$ of $A$ such that 
$f(B) | f(A \setminus B)$ and $B|a$ for all $a \in A \setminus B$. 
By the homogeneity of $(\mL;C)$ we can choose an element $c \in \mL \setminus A$ such that $c|B$ and
$(B \cup \{c\})|a$ for all $a \in A \setminus B$. 
By assumption, there exists a $g \in M$ 
such that $g$ behaves as $\rer_c$ on $A$.
By Lemma \ref{lem:rercqpreservation}, $g$ preserves $Q$ on $A$.

We claim that $g(B)|g(A \setminus B)$. First,
we show that $g(b_1)g(b_2)|g(a)$ 
for every $b_1, b_2 \in B$ 
and $a \in A \setminus B$. 
By the choice of $c$ we have $b_1b_2:ac$. Since $g$ preserves $Q$ on $A$, we have 
 $g(b_1)g(b_2):g(a)g(c)$. Since $g(c)|g(b_1)g(b_2)g(a)$, we have $g(b_1)g(b_2)|g(a)$. 
Next, we show that $g(b)|g(a_1)g(a_2)$ for every $b \in B$ and $a_1, a_2 \in A \setminus B$. 
By the choice of $c$ we have $bc:a_1a_2$. Since $g$ preserves $Q$ on $A$, we have $g(b)g(c):g(a_1)g(a_2)$. Since $g(c)|g(b)g(a_1)g(a_2)$, we have $g(b)|g(a_1)g(a_2)$. 

Let $\beta \colon g(A) \to f(A)$ be defined by
$\beta(x)=f(g^{-1}(x))$ for all $x \in g(A)$. 
Since both $g$ and $f$ preserve
$Q$, we have that $\beta$ preserves $Q$ by Lemma~\ref{lem:pres-q}. Since $\beta(g(B))|\beta(g(A \setminus B))$, the conditions of Lemma~\ref{lem:q-to-c}
apply to $\beta$ for $A_1 := g(B)$ and $A_2 := g(A \setminus B)$, and hence $\beta$ preserves $C$. 
By the homogeneity of $(\mL;C)$, 
there exists an $\gamma \in \Aut(\mL;C) \subseteq M$ that extends $\beta$. Then $e := \gamma \circ g \in \left<M \right>$ has the desired property.
\end{proof}

\ignore{
\begin{lemma}[Rerooting Lemma, general form]\label{lem:rer2}
Let $g \colon \mL \to \mL$ be such that
for every finite $X \subset \mL$ and $x \in X$
there exists an $\alpha \in \Aut(\mL;C)$ such that
$g$ preserves $Q$ on $\alpha(X)$,
and $g(\alpha(X \setminus \{x\}))|g(\alpha(x))$. 
Then $\{g\} \cup \Aut(\mL;C)$ generates $\End(\mL;Q)$. 
\end{lemma}
\begin{proof}
The proof of Lemma~\ref{lem:rer} actually shows the statement above. 
\end{proof}
}

\begin{corollary}\label{cor:rer}
Let $c \in \mL$,
let $A \subseteq \mL \setminus \{c\}$ be $c$-universal, and let $g$ be a function that behaves as 
$\rer_c$ on $A$. Then $\{g\} \cup \Aut(\mL;C)$ generates $\End(\mL;Q)$.
\end{corollary}
\begin{proof}
We claim that the conditions
of Lemma~\ref{lem:rer2} apply to $M := \left<\{g\}\cup \Aut(\mL;C)\right>$. Let $U \subset \mL$ be finite and arbitrarily choose $u \in U$. By $c$-universality of $A$ there exists an $\alpha \in \Aut(\mL;C)$ such that $\alpha(U) \subseteq A$ and $\alpha(u) = c$. Since $g$ behaves as $\rer_c$ on $\alpha(U)$, the function $g\circ \alpha$ behaves as $\rer_u$ on $U$. By Lemma~\ref{lem:rer2} $\left<\{g\} \cup \Aut(\mL;C)\right>$ generates $\End(\mL;Q)$, and so does $\set{g}\cup \Aut(\mL;C)$. 
\end{proof}

By Lemma \ref{lem:rercqpreservation}, a function with behavior $\rer_c$ on $A$ preserves $Q$ on $A$. We now characterize the situation where a function preserving $Q$ behaves as $\rer_c$. 

\begin{lemma}\label{Q:prop}
Let $X$ be a subset of $\mL$, arbitrarily choose $a\in X$, and let 
$f \colon \mL \to \mL$ be a function
 that preserves $Q$ on $X$ and has the property that $f(a)|f(X\setminus \{a\})$. Then $f$ preserves $C$ on $X\cap S^a_x$ for every $x\in X\setminus \{a\}$, and for any $x,y\in X\setminus \{a\}$ either $S^a_x=S^a_y$ or $f(X\cap S^a_x)|f(X\cap S^a_y)$.
\end{lemma}
\begin{proof}
Arbitrarily choose 
$x \in X\setminus \{a\}$ and pick pairwise distinct 
elements $u,v,w\in X \cap S^a_x$. By C8 we can assume that $uv|w$. Clearly, $aw:uv$ and it follows that $f(a)f(w):f(u)f(v)$. Since $f(a)|f(u) f(v) f(w)$ by the assumptions on $f$, we have $f(w)|f(u) f(v)$. This concludes the proof of the first assertion of the lemma.

To show the remaining assertion, let $x,y \in X \setminus \{a\}$ and assume that $S^a_x \neq S^a_y$. 
We claim that $f(X\cap S^a_x)|f(y)$. It suffices to show that $f(x')f(x)|f(y)$ for any $x'\in X \cap  S^a_x$. Clearly, we have $ay:xx'$ and, consequently, $f(a) f(y):f(x) f(x')$. Since $f(a)|f(x) f(x')f(y)$, it follows that $f(x) f(x')|f(y)$. This concludes the proof of the claim.
Similarly, it can be shown that $f(x)|f(X\cap S^a_y)$. 
To complete the argument, arbitrarily choose $x_1,x_2 \in X \cap S^a_x$
and $y_1,y_2 \in X \cap S^a_y$. 
The claim implies that $f(x_1)f(x)|f(y)$ and $f(x_2)f(x)|f(y)$.
Similarly, $f(x)|f(y_1)f(y)$ and $f(x)|f(y_2)f(y)$
because $f(x)|f(X\cap S^a_y)$. 
This implies that $f(x_1)f(x_2)|f(y_1)f(y_2)$, and, consequently, $f(X\cap S^a_x)|f(X\cap S^a_y)$.
\end{proof}

\begin{corollary}\label{Q:rer}
Let $X$ be a subset of $\mL$, arbitrarily choose $a\in X$, and let $f \colon \mL \to \mL$ be a function such that $f(a)|f(X\setminus \{a\})$. Then $f$ behaves as $\rer_a$ on $X$ if and only if $f$ preserves $Q$ on $X$.
\end{corollary}
\begin{proof}
By Lemma \ref{lem:rercqpreservation} if $f$ behaves as $\rer_a$ on $X$ then $f$ preserves $Q$ on $X$. Conversely, suppose that $f$ preserves $Q$ on $X$. By Lemma \ref{Q:prop}, it remains to show that for $x,y,z\in X\setminus \{a\}$ such that $x|yza\land y|za$ we have $f(a)|f(x) f(y)f(z)\land f(z)|f(x)f(y)$. Clearly, we have $xy:za$ and $f(x)f(y):f(z)f(a)$. Since $f(a)|f(x)f(y)f(z)$, it follows that $f(x)f(y)|f(z)$. 
\end{proof}

\subsubsection{The behavior $\tilde\rer_c$}
\label{sect:tilde-rer}
We finally study
functions with behavior
$\tilde\rer_{c}$ on a $c$-universal set. We need the following lemma for the proof of the next important lemma.
\begin{lemma}\label{lem:tildererqpreservation} 
Let $c\in \mL$ and $A\subseteq \mL\backslash \{c\}$ be a $c$-universal set. Let $f:\mL\to \mL$ be a function that behaves as $\tilde \rer_c$ on $A$. Then $f$ preserves $Q$ on $A$. 
\end{lemma}
\begin{proof}
Let $x,y,z,t$ be elements in $A$ such that $xy:zt$. It suffices to show that $f(x)f(y):f(z)f(t)$. Without loss of generality we assume that $xy|z\wedge xy|t$. The following cases are exhaustive.
\begin{itemize}
\item $c|xy$.\\
\hspace*{0.5cm} We will show that $f(x)f(y)|f(z)$. If $xyz|c$ holds, then by Item 2 of Definition \ref{def:cut-rer} we have $f(x)f(y)|f(z)$. If $xy|zc$ holds, then by the definition of $\tilde \rer_c$, we have $f(z)|f(x)f(c)$ and $f(z)|f(y)f(c)$. It follows that $f(x)f(y)|f(z)$. If $z|xyc$ holds, by Item 3 of Definition \ref{def:cut-rer}, we have $f(x)f(y)|f(z)$. These cases are exhaustive, thus $f(x)f(y)|f(z)$. By the same arguments we have $f(x)f(y)|f(t)$. Thus $f(x)f(y):f(z)f(t)$.
\item $\neg c|xy$.\\
\hspace*{0.5cm} We consider the case $cx|y$. The case $cy|x$ are proved similarly. Since $xy|z$ and $xy|t$, it follows that $cxy|z$ and $cxy|t$. By the definition of $\tilde \rer_c$, we have $f(c)f(z)|f(x)\wedge f(c)f(z)|f(y)\wedge f(c)f(t)|f(x)\wedge f(c)f(t)|f(y)$. It follows from $f(c)f(z)|f(x)\wedge f(c)f(t)|f(x)$ that $f(z)f(t)|f(x)$. It follows from $f(c)f(z)|f(y)\wedge f(c)f(t)|f(y)$ that $f(z)f(t)|f(y)$. It implies that $f(z)f(t)|f(x)\wedge f(z)f(t)|f(y)$. Thus $f(x)f(y):f(z)f(t)$.
\end{itemize}
\end{proof}
\begin{lemma}\label{lem:tilde-rer}
Let $c \in \mL$,
let $A \subseteq \mL \setminus \{c\}$ be $c$-universal,
and let $g \colon \mL \to \mL$ be a function 
that behaves as $\tilde\rer_c$ on $A$. 
Then $\{g\} \cup \Aut(\mL;C)$ generates 
$\lin$.
\end{lemma}

\begin{proof}
The proof has two steps. We first show
that $\{g\} \cup \Aut(\mL;C)$ generates $\End(\mL;Q)$, and then prove that $\{g\} \cup \End(\mL;Q)$ generates $\lin$.

For the first step it suffices to show that $\left<\{g\}\cup \Aut(\mL;C)\right>$ satisfies the conditions
of Lemma~\ref{lem:rer2}. Let $U$ be a non-empty finite subset of $\mL$ and arbitrarily choose an element $u \in U$. 
Let $c' \in \mL$ be such that $c' \neq u$ and for every $v \in U \setminus \{u\}$ we have $uc'|v$. Since $A$ is $c$-universal, there is an $\alpha \in \Aut(\mL;C)$ such that $\alpha(U) \subseteq A$ and $\alpha(c')=c$.  Arbitrarily choose two members $v_1,v_2$ of $U \setminus \{u\}$.
 By the choice of $c'$ we have that 
 either $uc'|v_1\land uc'v_1|v_2$, $uc'|v_2\land uc'v_2|v_1$, or $v_1 v_2|u c'$. 
  Since $\alpha$ preserves $C$ and $g$ behaves as $\tilde \rer_c$ on $A$, it follows that $g(\alpha(v_1))g(\alpha(v_2))|g(\alpha(u))$ in each of the three cases. This implies that $g(\alpha(U \backslash\set{u}))|g(\alpha(u))$. 
By Lemma \ref{lem:tildererqpreservation}, $g$ preserves $Q$ on $A$ so $g\circ \alpha$ preserves $Q$ on $U$. By Corollary \ref{Q:rer}  the function $g\circ \alpha$ behaves as $\rer_u$ on $U$. This concludes the first step.

To show that $\{g\} \cup \End(\mL;Q)$ generates $\lin$, we use Lemma~\ref{lem:cut}. Let $U\subseteq \mL$ be finite and arbitrarily choose $u \in U$. Let $v \in \mL$ be such that $U|v$. By the $c$-universality of $A$ there is an $\alpha \in \Aut(\mL;C)$ such that $\alpha(U\cup \{v\})\subseteq A\cup \{c\}$ and $\alpha(u)=c$. Since for every $x \in U\setminus \{ u \}$ we have $c\alpha(x)|\alpha(v)$, it follows that $g(c) g(\alpha(v))|g(\alpha(x))$ for every $x\in U\backslash \{u\}$. This property together with the homogeneity of $(\mL;Q)$ allows us to choose $\beta\in \Aut(\mL;Q)$ such that $\beta(g(c))\beta(g(\alpha(v)))|\beta(g(\alpha(U\setminus \{u\})))$. Let $h=\beta\circ g\circ \alpha$. Since $g$ behaves as $\tilde \rer_c$ on $A$, it follows from Lemma \ref{lem:tildererqpreservation} that $g$ preserves $Q$ on $A$, therefore $g$ preserves $Q$ on $\alpha(\{v\}\cup U\setminus \{u\})$. This implies that $h$ preserves $Q$ on $\{v\}\cup U\setminus \{u\}$. Since $v|U\setminus \{u\}$, $h(v)|h(U\setminus \{u\})$, and $h$ preserves $Q$ on $\{v\}\cup U\setminus \{u\}$, it follows from Lemma \ref{lem:q-to-c} that $h$ preserves $C$ on $U\setminus \{u\}$. Since $h(u)=\beta(g(c))$, we have that $h(u)|h(U\setminus \{u\})$. It follows that $h$ behaves as $\cut_{u}$ on $U$.
\end{proof}

\subsubsection{Classification of behaviors with one constant}
\label{sect:one-const-class}
The main result of this section is
Lemma~\ref{lem:canonical} below, which can be seen as a classification of canonical functions from $(\mL;C,\prec,c)$ to $(\mL;C,\prec)$. 
We first need two lemmata
about $c$-universal sets. 

\begin{lemma}\label{lem:trung-reordering}
Choose $c \in \mL$ and let $A \subseteq \mL \setminus \{c\}$ be a $c$-universal set. Then for every finite subset $X$ of $\mL$ there is an $\alpha \in \Aut(\mL;C,\prec)$ such that $\alpha(X) \subseteq A$ and $\alpha(X) | c$. 
\end{lemma}
\begin{proof}
Recall that the class of all ordered leaf structures has the ordering property (Theorem~\ref{thm:ordering}). By the formulation 
of the ordering property from Proposition~\ref{prop:ordering}, 
there exists a finite subset $Y$ of $\mL$ such
that for every $\beta \in \Aut(\mL;C)$ there exists a $\alpha \in \Aut(\mL;C,\prec)$ such that $\alpha(X) \subseteq \beta(Y)$. Let $y \in \mL$ be such that $y|Y$ holds. Since $A$ is $c$-universal, 
there exists a $\gamma \in \Aut(\mL;C)$ such that $\gamma(y) = c$ and $\gamma(Y) \subseteq A$. By the choice of $Y$ there exists an $\alpha \in \Aut(\mL;C,\prec)$ such that $\alpha(X) \subseteq \gamma(Y) \subseteq A$. Since $\gamma(Y)|c$, we have $\alpha(X)|c$ which concludes the proof. 
\end{proof}

\begin{lemma}\label{lem:equiclas:Cpreserving}
Let $c\in \mL$, and let $A$ be a subset of $\mL$ such that $A\subseteq \{x\in \mL:x\prec c\}$ or $A\subseteq \{x\in \mL:c\prec x\}$, and $A$ is $c$-universal. Let $e \colon \mL\to \mL$ be an injective function that is canonical as a function from $(\mL;C,\prec,c)$ to $(\mL;C)$ and preserves $C$ on $A\cap S^c_a$ for every $a\in A$. Then 
\begin{itemize}
  \item for every $a\in A$ we have $e(c)|e(A\cap S^c_a)$, and
  \item for every $a,b\in A$ we have either $S^c_a=S^c_b$ or $e(A\cap S^c_a)|e(A\cap S^c_b)$.
\end{itemize}  
\end{lemma}
\begin{proof}
To prove the first assertion of the lemma, it suffices to prove that for arbitrary $x,y \in A$ satisfying $xy|c$ we have $e(x)e(y)|e(c)$. Assume for contradiction that there are $x_1,x_2 \in A$ such that $x_1x_2|c$ and $e(x_1)|e(x_2)e(c)$. Let $x_3 \in \mL$ be such that $x_1\neq x_3\wedge x_1x_3|x_2$ holds. By Lemma~\ref{lem:trung-reordering} there is an $\alpha \in \Aut(\mL;C,\prec)$ such that $\alpha(\{x_1,x_2,x_3\}) \subseteq A$ and $\alpha(\{x_1,x_2,x_3\})|c$. 
For $i \in \{1,2,3\}$ let $x_i' := \alpha(x_i)$. 
The pairs $(x_1,x_2)$, $(x_1',x_2')$,
and $(x_3',x_2')$ are in the same orbit of $\Aut(\mL;C,\prec,c)$: we have $x_1 \prec x_2$ if and only if $x_1' \prec x_2'$ since $\alpha$ preserves $\prec$. Further, $x_1' \prec x_2'$ if and only if $x_3' \prec x_2'$ by convexity of $\prec$ since $x_1x_3|x_2$ holds and $\alpha$ preserves $C$.
Moreover, it holds that $x_1x_2|c$ by assumption, and $x_1'x_2'|c$ and $x_3'x_2'|c$ by the properties of $\alpha$. In case that $A \subseteq \{x \in \mL : x \prec c\}$ we have
$x_1,x_2,x_3,x_1',x_2',x_3' \prec c$, otherwise we have $c \prec x_1,x_2,x_3,x_1',x_2',x_3'$.
By the homogeneity
of $(\mL;C,\prec,c)$ we conclude that $(x_1,x_2)$, $(x_1',x_2')$,
and $(x_3',x_2')$ are indeed in the same orbit
of $\Aut(\mL;C,\prec,c)$.

By the canonicity of $e$, we have that $e(x_1')|e(x_2')e(c)$
and $e(x_3')|e(x_2')e(c)$. Since $e$ preserves
$C$ on $S_x \cap A$, we have $e(x_1')e(x_3')|e(x_2')$ which implies that 
$e(x_1')e(x_3')|e(x_2')e(c)$. Since $x_1'x_2'x_3'|c$ and $x_1',x_2',x_3'$ are pairwise distinct, it follows that either $(x_1',x_3')$ or $(x_3',x_1')$ is in the same orbit as $(x_1',x_2')$ in $\Aut(\mL;C,\prec,c)$. If $(x_1',x_3')$ and $(x_1',x_2')$ are in the same orbit in $\Aut(\mL;C,\prec,c)$, then by the canonicity of $e$, $(e(x_1'),e(x_3'),e(c))$ and $(e(x_1'),e(x_2'),e(c))$ are in the same orbit in $\Aut(\mL;C)$. It is impossible since $e(x_1')e(x_3')|e(x_2')e(c)$. The case $(x_3',x_1')$ and $(x_1',x_2')$ are in the same orbit of $\Aut(\mL;C,\prec,c)$ is proved similarly. The first assertion of the lemma therefore follows.   

It remains to show the second assertion of the lemma. We consider the case $A\subseteq \{x\in \mL:x\prec c\}$. The case $A\subseteq \{x\in \mL:c\prec x\}$ is argued similarly. We need the two following claims.

{\bf Claim 1.} \emph{For any $x_1,x_2,x_3 \in A$
satisfying $x_1|x_2x_3c$ and $x_2x_3|c$ we have $e(x_1)|e(x_2)e(x_3)$.}

{\textit{Proof of Claim 1.}}  For a contradiction we assume that $e(x_1)e(x_2)|e(x_3)$. 
Let $y_1,\dots,y_5 \in A$ be pairwise distinct such that $y_1|y_2y_3y_4y_5c$,
$y_2y_3y_4y_5|c$, and $y_2y_3|y_4y_5$. It follows from the convexity of $\prec$ that
$y_1 \prec y_i$ for $i \in \{2,\dots,5\}$. 
Since $y_2y_3|y_4y_5$, we have either $y_i \prec y_j$ for $i \in \{2,3\}$ and $j \in \{4,5\}$, 
or $y_i \prec y_j$ for $i \in \{4,5\}$ and $j \in \{2,3\}$. Without loss of generality we assume that $y_2 \prec y_3 \prec y_4 \prec y_5$. If $x_2 \prec x_3$, then the tuples $(x_1,x_2,x_3)$, $(y_1,y_2,y_4)$, and $(y_1,y_4,y_5)$ are in the same orbit of $\Aut(\mL;C,\prec,c)$. Thus $e(y_1)e(y_2)|e(y_4)$ and
$e(y_1)e(y_4)|e(y_5)$. Since $e$ preserves $C$ on $A\cap S^c_{y_2}$ and since $\{y_2,y_3,y_4,y_5\}\subseteq A\cap S^c_{y_2}$, we have $e(y_2)|e(y_4)e(y_5)$. These conditions are impossible to satisfy over $(\mL;C)$. 
If $x_3 \prec x_2$ then we consider the three tuples $(x_1,x_3,x_2)$, $(y_1,y_2,y_3)$, and $(y_1,y_3,y_4)$ that lie in the same orbit of 
$\Aut(\mL;C,\prec,c)$ and proceed analogously.\qed

{\bf Claim 2.} \textit{For any $x_1,x_2,x_3 \in A$ satisfying $x_1x_2|x_3c$ we have $e(x_1)e(x_2)|e(x_3)$.}

This claim can be shown similarly as for the claim above by choosing five distinct elements $y_1,\dots,y_5$ such that $y_1y_2y_3y_4|y_5c$ and $y_1y_2|y_3y_4$.

Let $a,b\in A$ such that $S_a^c\neq S_b^c$. This condition implies that $ca|b$ or $cb|a$. We consider the case $ca|b$. The case $cb|a$ is argued similarly. By the definition of cones, we have $S_a^c|S_b^c$, thus $A\cap S_a^c|A\cap S_b^c$. Let $x,y\in A\cap S_a^c$ and $z,t\in A\cap S_b^c$ be arbitrary. It follows from $ca|b\wedge xya|c\wedge ztb|c$ that $z|xyc\wedge xy|c\wedge t|xyc\wedge xy|c\wedge zt|xc\wedge zt|yc$. It follows from Claim 1 that $e(z)|e(x)e(y)\wedge e(t)|e(x)e(y)$, and it follows from Claim 2 that $e(z)e(t)|e(x)\wedge e(z)e(t)|e(y)$. Thus $e(x)e(y)|e(z)e(t)$. The second assertion follows.
\end{proof}

\begin{lemma}\label{lem:canonical}
Let $c\in \mL$, and let $A$ be a subset of $\mL$ such that $A\subseteq \{x\in \mL:x\prec c\}$ or $A\subseteq \{x\in \mL:c\prec x\}$, and $A$ is $c$-universal. Let $e \colon \mL \to \mL$ be an injective 
function that 
is canonical on $A$ as a function 
from $(\mL;C,\prec,c)$ to $(\mL;C,\prec)$. Then $\{e\} \cup \Aut(\mL;C)$ generates $\lin$
or $e$ behaves on $A$ as $\id_{c}$ or $\rer_{c}$. 
\end{lemma}
\begin{proof}
We consider the case $A\subseteq \{x\in \mL:x\prec c\}$. The case $A \subseteq \{x \in \mL : c \prec x\}$ is argued similarly. The canonicity of $e$ implies that either
\begin{itemize}
\item $e(x) \prec e(y)$ for all $x,y \in A$ such that $xy|c$ and $x \prec y$, or 
\item $e(y) \prec e(x)$ for all $x,y \in A$ such that $xy|c$ and $x \prec y$. 
\end{itemize}
If the second case applies, we continue the proof
with $- \circ e$ instead of $e$. Thus we assume in the following that the first case applies.

Since $A$ is $c$-universal, there is an $a \in A$
such that $S^c_a \cap A$ contains four distinct elements $x,y,u,v$ satisfying $xy|uv$. 
The function $e$ is canonical on $S^c_a \cap A$
as a function from $(\mL;C,\prec)$ to $(\mL;C,\prec)$. Lemma~\ref{lem:exhaust} shows that
$e$ behaves as $\id$, $\lin$, or $\lrl$ on $S^c_a \cap A$. The canonicity of $e$ implies that $e$ has the same behavior on all sets of the form $S^c_x \cap A$ for $x \in A$. 

If $e$ behaves as $\lin$ on all those sets, then we show that
$\Aut(\mL;C) \cup \{e\}$ generates $\lin$. Let $X$ be a finite subset of $\mL$. 
By Lemma~\ref{lem:trung-reordering} there is an $\alpha \in \Aut(\mL;C,\prec)$ such that 
$\alpha(X) \subseteq A$ and $\alpha(X) | c$. 
Since $e$ behaves as $\lin$ on $\alpha(X)$
and $\alpha$ preserves $\prec$, the function
$e \circ \alpha$ behaves as $\lin$ on $X$. This implies that $\{e\} \cup \Aut(\mL;C)$ generates $\lin$. 

If $e$ behaves as $\lrl$ on all sets of the form $S^c_x \cap A$ for $x \in A$, then by the same argument it can be shown that $\{e \} \cup \Aut(\mL;C)$ generates $\lrl$, and therefore $\lin$ by Proposition~\ref{prop:nil-lin}. We therefore assume in the following that $e$ preserves $C$ on each set of the form $S^c_x \cap A$, $x \in A$.  Lemma \ref{lem:equiclas:Cpreserving} implies that the preconditions of the behaviors are satisfied (Definition~\ref{def:cut-rer}). Let $u,v \in A$ be such that $u|vc$. Clearly, we have $u \prec v$. By C8, the following cases are exhaustive. 
\begin{enumerate}
\item $e(u)|e(v)e(c)$. Let $x,y,z \in A$ be such that $x|yzc$ and $y|zc$. Clearly, we have $x \prec y \prec z \prec c$. Since $(x,y)$ and $(y,z)$ are in the same orbit of $\Aut(\mL;C,\prec,c)$ as $(u,v)$, we have $e(x)|e(y)e(c)$ and $e(y)|e(z)e(c)$ implying that $e(x)|e(y)e(z)e(c)$. Thus, $e$ behaves as $\id_c$ on $A$. 
\item $e(v)|e(u)e(c)$. Let $x,y,z \in A$ be such that
$x|yzc$ and $y|zc$. Since $(x,y)$ and $(y,z)$ are in the same orbit as $(u,v)$ in $\Aut(\mL;C,\prec,c)$, we have $e(y)|e(x)e(c)$ and $e(z)|e(y)e(c)$, and therefore $e(z)|e(x)e(y)e(c)$. Thus, $e$ behaves as $\tilde \rer_c$ on $A$, and $\{e\} \cup \Aut(\mL;C)$ generates $\lin$ by Lemma~\ref{lem:tilde-rer}. 
\item $e(u)e(v)|e(c)$. Canonicity of $e$ on $A$ implies that for any $x,y \in A$ we have $e(x)e(y)|e(c)$ and
thus $e(A)|e(c)$. Let $a_1,a_2,a_3 \in A$ be three distinct elements such that $a_1|a_2a_3c$ and $a_2|a_3c$. By the convexity of $\prec$ we have $a_1 \prec a_2 \prec a_3$ so
$(a_1,a_2)$ and $(a_2,a_3)$ are in the same orbit  of $\Aut(\mL;C,\prec)$ as $(u,v)$. The canonicity of $e$ implies that either $e(a_1) \prec e(a_2) \prec e(a_3)$ or $e(a_3) \prec e(a_2) \prec e(a_1)$ holds.
It follows from the convexity of $\prec$ that either $e(a_1)|e(a_2)e(a_3)$ or $e(a_1)e(a_2)|e(a_3)$ holds. If the first case holds then $e$ behaves as $\cut_c$, and if the second case holds then $e$ behaves as $\rer_c$ on $A$. 
\end{enumerate}
We conclude that unless $\{e\} \cup \Aut(\mL;C)$ generates $\lin$, it behaves as $\id_c$ or $\rer_c$ on $A$. 
\end{proof}

%% file: two-constants.tex
\subsection{Canonical behavior with two constants}
\label{sect:two-constants}
In this section we analyze canonical functions from $(\mL;C;\prec,c_1,c_2)$ to $(\mL;C;\prec)$.
For our purposes, it suffices to 
treat some special behaviors (Lemma~\ref{lem:two-constants-up}); the motivation for those behaviors
will become clear in the proof of Proposition~\ref{prop:rooted}. 
Then, we prove that certain behaviors of $f$ imply that $\{f\} \cup \Aut(\mL;C)$ generates $\lin$ (Lemma~\ref{lem:destruction} and Lemma~\ref{lem:destruction-2}).

\begin{definition}
Let $c_1,c_2 \in \mL$ be distinct. 
Then $A \subset \mL \setminus \{c_1,c_2\}$ is called 
\emph{$(c_1,c_2)$-universal} if 
for every finite $U \subset \mL$ and 
$u_1,u_2 \in U$ there is an $\alpha \in \Aut(\mL;C)$ such that
$\alpha(U) \subseteq A \cup \{u_1,u_2\}$, $\alpha u_1 = c_1$, and $\alpha u_2 = c_2$.
\end{definition}

Note that when $A$ is $(c_1,c_2)$-universal then this implies
in particular that $\{x \in A : xc_1|c_2\}$ and $\{x \in A : x|c_1c_2\}$ are $c_1$-universal.

\begin{lemma}\label{lem:two-constants-up}
Let $c_1,c_2 \in \mL$ be distinct,
and let $A$ be $(c_1,c_2)$-universal such that
all elements in $A_1 := \{x \in A : xc_1|c_2\}$ are in the same orbit in $(\mL;C;\prec,c_1,c_2)$,
and all elements in $A_2 := \{x \in A : x|c_1c_2\}$ are in the same orbit in $(\mL;C;\prec,c_1,c_2)$. 
Let $f \colon \mL \to \mL$ be canonical on $A$
as a function from $(\mL;C;\prec,c_1,c_2)$ 
to $(\mL;C;\prec)$.
Then $\{f\} \cup \Aut(\mL;C)$ generates $\lin$ or $\End(\mL;Q)$, or $f$ preserves $C$ on $\{c_1\} \cup A_1 \cup A_2$. 
\end{lemma}
\begin{proof}
It follows from the assumption on $A_1$ and $A_2$ that $f$ is canonical on $A_1$ and $A_2$ as a function from $(\mL;C,\prec,c_1)$ to $(\mL;C,\prec)$. By Lemma \ref{lem:canonical}, if $f$ does not preserve $C$ on $A_1\cup \set{c_1}$ 
 then $\set{f} \cup \Aut(\mL;C)$ 
 generates $\lin$ and we are done, or it behaves as $\rer_{c_1}$ 
 on $A_1\cup \set{c_1}$
 in which case
 $\set{f} \cup \Aut(\mL;C)$ 
 generates $\End(\mL;Q)$ by Corollary~\ref{cor:rer}, and we are again done. 
The same argument applies if $f$ does not preserve $C$ on $A_2\cup \set{c_2}$.

 Thus, it remains to consider the case when $f$ preserves $C$ on both $A_1\cup \set{c_1}$ and on $A_2\cup \set{c_1}$. Let $a_1\in A_1$ and $a_2 \in A_2$. We distinguish the following cases:
\begin{itemize}
  \item $f(c_1)|f(a_1) f(a_2)$. 
It follows from the canonicity of $f$ that $f(c_1)|f(a_1) f(x)$ for all $x\in A_2$ so $f(c_1)|f(A_2)$. This is impossible since $f$ preserves $C$ on $\set{c_1}\cup A_2$.
 \item $f(a_2) f(c_1)|f(a_1)$. 
It follows from the canonicity of $f$ that $f(x)f(c_1)|f(a_1)$ for all $x\in A_2$, therefore $f(a_1)|f(A_2)$.
 This implies that $f$ behaves as $\cut_{a_1}$ on $A_2$. 
 Since $c_1a_1|x$ for all $x \in A_2$ and $A_2$ is $c_1$-universal, we have
 that $A_2$ is also $a_1$-universal, and,  
 by Corollary~\ref{cor:cut}, $\set{f}\cup \Aut(\mL;C)$ generates $\lin$.
\item $f(c_1)f(a_1)|f(a_2)$. 
It follows from the canonicity of $f$ that $f$ preserves $C$ on $\set{c_1}\cup A_1\cup A_2$, and we are done.
\end{itemize}  
Since these three cases are exhaustive, the statement follows. 
\end{proof}

\begin{lemma}\label{lem:destruction}
Let $c_1,c_2 \in \mL$ 
and $A \subseteq \mL \setminus \{c_1,c_2\}$ be
$(c_1,c_2)$-universal. 
Let $A_1 = \{x \in A \, : \, xc_1|c_2\}$, 
$A_2 = \{x \in A \, : \, c_1|xc_2\}$,
and $A_3 = \{x \in A \, : \, c_1c_2|x\}$. 
Let $g \colon \mL \to \mL$ be an injection such that 
\begin{itemize}
\item  $g(A_1 \cup \{c_1\})|g(c_2)$, 
\item $g$ preserves $C$ on $\{c_1\} \cup A_1$
and on $\{c_2\} \cup A_2 \cup A_3$, and
\item $g(c_1)g(c_2)|g(x)$ for every $x \in A_2 \cup A_3$. 
\end{itemize}
Then $\{g\} \cup \Aut(\mL;C)$ generates $\nil$.
\end{lemma}
\begin{proof}
We need to show that for all $k$ 
and all $x_1,\dots,x_k \in \mL$
 there is an $f \in M := \left< \Aut(\mL;C) \cup \{g\} \right>$
 such that $f(x_j)=\nil(x_j)$ for all $j \leq k$. 
This is clearly true for $k \leq 2$.
  To prove it for $k \geq 3$, 
suppose without loss of generality that $x_1 \prec \cdots \prec x_k$.   
We first 
 prove by induction on $i \in \{1,\dots,k-1\}$ that there exists
an $h \in M$ with the following properties.
\begin{enumerate}
\item $h(x_1) \cdots h(x_i) | h(x_j)$ for every $j \in \{i+1,\dots,k\}$
\item $h$ preserves $C$ on $\{x_i,\dots,x_k\}$
\item for $i \geq 2$ we additionally require that 
$h(x_1) \cdots h(x_j) | h(x_{j+1})$ 
for every $j \in \{1,\dots,i-1\}$ 
\end{enumerate}
For $i=1$ the identity function has the properties that we
require for $h \in M$. 
For $i \geq 2$, we inductively assume the
 existence of a function $h' \in M$ such that 
\begin{itemize}
\item $h'(x_1)\cdots h'(x_{i-1}) | h'(x_j)$ for every $j \in \{i,\dots,k\}$,
\item $h'$ preserves $C$ on $\{x_{i-1},\dots,x_k\}$, and
\item if $i \geq 3$ we additionally have 
$h(x_1) \cdots h(x_j) | h(x_{j+1})$
for every $j \in \{1,\dots,i-2\}$.  
\end{itemize}
By $(c_1,c_2)$-universality of $A$, there exists an $\alpha \in \Aut(\mL;C)$ that maps $h'(x_1)$ to $c_1$, $h'(x_i)$ to $c_2$, 
and such that $\alpha h'(\{x_1,\dots,x_k\}) \subseteq A \cup \{c_1,c_2\}$. 

\vspace{.2cm}
{\bf Observation.}
 $\alpha h'(\{x_1,\dots,x_{i-1}\}) \subseteq \{c_1\} \cup A_1$ and $\alpha h'(\{x_i,\dots,x_k\}) \subseteq  \{c_2\} \cup A_2 \cup A_3$. 
\begin{proof}[Proof of the observation.]
The first property of $h'$ 
implies that $h'(x_1) h'(x_{i-1})|h'(x_i)$.
Therefore, $\alpha h'(x_1) \alpha h'(x_{i-1}) | \alpha h'(x_i)$ and $c_1 x|c_2$ for every
$x \in \alpha h'(\{x_1,\dots,x_{i-1}\})$ which concludes the proof of the first part of the observation.

To show the second part, arbitrarily choose $j \in \{i,\dots,k\}$. If $j=i$ then $\alpha h'(x_j) = \alpha h'(x_i) = c_2$ and there is nothing to show. 
Since $x_{i-1} \prec x_i \prec x_j$, we distinguish the
cases that $x_{i-1}|x_ix_j$ and $x_{i-1}x_i|x_j$. By the inductive assumption, $h'$ preserves $C$ 
on $\{x_{i-1},\dots,x_k\}$ so we have 
$$h'(x_{i-1})|h'(x_i)h'(x_j)\text{ or }h'(x_{i-1})h'(x_i)|h'(x_j).$$
First consider the case $h'(x_{i-1})|h'(x_i)h'(x_j)$.
By the first property of $h'$, we also have
$$h'(x_{1})h'(x_{i-1})|h'(x_j)\text{ and }h'(x_1) h'(x_{i-1})|h'(x_i)h'(x_j).$$ 
Consequently, $\alpha h'(x_1) | \alpha h'(x_i) \alpha h'(x_j)$, 
and thus $c_1|c_2 \alpha h'(x_j)$. Hence $\alpha h'(x_j) \in A_2$. 
Now consider the case $h'(x_{i-1})h'(x_i)|h'(x_j)$.
Since $h'(x_{1})h'(x_{i-1})|h'(x_j)$, we have
$h'(x_1) h'(x_{i})|h'(x_j)$. Thus 
$c_1 c_2 |\alpha h'(x_j)$, and $\alpha h'(x_j) \in A_3$.
\end{proof}

We claim that $h := g \circ \alpha \circ h'$ satisfies the inductive claim so we have to verify the three
properties from the inductive statement.
\begin{enumerate}
\item[Ad 1.]
By the observation above with the facts that
$\alpha h'(x_i) = c_2$ and
$\alpha h'$ is injective, it follows that
$\alpha h'(\{x_{i+1},\dots,x_k\}) \subseteq A_2 \cup A_3$.
Since $g(c_1)g(c_2)|g(x)$ for every $x \in A_2 \cup A_3$ and $g(A_1 \cup \{c_1\})|g(c_2)$, we have $g(\{c_1,c_2\} \cup A_1)|g(x)$ for every $x \in A_2 \cup A_3$. 
Therefore, $(g \circ \alpha \circ h')(\{x_1,\dots,x_i\}) | (g \circ \alpha \circ h')(x_j)$ for every $j \in \{i+1,\dots,k\}$, or, equivalently, $h(x_1) \cdots h(x_{i}) | h(x_j)$,
which is what we had to show. 
\item[Ad 2.]
By the second property of $h'$, the restriction of $h'$
to $\{x_{i-1},\dots,x_k\}$ preserves $C$. 
Since $\alpha h'(\{x_i,\dots,x_k\}) \subseteq \{c_2\} \cup A_2 \cup A_3$ and $g$ preserves $C$ over $\{c_2\} \cup A_2 \cup A_3$, the restriction of $h = g \circ \alpha \circ h'$ to $\{x_i,\dots,x_k\}$ preserves $C$ as well. 
\item[Ad 3.] We assume that $i \geq 3$ since otherwise
there is nothing to show. 
Since $g$ preserves $C$ over $A_1 \cup \{c_1\}$ and $\alpha h'(\{x_1,\dots,x_{i-1}) \subseteq A_1 \cup \{c_1\}$, the third property of $h'$ implies that 
$g \circ \alpha \circ h'(\{x_1,\dots,x_j\}) | g \circ \alpha \circ h'(x_{j+1})$ for all $j \in \{1,\dots,i-2\}$. 
Equivalently, $h(x_1) \cdots h(x_{j})|h(x_{j+1})$
for all $j \in \{1,\dots,i-2\}$. 
It remains to show
that $h(x_{i-2})h(x_{i-1})|h(x_{i})$. 
This 
follows directly from the fact that $g(A_1 \cup \{c_1\})|g(c_2)$ and $\alpha h'(x_{i}) = c_2$. 
\end{enumerate}
This concludes the induction. 
For $i = k$ the third property of $h$ 
implies that
$h(x_1) \cdots h(x_{j}) | h(x_{j+1})$ for all 
$j \in \{1,\dots,k-1\}$. This property and the homogeneity of $(\mL;C)$ imply the existence of $\beta \in \Aut(\mL;C)$ such that
$\beta h(x) = \nil(x)$ for all $x \in X$, and
hence $f:=\beta \circ h \in M$ is a function with the desired properties.  
\end{proof}

\begin{lemma}\label{lem:destruction-2}
Let $c_1,c_2 \in \mL$ 
and $A \subseteq \mL \setminus \{c_1,c_2\}$ be
$(c_1,c_2)$-universal.
Let $A_1 = \{x \in A \, : \, xc_1|c_2\}$, 
$A_2 = \{x \in A \, : \, c_1|xc_2\}$,
and $A_3 = \{x \in A \, : \, c_1c_2|x\}$. 
Let $g \colon \mL \to \mL$ be an injection such that 
\begin{itemize}
\item for all $a_1 \in A_1$, $a_2 \in A_2$ we either
have $g(c_1)g(a_1)|g(a_2)$ or $g(c_1)g(a_2)|g(a_1)$;
\item $g$ preserves $C$ on $\{c_1\} \cup A_1 \cup A_3$ and on $\{c_2\} \cup A_2 \cup A_3$, and
\item $g(c_1)g(c_2)|g(x)$ for every $x \in A$. 
\end{itemize}
Then $\{g\} \cup \Aut(\mL;C)$ generates $\lin$. 
\end{lemma}
\begin{proof}
We first show that $\{g\} \cup \Aut(\mL;C)$
generates a function $f$ with the property that
there are no $a,b,c,d \in \mL$ such that
$f(a)f(b)|f(c)f(d)$. For this, it suffices by a standard application of K\"onig's tree lemma (see e.g.~Section 3.1 in~\cite{Bodirsky-HDR}) 
to show that 
for all finite $S = \{x_1,\dots,x_k\} \subset \mL$ there
is an $h \in M := \left< \Aut(\mL;C) \cup \{g\} \right>$
such that there are no $a,b,c,d \in S$ with
$h(a)h(b)|h(c)h(d)$.

This is clearly true for $k \leq 1$.
  To prove it for $k \geq 2$, 
 we prove by induction on $i \in \{1,\dots,k-1\}$ 
 that there exists
an $h \in M$ with the following property.
\begin{enumerate}
\item[$P_{h,i}$] The equivalence relation $\sim_h$ defined on $\{x_2,\dots,x_k\}$ by $u \sim_h v$ iff $h(x_1)|h(u)h(v)$ has
at least $i$ equivalence classes. 
\end{enumerate}
For $i=1$ the statement is trivial and we let $h \in M$ be the identity function.  
For $i \geq 2$, we inductively assume 
the existence of a function $h' \in M$ satisfying 
$P_{h',i-1}$. If $P_{h',i}$ holds,  
then there is nothing to be shown so we assume that there are 
distinct $p,q \leq k$ such that $h'(x_1)|h'(x_p)h'(x_q)$. 
By $(c_1,c_2)$-universality of $A$, there exists an $\alpha \in \Aut(\mL;C)$ that maps $h'(x_1)$ to $c_1$, $h'(x_p)$ to $c_2$, 
and such that $\alpha h(\{x_1,\dots,x_k\}) \subseteq A \cup \{c_1,c_2\}$. We claim that $h := g \circ \alpha \circ h'$ satisfies $P_{h,i}$. To show that, we first prove that 
$\sim_h$ has at least as
many equivalence classes as $\sim_{h'}$. 

Observe that when $r,s \in \{2,\dots,k\}$ are such that 
$x_r \nsim_{h'} x_s$, then
$x_r \nsim_{h} x_s$:
when both $\alpha(h'(x_r)),\alpha(h'(x_s)) \in A_1 \cup A_3$ then this follows the assumption that $g$ preserves $C$ on $\{c_1\} \cup A_1 \cup A_3$; a similar argument applies when both $\alpha(h'(x_r)),\alpha(h'(x_s)) \in A_2 \cup A_3$. When $\alpha(h'(x_r)) \in A_1$
and $\alpha(h'(x_s)) \in A_2$ then either 
$g(c_1)g(\alpha(h'(x_r)))|g(\alpha(h'(x_s)))$ or 
$g(c_1)g(\alpha(h'(x_s)))|g(\alpha(h'(x_r)))$, so $x_r \nsim_{h} x_s$. 

Next, consider the case that one of $r,s$, say $r$, equals $p$,
that is, $\alpha(h'(x_r)) = c_2$. In this case we have
by the third assumption on $g$ in the statement of the lemma that
$g(\alpha(h'(x_1)))g(\alpha(h'(x_r)))|g(x)$ for all $x \in A$,
and in particular that $g(\alpha(h'(x_1)))g(\alpha(h'(x_r)))|g(\alpha(h'(x_s)))$. Hence, $x_r \nsim_{h} x_s$. 

To see that $\sim_h$ has \emph{strictly} more equivalence
classes than $\sim_{h'}$, observe that $x_p \sim_{h'} x_q$ but
$x_p \nsim_{h} x_q$ as we have just seen.
This concludes the inductive proof. 

Note that $P_{h,k-1}$ 
implies that for all $p< q \leq k$ we have
$h(x_p)h(x_1)|h(x_q)$ or $h(x_q)h(x_1)|h(x_p)$,
and in particular there cannot be $a,b,c,d \in S$
such that $h(a)h(b)|h(c)h(d)$. This concludes our
proof of the existence of $f$.

Since $(\mL;C,\prec)$ is homogeneous, Ramsey, 
and $\omega$-categorical, Theorem~\ref{thm:can}
asserts the existence of a function $f' \in \overline{ \{\alpha_1 f \alpha_2 \; | \; \alpha_1,\alpha_2 \in \Aut(\mL;C,\prec)\}}$ which is canonical as a function from $(\mL;C,\prec)$ to
$(\mL;C)$. Clearly, $f'$ also has the property
that there are no $a,b,c,d \in \mL$
such that $f'(a)f'(b)|f'(c)f'(d)$, and the behavior of $f'$ is either $\lin$ or $\nil$ by Lemma~\ref{lem:exhaust}. In both cases, $\{f'\} \cup \Aut(\mL;C)$ generates $\lin$ by Proposition~\ref{prop:nil-lin}. 
\end{proof}

%% file: classification.tex
\subsection{Proof of the main result}
\label{sect:classification}
In this section we finish the proof of Theorem~\ref{thm:endos} and Corollary~\ref{cor:ex-reducts}. We begin by proving two auxiliary results
in Propositions~\ref{prop:rooted} and \ref{prop:unrooted}.



\begin{proposition}\label{prop:rooted}
Let $\Gamma$ be a reduct of $(\mL;C)$. Then one of the following applies.
\begin{enumerate}
\item $\End(\Gamma) = \End(\mL;C)$;
\item $\End(\Gamma)$ contains a constant operation;
\item $\End(\Gamma)$ contains $\lin$;
\item $\End(\Gamma)$ contains $\End(\mL;Q)$. 
\end{enumerate}
\end{proposition}
\begin{proof}
If $\Gamma$ has a non-injective endomorphism, then $\Gamma$ also
has a constant endomorphism by Lemma~\ref{lem:const} and the second item of the statement of the proposition applies.
Therefore we 
suppose in the following that all endomorphisms are injective. 
If all endomorphisms preserve $C$, then
the first item applies and we are done. 
Hence, suppose that 
there is an injective endomorphism $e$
that violates the rooted triple relation,
that is, there are $c_1,c_2,c_3$
such that $c_1|c_2c_3$ and not $e(c_1)|e(c_2)e(c_3)$. 
Under this assumption, we claim that there are $d_1,d_2,d_3 \in \mL$ such that $d_1|d_2d_3$
and $e(d_1)e(d_2)|e(d_3)$. 
By injectivity of $e$, we either have 
$e(c_1)e(c_3)|e(c_2)$ or $e(c_1)e(c_2)|e(c_3)$. 
In the first case, choose $(d_1,d_2,d_3) := (c_1,c_3,c_2)$ and in the second case
choose $(d_1,d_2,d_3) := (c_1,c_2,c_3)$.

By convexity of $\prec$ we have either $d_1 \prec d_2 \prec d_3$, $d_1 \prec d_3 \prec d_2$, $d_2 \prec d_3 \prec d_1$, or $d_3 \prec d_2 \prec d_1$.  
In each case, by the homogeneity of $(\mL;C)$, there exists an $\alpha \in \Aut(\mL;C)$ such that 
$\alpha d_1 \prec \alpha d_2 \prec \alpha d_3$.
After replacing $d_1,d_2,d_3$ by $\alpha d_1,\alpha d_2,\alpha d_3$ and
$e$ by 
$x \mapsto e(\alpha^{-1} x)$, we still have 
$d_1|d_2d_3$ and $e(d_1)e(d_2)|e(d_3)$.
So we assume in the following that
$d_1 \prec d_2 \prec d_3$.
There also exists $\beta \in \Aut(\mL;C)$ 
such that $\beta(e(d_1)) \prec \beta(e(d_2)) \prec \beta(e(d_3))$. 
By replacing $e$ with the function $x \mapsto \beta e(x)$,
we may henceforth assume that $e(d_1) \prec e(d_2) \prec e(d_3)$.

Recall our strategy described at the beginning of this section:  we explained that 
one can 
additionally assume (by Corollary~\ref{cor:const-can}) that $e$ is canonical as a function from
$(\mL;C,\prec,d_1,d_2,d_3)$ to $(\mL;C,\prec)$. 
Define
\begin{align*}
A_1 := & \; \{a :  a \prec d_1, ad_1|d_2d_3\} \\
A_2 := & \; \{a : d_1 \prec a \prec d_2, d_1|ad_2d_3\land a|d_2 d_3\} \\
A_3 := & \; \{a : a \prec d_1, a|d_1d_2d_3 \}
\end{align*}
Note that $A := A_1 \cup A_2 \cup A_3$ is $(d_1,d_2)$-universal: for every finite $X \subseteq \mL$ and arbitrary $x_1,x_2 \in X$  there exists
an $\alpha \in \Aut(\mL;C)$ such that 
\begin{itemize}
\item $\alpha x_1 = d_1$ and $\alpha x_2 = d_2$,
\item $\{\alpha x : x \in X \setminus \{x_1,x_2\}, xx_1|x_2\} \subseteq A_1$, 
\item $\{\alpha x : x \in X \setminus \{x_1,x_2\}, x_1|x_2x\} \subseteq A_2$, and 
\item $\{\alpha x : x \in X \setminus \{x_1,x_2\}, x|x_1x_2\} \subseteq A_3$. 
\end{itemize}
We observe that if $A_i$ is $d_j$-universal, for
$1 \leq i \leq 3$ and $1 \leq j \leq 3$, 
and $e$ is canonical on $A_i$ as a function 
from $(\mL;C,\prec,d_j)$ to $(\mL;C,\prec)$, 
then Lemma \ref{lem:canonical} 
implies that $e$ behaves as $\id_{d_j}$ or $\rer_{d_j}$ on $A_i$ unless $\{e\}\cup \Aut(\mL;C)$ generates $\lin$. 
If $\{e\}\cup \Aut(\mL;C)$ generates $\lin$ then the third item of the statement of the proposition holds and we are done. If $e$ behaves
as $\rer_{d_j}$ on $A_i$ then Corollary \ref{cor:rer} implies
that $\{e\}\cup \Aut(\mL;C)$ generates $\End(\mL;Q)$;
in this case the fourth item of the statement holds. 
Therefore we assume in the following that $e$ behaves as $\id_{d_j}$ on $A_i$. Note that this assumption implies that $e$ preserves $C$ on $\{d_j\}\cup A_i$.
 
Now, pick $r \in A_2$ arbitrarily.  
By the injectivity of $e$, the following cases are exhaustive.
\begin{itemize}
\item $e(d_1)e(d_2)|e(r)e(d_3)$. 
This is in contradiction with the assumption that 
$e$ behaves as $\id_{d_2}$ on $A_2$. To see this, choose an element $a \in A_2$ and note 
that $e(d_2)|e(a)e(d_3)$ by the canonicity of $e$ on $A_2$. This implies that $e(d_2)|e(A_2)$. 
\item $e(d_1)e(d_2)e(r)|e(d_3)$. 
This is in contradiction with the assumption that 
$e$ behaves as $\id_{d_3}$ on $A_2$. To see this, choose an element $a \in A_2$ and note
that $e(d_2)e(a)|e(d_3)$ by the canonicity of $e$ on $A_2$. This implies that $e(d_3)|e(A_2)$. 
\item $e(d_1)e(d_2)e(d_3)|e(r)$.
This is the remaining case that we will consider in the rest of the proof. 
\end{itemize}
Lemma~\ref{lem:two-constants-up} applied to $f := e$, $c_1:=d_1$, $c_2:=d_2$, and $A$ shows that  $e$ preserves $C$ on $\{d_1\} \cup A_1 \cup A_3$, unless $\{e\} \cup \Aut(\mL;C)$ generates $\lin$ or $\End(\mL;Q)$. The same argument can be applied when we exchange  $d_2$ with $d_1$ and $A_2$ with $A_1$ so we assume that $e$ preserves $C$ on $\{d_1\}\cup A_1\cup A_3$ and on $\{d_2\} \cup A_2 \cup A_3$. 
 
If there were a $u \in A_3$ such that $e(d_1)e(u)|e(d_3)$
or $e(d_1)|e(u)e(d_3)$ then $e$ would not behave
as $\id_{d_3}$ or $\id_{d_1}$ on $A_3$ since $e(d_3)|e(A_3)$ or $e(d_1)|e(A_3)$ by the canonicity of $e$, respectively. Hence, we have
$e(d_1)e(d_2)e(d_3)|e(A_3)$. 
If there were a $u \in A_1$ such that $e(d_1)|e(u)e(d_2)$
then by the canonicity of $e$ we would have $e(d_1)|e(A_1)$, and $e$ 
would not behave as $\id_{d_1}$ on $A_1$. Thus $e(u)e(d_1)|e(d_2)$ or $e(u)|e(d_1)e(d_2)$. This implies that either $e(u)e(d_1)|e(d_2)$ for all $u\in A_1$ or $e(u)|e(d_1) e(d_2)$ for all $u\in A_1$.

In the former case,  we have $e(A_1 \cup \{d_1\})|e(d_2)$ and  Lemma~\ref{lem:destruction} applied to $c_1=d_1$ and $c_2 = d_2$ shows that $\{e\} \cup \Aut(\mL;C)$ generates $\nil$, and therefore $\lin$ by Proposition~\ref{prop:nil-lin}. 

In the latter case we show that the conditions in Lemma~\ref{lem:destruction-2} are satisfied for $A$, $g:=e$, $c_1:=d_1$, and $c_2:=d_2$. Clearly, the second and the third conditions are satisfied. It remains to show that the first condition is satisfied. Arbitrarily choose $a_1\in A_1$ and $a_2\in A_2$. If $e(d_1)|e(a_1) e(a_2)$ then for all $u\in A_1$ we have $e(d_1)|e(u) e(a_2)$ by the canonicity of $e$. 
This implies that $e(d_1)|e(A_1)$ which leads to a contradiction since $e$ behaves as $\id_{d_1}$ on $A_1$. Thus either $e(d_1) e(a_1)|e(a_2)$ or $e(d_1)e(a_2)|e(a_1)$ holds.
Hence, Lemma~\ref{lem:destruction-2} shows that 
$\{e\} \cup \Aut(\mL;C)$ generates $\lin$. 
\end{proof}

Proposition~\ref{prop:rooted} leaves us with the task of further analyzing 
the reducts of $(\mL;Q)$.
We first need the following lemma.

\begin{lemma}\label{lem:splitting}
Let $U \subset \mL$ be finite and arbitrarily choose $c \in \mL \setminus U$. Then there are $U_1,\dots,U_k \subseteq U$
such that $U_1 \cup \cdots \cup U_k = U$ and
$(\{c\} \cup \bigcup_{j = 1}^{i-1} U_j)|U_i$ for all
$i \leq k$. 
\end{lemma}
\begin{proof}
By induction on the size of $U$.
If $\{c\}|U$, then $k:=1$ and $U_1 := U$ satisfies
the statement. Otherwise, $|U| \geq 2$, and
by Lemma~\ref{lem:split}, there are two non-empty 
subsets $V,W$ of $U$ such that $V \cup W = U$
and $V|W$. We either have $(\{c\} \cup V)|W$
or $V|(W \cup \{c\})$. In the first case, we inductively have $U_1,\dots,U_{k-1}$ such that $U_1 \cup \cdots \cup U_{k-1} = V$ and
$(\{c\} \cup \bigcup_{j = 1}^{i-1} U_j)|U_i$ for all
$i \leq k-1$. Set $U_k := W$. Then $U_1,\dots,U_{k-1},U_k$ satisfy the requirements from the statement. 
The case when $V|(W \cup \{c\})$ can be shown analogously. 
\end{proof}

\begin{proposition}\label{prop:unrooted}
Let $\Gamma$ be a reduct of $(\mL;Q)$.
Then one of the following cases applies.
\begin{enumerate}
\item All endomorphisms of $\Gamma$ preserve $Q$;
\item $\Gamma$ has a constant endomorphism;
\item $\Gamma$ is preserved by $\lin$.
\end{enumerate}
\end{proposition}

\begin{proof}
If all endomorphisms of $\Gamma$ preserve $Q$, then we are in case one of the statement of the proposition; in the following
we therefore assume that $\Gamma$ has an endomorphism $f$ that violates $Q$. We can then 
choose four elements $d_1,d_2,d_3,d_4 \in \mL$ such that $d_1 d_2:d_3 d_4$ and $f(d_1) f(d_3) :f(d_2) f(d_4)$. By the homogeneity of $(\mL;Q)$ there are $\gamma, \delta \in\Aut(\mL;Q)$ such that 
\begin{itemize}
\item $\gamma(d_1)\prec \gamma(d_2)\prec \gamma(d_3)\prec \gamma(d_4)$, 
\item $\gamma(d_1)\gamma(d_2)|\gamma(d_3)\gamma(d_4)$, 
\item $\delta(f(d_1))\prec \delta(f(d_3))\prec \delta(f(d_2))\prec \delta(f(d_4))$, and
\item $\delta(f(d_1))\delta(f(d_3))|\delta(f(d_2))\delta(f(d_4))$.
\end{itemize}
(Here, the order $\prec$ is still the order as defined in Section~\ref{sect:ordering}.)
 By replacing $f$ by $ \delta \circ f\circ \gamma^{-1}$, we can assume that $d_1\prec d_2\prec d_3\prec d_4$, $d_1d_2|d_3 d_4$, $f(d_1)\prec f(d_3)\prec f(d_2)\prec f(d_4)$, and $f(d_1) f(d_3)|f(d_2)f(d_4)$. 
Corollary~\ref{cor:const-can}
asserts the existence of a function 
$$g \in \overline{\{ \alpha_2 f \alpha_1 : \alpha_1 \in \Aut(\mL;C,\prec,d_1,\dots,d_4),
\alpha_2 \in \Aut(\mL;C,\prec) \} }$$
 which is canonical as a function from 
 $(\mL;C,\prec,d_1,\dots,d_4)$ to $(\mL;C,\prec)$.
Note that there exists an $\alpha \in \Aut(\mL;C,\prec)$ 
such that
$g(d_i) = \alpha(f(d_i))$ for all $i \in \{1,\dots,4\}$,
and, in particular, $g(d_1)g(d_3) | g(d_2) g(d_4)$. 

Let $S=\set{x \in \mL:d_1 d_2 |x\land d_1 d_2 x|d_3 d_4 \land d_2\prec x}$. 
Note that $S$ is $d_1$-universal and $d_2$-universal. 
By Lemma~\ref{lem:canonical}, 
either $g$ behaves on $S$ as $\id_{d_1}$ or $\rer_{d_1}$, or $\{g\} \cup \Aut(\mL;C)$ generates $\lin$. 
Similarly, either $g$ behaves on $S$ as $\id_{d_2}$ or $\rer_{d_2}$, or $\{g\} \cup \Aut(\mL;C)$ generates $\lin$. In the latter cases we are done, so assume
that $g$ behaves on $S$ as $\id_{d_1}$ or $\rer_{d_1}$, and as $\id_{d_2}$ or $\rer_{d_2}$. 

We then show that the conditions of Lemma~\ref{lem:cut} apply to $M := \left<\Aut(\mL;Q) \cup \{g\}\right>$. Let $U \subset \mL$ be finite and arbitrarily choose $u \in U$. 
By Lemma~\ref{lem:splitting} there exists a partition
$U_1 \cup \cdots \cup U_k$ of $U \setminus \{u\}$ such that 
$(\{u\} \cup \bigcup_{j=1}^{i-1} U_j)|U_i$ for all $i \in \{1,\dots,k\}$. 
By $d_1$-universality of $S$, there are subsets $X_1,\dots,X_k$ of $S$ such that 
$X_i | (\bigcup_{j=i+1}^k X_j \cup \{d_1\})$ for all $i \in \{1,\dots,k\}$, and $(X_i;C)$ is isomorphic
to $(U_i;C)$. By the homogeneity of $(\mL;Q)$ 
there is an $\alpha \in \Aut(\mL;Q)$ such that $\alpha(u) = d_3$, $\alpha(U_i) = X_i$, and $\alpha$ preserves $C$ on each $U_i$. 

First consider the case that $g$ behaves 
on $S$ as $\rer_{d_1}$.
We claim that
$g(d_3)|g(S)$ or $g(d_4)|g(S)$. Since $g(d_1)g(d_3)|g(d_2)g(d_4)$ and by the canonicity of $g$ as a function from $(\mL;C,\prec,d_1,\dots,d_4)$ to $(\mL;C,\prec)$, either 
$$g(d_1)g(d_3)g(S)|g(d_2)g(d_4)\text{,}g(d_1)g(d_3)|g(S)g(d_2)g(d_4)\text{, or } g(d_1)g(d_3)g(d_2)g(d_4)|g(S).$$ In the first case $g(d_4)|g(S)$ holds while in the second and third case $g(d_3)|g(S)$ holds. 
We first consider the case
$g(d_3)|g(S)$. Let 
$h := g \circ \alpha$. Clearly, $h$ is in $M$. We will show that $h$ behaves as $\cut_u$ on $U$. Since $g$ behaves as $\rer_{d_1}$ on $S$ and $X_i|(\bigcup_{j=i+1}^k X_j\cup \{d_1\})$ for all $i\in \{1,\ldots,k\}$, it follows from the definition of $\rer_{d_1}$ that $\bigcup_{j=1}^{i} g(X_j)|g(X_{i+1})$ for all $i\in \{1,\ldots,k-1\}$. It implies that $\bigcup_{j=1}^{i}h(U_j)|h(U_{i+1})$ for all $i\in \{1,\ldots,k-1\}$. Since $\alpha$ preserves $C$ on each $U_i$ and $g$ preserves $C$ on each $X_i$, $h$ preserves $C$ on each $U_i$. Since $\bigcup_{j=1}^{i} U_j|U_{i+1}$, $\bigcup_{j=1}^{i}h(U_j)|h(U_{i+1})$ for all $i\in \{1,\ldots,k-1\}$ and $h$ preserves $C$ on each $U_i$, it follows that $h$ preserves $C$ on $\bigcup_{i=1}^k U_i$. Since $g(d_3)|g(S)$, we have that $h(u)|h(\bigcup_{i=1}^k U_i)$. Thus $h$ behaves as $\cut_u$ on $U$. By Lemma~\ref{lem:cut}, $\Gamma$ is preserved by $\lin$. 
The case $g(d_4)|g(S)$ can be treated similarly
(by choosing $\alpha \in \Aut(\mL;Q)$ such that $\alpha(u) = d_4$ instead of $\alpha(u) = d_3$).

Finally, we consider the case when $g$ behaves on $S$ as $\id_{d_1}$. 
By the canonicity of $g$ as a function
from $(\mL;C,\prec,d_1,\dots,d_4)$ to $(\mL;C,\prec)$
and since all the elements of $S$ lie in the same orbit
of $(\mL;C,\prec,d_1,\dots,d_4)$, either $g(d_1)|g(d_2)g(S)$,  $g(S)g(d_1)|g(d_2)$,  or $g(d_1)g(d_2)|g(x)$ for all $x \in S$.
The first case is impossible because 
$g$ behaves as $\id_{d_1}$ on $S$. 
The second case is impossible, too: to see this,
pick $a,b,c \in S$ such that $d_1a|b$ and $d_1ab|c$. Since
$g$ behaves on $S$ as $\id_{d_1}$, we 
have 
$g(d_1)g(a)g(b)|g(c)$ and $g(d_1)g(a)|g(b)$,
and $g(d_1)g(a)g(b)g(c) | g(d_2)$ by assumption.  
In case that $g$ behaves as $\id_{d_2}$ 
we would have $g(d_2)g(a)|g(b)$ which is inconsistent with the above.
In case that $g$ behaves as $\rer_{d_2}$ we would have $g(d_2) g(a) : g(b) g(c)$ which is inconsistent with the above, too. 

In the third and last case, we first show that $g$ does not behave as $\rer_{d_2}$ on $S$. Let $a,b,c\in S$ such that $d_1ab|c\wedge d_1 a|b$. Since $g$ behaves as $id_{d_1}$, we have $g(d_1)g(a)g(b)|g(c)$. By the assumption $g(d_1)g(d_2)|g(x)$ for all $x\in S$, we have $g(d_1)g(d_2)|g(c)$. It follows from $g(d_1)g(a)g(b)|g(c)$ that $g(d_2)g(a)g(b)|g(c)$. Therefore $g$ does not behave as $\rer_{d_2}$ on $S$. Hence $g$ behaves as $\id_{d_2}$ on $S$. Thus, $\big(\{g(d_1),g(d_2)\} \cup \bigcup_{j=i+1}^k g(X_j)\big) | g(X_i)$ for all $i \in \{1,\dots,k\}$. 
Since $g(d_1)g(d_3) | g(d_2)g(d_4)$,
we therefore must have that $\big(\{g(d_3)\} \cup \bigcup_{j=i+1}^k g(X_j)\big) | g(X_i)$ for all $i \in \{1,\dots,k\}$. 
Since $(\mL;C)$ embeds all finite leaf structures,
there are subsets $Z_1,\dots,Z_k$ of $\mL$ and $z \in \mL$ such that $z | \bigcup_{j=1}^{k} Z_j$, 
$(\bigcup_{j=1}^{i-1} Z_j) | Z_i$, and $(Z_i;C)$ is isomorphic to $(g(X_i);C)$ for all $i \in \{1,\dots,k\}$. The situation is illustrated in Figure~\ref{fig:unrooted}. 
\begin{figure}
\begin{center}
\includegraphics[scale=.7]{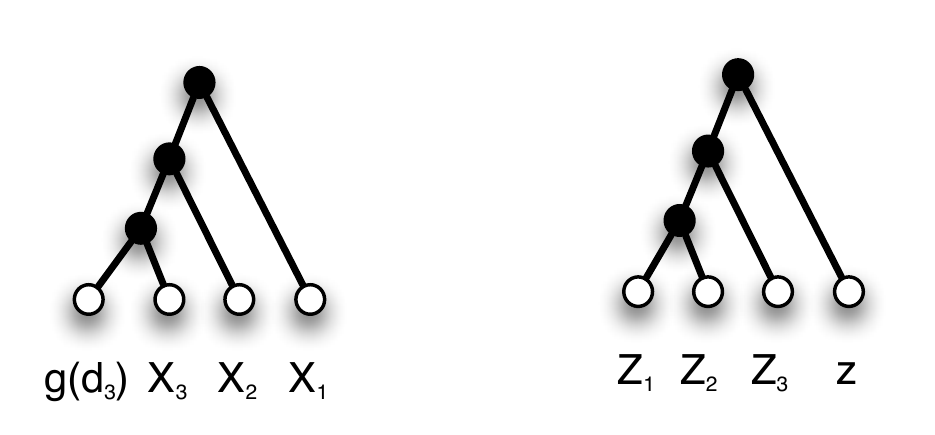} 
\end{center}
\caption{Illustration for the last proof step for Proposition~\ref{prop:unrooted}.}
\label{fig:unrooted}
\end{figure}
By the homogeneity of $(\mL;Q)$
there is a $\beta \in \Aut(\mL;Q)$ such that $\beta(g(d_3)) = z$, $\beta$ preserves $C$ on each $g(X_i)$,
and $\beta(g(X_i)) = Z_i$ for all $i \in \{1,\dots,k\}$. 
Then $\beta \circ g \circ \alpha$ is in $M$ 
and behaves as $\cut_u$ on $U$. Again, 
Lemma~\ref{lem:cut} implies that $\Gamma$ is preserved by $\lin$. 
\end{proof}

We can now prove Theorem~\ref{thm:endos}
and Corollary~\ref{cor:ex-reducts} that have already
been stated in Section~\ref{sect:results}. 

\begin{proof}[Proof of Theorem~\ref{thm:endos}]
Let $\Gamma$ be a reduct of $(\mL;C)$. 
We apply Proposition~\ref{prop:rooted} 
and consider the following cases. 
\begin{itemize}
\item
All endomorphisms of $\Gamma$ preserve $C$. Then $\End(\Gamma) \subseteq \End(\mL;C)$; we claim that the opposite inclusion holds as well. Since $\End(\Gamma)$ is closed, it suffices to show that for every $e \in \End(\mL;C)$ and every finite $S \subset \mL$ there exists an $f \in \End(\Gamma)$ 
such that $f(s) = e(s)$ for all $s \in S$. Since $e$ preserves $C$, $e|_{S}$ is a partial isomorphism from $(S;C)$ to $(e(S);C)$ by Lemma~\ref{lem:pres-c}. 
By homogeneity, $e|_{S}$ can be extended to an automorphism $f\in \Aut(\mL;C)$. Since $\Aut(\mL;C)\subseteq \End(\Gamma)$, we have $f\in \End(\Gamma)$.

\item $\Gamma$ has a constant endomorphism. Then 
there is nothing to show since the second item of the statement applies. 
\item $\Gamma$ is preserved by
$\lin$. Then Lemma~\ref{lem:lin} shows that
$\Gamma$ is homomorphically equivalent to 
a reduct of 
$(\mL;=)$ and the third item of the statement applies. 
\item $\End(\Gamma)$ contains $\End(\mL;Q)$. This case implies that $\Gamma$ is a reduct of $(\mL;Q)$. If $\Gamma$ has a constant endomorphism or $\Gamma$ is preserved by $\lin$, then we are done as the above cases. Otherwise, by Proposition \ref{prop:unrooted} all endomorphisms of $\Gamma$ preserve $Q$. This means that $\End(\Gamma)\subseteq \End(Q)$. Hence we have $\End(\Gamma)=\End(Q)$, therefore the fourth item applies.
\end{itemize}
By Proposition~\ref{prop:rooted}, these four cases are exhaustive.
\end{proof}

One may observe at this point that the proof of Theorem~\ref{thm:endos} does not rely on any results
concerning Jordan permutation groups.
We finally show that every reduct of $\Gamma$
is existentially interdefinable with $(\mL;C)$, with
$(\mL;Q)$, or with $(\mL;=)$.

\begin{proof}[Proof of Corollary~\ref{cor:ex-reducts}]
Let $\Gamma$ be a reduct of $(\mL;C)$. 
Let $\Gamma'$ be the expansion of $\Gamma$ by the relations defined by negations of atomic formulas over $\Gamma$, including the equality relation (for example, when $R$ is a ternary relation of $\Gamma$, the structure $\Gamma'$
contains the binary relation defined by $\neg R(x,x,y)$). 
We apply Theorem~\ref{thm:endos}
to $\Gamma'$. 
Since for every atomic formula $\phi$ over $\Gamma'$ the signature of $\Gamma'$ also
contains a relation symbol for $\neg \phi$,
all endomorphisms of $\Gamma'$ must be
embeddings, and therefore item 2 of Theorem~\ref{thm:endos}
is impossible. 
If $\Gamma'$ has the same endomorphisms
as $(\mL;C)$ or $(\mL;Q)$, then by Proposition~\ref{prop:galois} the structure
$\Gamma'$ is existentially positively interdefinable with $(\mL;C)$ 
or with $(\mL;Q)$; hence, $\Gamma$ is existentially interdefinable with one of those structures and we are done. 
Otherwise, $\Gamma'$ is homomorphically
equivalent with a reduct $\Delta$ of $(\mL;=)$. 
Again, the homomorphism from $\Gamma'$ to
$\Delta$ must in fact be an embedding.
Hence, $\Gamma'$ is isomorphic to 
a substructure of $\Delta$. Since $\Delta$
is preserved by all permutations, so is this substructure, and so is $\Gamma'$. 
It follows that $\Gamma$ is preserved
by all permutations, so $\Gamma$ is a reduct of $(\mL;=)$ by Proposition~\ref{prop:equality}. 
In fact, $\Gamma$ is even preserved by all injective maps from
$\mL$ to $\mL$ and therefore by all 
self-embeddings of $(\mL;=)$. Hence,
Proposition~\ref{prop:galois} shows that 
$\Gamma$ has an existential definition over $(\mL;=)$. 
Conversely, $(\mL;=)$
has an existential definition in every structure with domain $\mL$,
so $\Gamma$ is existentially interdefinable with $(\mL;=)$. 
\end{proof}